\documentclass[12pt]{article}  
\def\sq{\hbox {\rlap{$\sqcap$}$\sqcup$}}
\def\sq{\hbox {\rlap{$\sqcap$}$\sqcup$}}
\def\R{ {\rm R \kern -.31cm I \kern .15cm}}
\def\C{ {\rm C \kern -.15cm \vrule width.5pt \kern .12cm}}
\def\Z{ {\rm Z \kern -.27cm \angle \kern .02cm}}
\def\N{ {\rm N \kern -.26cm \vrule width.4pt \kern .10cm}}
\def\1{{\rm 1\mskip-4.5mu l} }
\def\lsim{\raise0.3ex\hbox{$<$\kern-0.75em\raise-1.1ex\hbox{$\sim$}}}
\def\gsim{\raise0.3ex\hbox{$>$\kern-0.75em\raise-1.1ex\hbox{$\sim$}}}
\def\noi{\noindent}

\def\beq{\begin{equation}}   \def\eeq{\end{equation}}
\def\bea{\begin{eqnarray}}  \def\eea{\end{eqnarray}}
\def\nn{\nonumber}
\def\noi{\noindent}
\def\beeq{\begin{eqnarray}} \def\eeeq{\end{eqnarray}}
\newcommand\mysection{\setcounter{equation}{0}\section}

\newcounter{hran}

\begin{document} 
\centerline{\Large\bf Long Range Scattering and Modified } 
 \vskip 3 truemm \centerline{\Large\bf  Wave Operators for the Maxwell-Schr\"odinger
System} \vskip 3 truemm 
\centerline{\Large\bf II. The general case}
\vskip 0.5 truecm

\centerline{\bf J. Ginibre}
\centerline{Laboratoire de Physique Th\'eorique\footnote{Unit\'e Mixte de
Recherche (CNRS) UMR 8627}}  \centerline{Universit\'e de Paris XI, B\^atiment
210, F-91405 ORSAY Cedex, France}
\vskip 3 truemm
\centerline{\bf G. Velo}
\centerline{Dipartimento di Fisica, Universit\`a di Bologna}  \centerline{and INFN, Sezione di
Bologna, Italy}

\vskip 1 truecm

\begin{abstract}
We study the theory of scattering for the Maxwell-Schr\"odinger system in space dimension 3,
in the Coulomb gauge. We prove
the existence of modified wave operators for that system with no size restriction on the
Schr\"odinger and Maxwell asymptotic data and we determine the asymptotic behaviour in time of solutions in the
range of the wave operators. The method consists in partially solving the Maxwell equations
for the potentials, substituting the result into the Schr\"odinger equation, which then
becomes both nonlinear and nonlocal in time. The Schr\"odinger function is then parametrized in terms of an amplitude and a phase satisfying a suitable auxiliary system, and the Cauchy problem for that system, with prescribed asymptotic behaviour determined by the asymptotic data, is solved by an energy method, thereby leading to solutions of the original system with prescribed asymptotic behaviour in time. This paper is the generalization of a previous paper with the same title. However it is entirely self contained and can be read without any previous knowledge of the latter.
\end{abstract}

\vskip 1.5 truecm
\noi MS Classification :    Primary 35P25. Secondary 35B40, 35Q40, 81U99.\par \vskip 2 truemm

\noi Key words : Long range scattering, modified wave operators, Maxwell-Schr\"odinger system. \par 
\vskip 1 truecm

\noindent LPT Orsay 06-37\par
\noindent June 2006\par \vskip 3 truemm

\newpage
\pagestyle{plain}
\baselineskip 18pt

\mysection{Introduction}
\hspace*{\parindent} This paper is devoted to the theory of scattering and more precisely to the existence of
modified wave operators for the Maxwell-Schr\"odinger system (MS) in $3 + 1$ dimensional
space time. This system describes the evolution of a charged nonrelativistic quantum
mechanical particle interacting with the (classical) electromagnetic field it generates. It can be written as follows~:
\beq
\label{1.1e}
\left \{ \begin{array}{l} i\partial_t u = -(1/2) \Delta_A u + 
A_{e} u \\ \\ \sq A_{e} - \partial_t \left ( \partial_t A_{e} + \nabla \cdot A\right )  = J_0 \\ \\ \sq A + \nabla \left ( \partial_t A_{e} + \nabla \cdot A\right )  = J\end{array} \right . \eeq

\noi where 
\beq
\label{1.2e}
J_0 = |u|^2 \qquad , \qquad J = \  {\rm Im}\ \bar{u} \ \nabla_Au\ .
\eeq

\noi Here $(A, A_e)$ is an ${I\hskip-1truemm R}^{3+1}$ valued function defined in space time ${I\hskip-1truemm R}^{3+1}$, $\nabla_A = \nabla - iA$ and $\Delta_A = \nabla_A^2$ are the covariant gradient and covariant Laplacian respectively and $\sq = \partial_t^2 - \Delta$ is the d'Alembertian. An important property of that system is its gauge invariance, namely the invariance under the transformation
\beq
\label{1.3e}
\left ( u, A, A_e \right ) \rightarrow \left (u \exp (- i \theta), A - \nabla \theta, A_e + \partial_t \theta \right ) \ ,
\eeq

\noi where $\theta$ is an arbitrary real function defined in ${I\hskip-1truemm R}^{3+1}$. As a consequence of that invariance, the system (\ref{1.1e}) is underdetermined as an evolution system and has to be supplemented by an additional equation, called a gauge condition. Here we shall use exclusively the Coulomb gauge condition, namely $\nabla \cdot A = 0$, which experience shows to be the most convenient one for the purpose of analysis.\par

The MS system (\ref{1.1e}) is known to be locally well posed both in the Coulomb gauge and in the Lorentz gauge $\partial_t A_e + \nabla \cdot A = 0$ in sufficiently regular spaces \cite{19r} \cite{20r} and to have weak global solutions in the energy space \cite{13r}. \par

A large amount of work has been devoted to the theory of scattering for nonlinear equations and systems centering on the Schr\"odinger equation, in particular for nonlinear Schr\"odinger (NLS) equations, Hartree equations, Klein-Gordon-Schr\"odinger (KGS) systems, Wave-Schr\"odinger (WS) systems and Maxwell-Schr\"odinger (MS) systems. As in the case of the linear Schr\"odinger equation, one must distinguish the short range case from the long range case. In the former case, ordinary wave operators are expected and in a number of cases proved to exist, describing solutions where the Schr\"odinger function behaves asymptotically like a solution of the free Schr\"odinger equation. In the latter case, ordinary wave operators do not exist and have to be replaced by modified wave operators including a suitable phase in their definition. In that respect, the MS system (\ref{1.1e}) in ${I\hskip-1truemm R}^{3+1}$ belongs to the borderline (Coulomb) long range case, because of the $t^{-1}$ decay in $L^{\infty}$ norm of solutions of the wave equation. Such is the case also for the Hartree equation with $|x|^{-1}$ potential and for the WS system in ${I\hskip-1truemm R}^{3+1}$.\par

The construction of the wave operators and more precisely the local Cauchy problem at infinity in the long range cases of the previous nonlinear equations and systems has been treated essentially by two methods. The first method is rather direct, starting from the original equation or system. It is intrinsically restricted to the case of small Schr\"odinger data and to the borderline long range case. It was initiated in \cite{23r} in the case of the NLS equation in ${I\hskip-1truemm R}^{1+1}$ and was subsequently extended to the NLS equation in ${I\hskip-1truemm R}^{n+1}$ for $n = 2, 3$, to the Hartree equation in ${I\hskip-1truemm R}^{n+1}$ for $n \geq 2$, to the KGS system in ${I\hskip-1truemm R}^{2+1}$ and to the WS and MS systems in ${I\hskip-1truemm R}^{3+1}$. In particular the WS and MS systems in ${I\hskip-1truemm R}^{3+1}$ are treated in \cite{25r} \cite{10r} and in \cite{28r} \cite{26r} \cite{11r} respectively by that method. We refer to \cite{12r} for a review of that method as applied to the previous systems and for additional references. The second method is more complex and starts with a phase/amplitude separation of the Schr\"odinger function, inspired by previous work on the Hartree equation \cite{14r} \cite{15r}. The main interest of that method is to eliminate the smallness condition on the Schr\"odinger function and the restriction to the borderline long range case in the case of the Hartree equation. It has been applied to the Hartree equation in ${I\hskip-1truemm R}^{n+1}$ \cite{4r}, with subsequent improvements eliminating a loss of regularity between the asymptotic data and the solution \cite{21r} \cite{22r}, to the WS system in ${I\hskip-1truemm R}^{3+1}$ \cite{5r} \cite{6r} \cite{8r} and to the MS system in ${I\hskip-1truemm R}^{3+1}$ in the special case of vanishing asymptotic data for the Maxwell field \cite{7r}.\par 

The present paper is devoted to the extension of the results of \cite{7r} (hereafter referred to as I) to the general case of non vanishing asymptotic data for the Maxwell field, and actually to the case of arbitrarily large asymptotic data both for the Schr\"odinger and for the Maxwell field. Thus this paper is a generalization of I and its general structure closely follows that of I. However it is entirely self contained and can be read without any previous knowledge of I, from which it differs significantly in technical detail. Those technical differences being put aside for the moment, the method used in this paper follows rather closely that of I. One first replaces the Maxwell equation for the vector potential by the associated integral equation and substitutes the
latter into the Schr\"odinger equation, thereby obtaining a new Schr\"odinger equation which is both nonlinear and nonlocal in time. The latter is then treated as in I, namely $u$ is expressed in terms of an amplitude $w$ and a phase $\varphi$ satisfying an auxiliary system similar to that introduced in I. In contrast with I, however, the representation of $u$ in terms of $(w, \varphi )$ is now chosen to involve a change of $t$ into $1/t$ in order to simplify the treatment of time derivatives (see Remark 2.1 below for details). One then obtains an auxiliary system of equations for $(w, \varphi )$, now to be considered in a neighborhood of $t=0$ instead of $t = \infty$. One solves the Cauchy problem at $t=0$ for that system in a roundabout way, since that system is singular at $t=0$. Returning to the original variable $u$ then allows one to construct modified wave operators for the original system (\ref{1.1e}). The detailed construction is too complicated to allow for a more precise description at this stage, and will be described in heuristic terms in
Section 2 below. \par

The main technical difference between this paper and I is that here we systematically use time derivatives (and to some extent covariant space derivatives), both in the definition of the basic function spaces and in the derivation of the relevant estimates, whereas in I we used only ordinary space derivatives. Since one time derivative is homogeneous to two space derivatives for the Schr\"odinger equation, one can thereby reach a given level of regularity in the space variable for the Schr\"odinger function by using lower total order derivatives than would be needed if one used space derivatives only. That property has been used extensively in \cite{28r} \cite{26r} \cite{8r} \cite{9r} \cite{10r} \cite{11r}. In the present case, it turns out that one needs the amplitude $w$ for the Schr\"odinger function mentioned above to have $H^k$ regularity with $k > 5/2$ in order to handle the nonlinear character of the MS system. Here, restricting our attention to integer order derivatives, we reach that level by using only one time and one space derivative, which yields a regularity corresponding to $H^k$ with $k = 3$. In addition, the use of time derivatives naturally leads to a systematic use of covariant instead of ordinary space derivatives of the Schr\"odinger function, thereby allowing to exploit the algebraic structure of the MS system more efficiently.\par

When dealing with systems containing a Schr\"odinger equation and a wave equation coupled together, such as the WS and MS systems, one always encounters the difficulty that the two constituent equations have different propagation properties, so that the product of solutions of the two equations decays rather weakly in the light cone directions. That difficulty can be circumvented by assuming a support condition for the Schr\"odinger asymptotic data (more precisely for the Fourier transform thereof), which implies that the corresponding solution of the free Schr\"odinger equation decays sufficiently fast in the light cone directions. Such an assumption appears in early works on the subject \cite{28r} \cite{5r}. In the framework of the first method mentioned above, that assumption can be eliminated by using a suitably improved asymptotic form for the Schr\"odinger function \cite{25r} \cite{26r}  (see \cite{12r} for a review). In the framework of the more complex method that we are using here, that can also be done for the WS system \cite{6r} \cite{8r}. In the case of the MS system however, the situation turns out to be less favourable. Even in the simpler case of a linear Schr\"odinger equation minimally coupled to a free external electromagnetic field, the support assumption has so far been eliminated only at the expense of an intermediate construction which does not allow for an explicit characterization of the admissible asymptotic Schr\"odinger data \cite{9r}. That construction could certainly be adapted to the MS system. However since it is complicated and since the result is not entirely satisfactory, we shall refrain from pursuing that line any further and we shall for simplicity stick to the previously used support condition in the present paper, awaiting a satisfactory treatment of that point in the linear case.\par

We now give a brief outline of the contents of this paper. A more detailed description of the technical parts will be given at the end of Section 2. After collecting some notation and preliminary estimates in Section 3, we start studying the dynamics for the auxiliary system in Section 4 and we solve the local Cauchy problem at $t=0$ for that system in Sections 5 and 6, which contain the main technical results of this paper.  We finally come back from the auxiliary system to the original one and construct the modified wave operators for the latter in Section 7, where the final result is stated in Proposition 7.2. A simplified version of that result is stated below as Proposition 2.1.

\mysection{Heuristics and formal computations}
\hspace*{\parindent} In this section we discuss in heuristic terms the construction of the modified wave operators for the MS system as it will be performed in this paper and we derive the equations needed for that purpose. We first recast the MS system in the Coulomb gauge in a slightly different form. Under the Coulomb gauge condition $\nabla \cdot A = 0$, the MS system becomes
\bea \label{2.1e} 
&&i\partial_t u = - (1/2) \Delta_A u + A_e u\\
\label{2.2e}
&&\Delta A_e = - J_0 \\
&&\sq A + \nabla \left (\partial_t A_e \right ) = J \ .
\label{2.3e}
\eea

\noi We replace that system by a formally equivalent one in the following standard way. We solve (\ref{2.2e}) for $A_e$ as
\beq
\label{2.4e}
A_e = - \Delta^{-1} \ J_0 = (4 \pi |x|)^{-1} \ * \  |u|^2 \equiv g(u)
\eeq

\noi where $*$ denotes the convolution in ${I\hskip-1truemm R}^3$, so that by the current
conservation $\partial_t J_0 + \nabla \cdot J = 0$,  
\beq
\label{2.5e}
\partial_t \ A_e = \Delta^{-1} \nabla \cdot J \ .
\eeq

\noi Substituting (\ref{2.4e}) into (\ref{2.1e}) and (\ref{2.5e}) into (\ref{2.3e}), we obtain
the new system
 \bea
\label{2.6e}
&&i \partial_t u = - (1/2) \Delta_A u + g(u) u \\
\label{2.7e}
&&\sq A = PJ \equiv P\ {\rm Im}\ \bar{u}\nabla_A u
\eea

\noi where $P = \1 - \nabla \Delta^{-1} \nabla$ is the projector on divergence free vector fields. The system (\ref{2.6e}) (\ref{2.7e}) is the starting point of our investigation. We want to address the problem of classifying the asymptotic behaviours in time of the solutions of the system (\ref{2.6e}) (\ref{2.7e}) by relating them to a set of model functions ${\cal V} = \{ v = v(v_+)\}$ parametrized by some data $v_+$ and with suitably chosen and preferably simple asymptotic behaviour in time. For each $v \in {\cal V}$, one tries to construct a solution $(u, A)$ of the system (\ref{2.6e}) (\ref{2.7e}) defined at least in a neighborhood of infinity in time and such that $(u, A)(t)$ behaves as $v(t)$ when $t \to \infty$ in a suitable sense. We then define the wave operator as the map $\Omega : v_+ \to (u, A)$ thereby obtained. A similar question can be asked for $t \to - \infty$. We restrict our attention to positive time. The more standard definition of the wave operator is to define it as the map $v_+ \to (u,A)(0)$, but what really matters is the solution $(u, A)$ in the neighborhood of infinity in time, namely in some interval $[T, \infty )$. Continuing such a solution down to $t = 0$ is a somewhat different question which we shall not touch here. \par 

In the case of the MS system, which is long range, it is known that one cannot take for ${\cal V}$ the set of solutions of the linear problem underlying (\ref{2.6e}) (\ref{2.7e}), namely of the free Schr\"odinger equation for $u$ and of the free wave equation for $A$, and one of the tasks that will be performed in this paper will be to construct a better set ${\cal V}$ of model asymptotic functions. The same situation prevails for long range Hartree equations and for the WS system in ${I\hskip-1truemm R}^{3+1}$.\par

Constructing the wave operators essentially amounts to solving the Cauchy problem with infinite initial time. The system (\ref{2.6e}) (\ref{2.7e}) in this form is not well suited for that purpose, and we now perform a number of transformations leading to an auxiliary system for which that problem can be handled. We first replace the equation (\ref{2.7e}) by the associated integral equation namely 
\beq
\label{2.8e}
A = A_0 + A'(u, A)
\eeq

\noi where
\beq
\label{2.9e}
A_0 = (\cos \omega t)A_+ + \omega^{-1} (\sin \omega t) \dot{A}_+\ ,
\eeq
\beq
\label{2.10e}
A'(u,A) = - \int_t^{\infty} dt' \omega^{-1} \sin (\omega (t-t')) P J (u, A) (t')
\eeq

\noi with $\omega = (- \Delta )^{1/2}$.\par

In particular, $A_0$ is a solution of the free (vector valued) wave equation with initial data $(A_+ ,\dot{A}_+)$ at $t = 0$, and $(A_+ ,\dot{A}_+)$ is naturally interpreted as the asymptotic state for $A$. In order to ensure the condition $\nabla\cdot A = 0$, we assume that $\nabla \cdot A_+ = \nabla \cdot \dot{A}_+ = 0$. \par

We next perform a change of variables which is well adapted to the study of the asymptotic behaviour in time of solutions of the Schr\"odinger equation. The unitary group 
\beq 
\label{2.11e}
U(t) = \exp (i(t/2)\Delta )
\eeq
\noi which solves the free Schr\"odinger equation can be written as
\beq 
\label{2.12e}
U(t) = M(t) \ D(t) \ F \ M(t)
\eeq
\noi where $M(t)$ is the operator of multiplication by the function
\beq 
\label{2.13e}
M(t) = \exp  \left ( i x^2/2t \right ) \ ,
\eeq
\noi $F$ is the Fourier transform and $D(t)$ is the dilation operator
\beq 
\label{2.14e}
(D(t)\ f)(x) = (it)^{-3/2} \ f(x/t)
\eeq
\noi normalized to be unitary in $L^2$. We shall also need the operator $D_0(t)$ defined by
\beq 
\label{2.15e}
\left ( D_0(t) f\right )(x) = f(x/t) \ .
\eeq
\noi We parametrize the Schr\"odinger function $u$ in terms of an amplitude $w$ and of a real phase $\varphi$ as 
\beq 
\label{2.16e}
u(t) = M(t) \ D(t) \exp ( i \varphi (1/t)) \overline{w}(1/t) \ .
\eeq

\noi Correspondingly we change the variable for the vector potential from $A$ to $B$ according to 
\beq
\label{2.17e} 
A(t) =-  t^{-1} \ D_0(t) \ B(1/t)
\eeq

\noi and similarly for $A_0$ and $A'$. \\

\noi {\bf Remark 2.1.} The change of variables (\ref{2.16e}) (\ref{2.17e}) is slightly different from that made in I. Denoting by $w_*$, $\varphi_*$, $B_*$ the variables introduced in I, the correspondence is
$$\left ( w_* \exp (- i \varphi_*)\right ) (t) = \overline{w \exp (-i \varphi}) (1/t)$$
$$B_*(t) = - B(1/t)\ .$$

\noi The main practical consequence is that we shall have to study the system of equations for $(w, \varphi , B)$ (see (\ref{2.41e}) below) in the neighborhood of $t=0$ whereas in I the system for $(w_*, \varphi_*, B_*)$ was studied in the neighborhood of $t = \infty$. The present choice is more convenient when dealing with time derivatives, as we shall do systematically in this paper. Note also that the change of variables (\ref{2.16e}) reduces to the pseudoconformal inversion for $\varphi = 0$.\\

We now perform the change of variables (\ref{2.16e}) (\ref{2.17e}) on the system (\ref{2.6e}) (\ref{2.7e}). Substituting (\ref{2.16e}) (\ref{2.17e}) into (\ref{2.6e}) and commuting the Schr\"odinger operator with MD, we obtain
\bea
\label{2.18e}
&&\left \{ \left ( i \partial_t + (1/2) \Delta_A - g(u)\right ) u \right \} (t) \nn \\
&&= t^{-2} M(t) D(t) \left \{ \overline{\left ( i \partial_t + (1/2) \Delta_B - {\check B} - t^{-1} g(w) \right ) \exp (- i \varphi ) w}\right \} (1/t)
\eea

\noi where
\beq
\label{2.19e}
{\check B}(t) = t^{-1} x \cdot B(t)
\eeq

\noi and more generally for any ${I\hskip-1truemm R}^{3}$ valued function $f$ of space time 
\beq
\label{2.20e}
{\check f}(t) = t^{-1} x \cdot f(t) \ .
\eeq

\noi Therefore (\ref{2.6e}) becomes
\beq
\label{2.21e}
\left ( i \partial_t + \partial_t \varphi + (1/2) \Delta_K - {\check B}(t) - t^{-1} g(w) \right ) w = 0
\eeq

\noi where
\bea
\label{2.22e}
&&K = B + s \qquad , \qquad s = \nabla \varphi \ , \\
&&\Delta_K = \nabla_K^2 = (\nabla - i K)^2 \ .
\label{2.23e}
\eea

We next turn to the Maxwell equation in the integral form (\ref{2.8e})-(\ref{2.10e}). Substituting (\ref{2.16e}) (\ref{2.17e}) into the definition of $J$, we obtain
\beq
\label{2.24e}
J(t) = t^{-3} D_0 (t) \ N(t)
\eeq

\noi where
\bea
\label{2.25e}
&&N(t) = - M_1(1/t) - t^{-1} M_2(1/t)\ , \\
\label{2.26e}
&&M_1 = - x |w|^2 \ , \\
&&M_2 = {\rm Im}\ \overline{w} \nabla_K w \ .
\label{2.27e}
\eea

\noi Substituting (\ref{2.24e}) into (\ref{2.10e}) and letting $t' = t \nu$ yields
\bea
\label{2.28e}
&&A'(t) = - \int_t^{\infty} dt' \omega^{-1} \sin (\omega (t - t')) t{'}^{-3} P D_0 (t') N(t')\nn\\
&&= t^{-1} D_0(t) \int_1^{\infty} d\nu\ \nu^{-3} \sin (\omega (\nu - 1)) P D_0 (\nu )N(t\nu )\nn \\
&&= - t^{-1} D_0 (t) B' (1/t)
\eea

\noi in accordance with (\ref{2.17e}), where
\beq
\label{2.29e}
B' = B_1 + B_2 \qquad , \qquad B_j = t^{j-1} F_j (PM_j) \quad \hbox{for $j = 1,2$}\ ,
\eeq

\noi with $M_j$ defined by (\ref{2.26e}) (\ref{2.27e}) and $F_j$ defined by
\beq
\label{2.30e}
F_j (M) = \int_1^{\infty} d\nu\ \nu^{-2-j} \omega^{-1} \sin (\omega (\nu - 1)) D_0(\nu )\ M(t /\nu )\ .
\eeq

\noi Note that (\ref{2.26e}) (\ref{2.29e}) (\ref{2.30e}) yield $B_1$ as an explicit function of $w$, namely
\beq
\label{2.31e}
B_1 \equiv B_1 (w) \equiv B_1 (w, w)
\eeq

\noi where 
\beq
\label{2.32e}
B_1(w_1, w_2) = - F_1 (Px \ {\rm Re} \ \overline{w}_1 w_2)\ .
\eeq

\noi On the other hand (\ref{2.27e}) (\ref{2.29e}) (\ref{2.30e}) yield only an implicit equation for $B_2$, since $M_2$ still contains $B_2$ through $K$. We shall take the following point of view. We shall take $B_2$ as the dynamical variable for the Maxwell field, we shall regard the decomposition
\beq
\label{2.33e}
B = B_0 + B_1 + B_2
\eeq

\noi as a change of dynamical variable from $B$ to $B_2$, and the relations (\ref{2.29e}) (\ref{2.30e}) with $j = 2$ as the equation for $B_2$, in the form
\beq
\label{2.34e}
B_2 = {\cal B}_2 (w, w, s + B)
\eeq

\noi where we define
\beq
\label{2.35e}
{\cal B}_2 (w_1, w_2, K) \equiv t F_2\left ( P\ {\rm Im} \ \overline{w}_1 \nabla_K w_2\right )\ .
\eeq

We now come back to the transformed Schr\"odinger equation (\ref{2.21e}). Here we have parametrized $u$ in terms of an amplitude $w$ and a phase $\varphi$ and we have only one equation for two functions $(w, \varphi )$. We then arbitrarily impose a second equation, namely an equation for the phase $\varphi$, thereby splitting (\ref{2.21e}) into a system of two equations, the other one of which being an equation for $w$. There is a large amount of freedom in the choice of the equation for the phase. The role of the phase is to cancel the long range terms in (\ref{2.21e}) coming from the interaction. The terms coming from the covariant Laplacian are expected (and will turn out) to be short range. Such is also the case for the contribution of $B_2$ to ${\check B}$ because of the factor $t$ in (\ref{2.35e}). The term $t^{-1} g(w)$ is clearly long range (of Hartree type) and is therefore included in the $\varphi$ equation. The term ${\check B}_1$ in ${\check B}$ is also long range, but since it is less regular than the previous one, it is convenient to split it into a short range and a long range part. Let $\chi \in {\cal C}^{\infty} ({I\hskip-1truemm R}^3 , {I\hskip-1truemm R})$, $0 \leq \chi \leq 1$, $\chi (\xi ) = 1$ for $|\xi | \leq 1$, $\chi (\xi ) = 0$ for $|\xi | \geq 2$, and let $0 < \beta < 1$. We define
\beq
\label{2.36e}
\left \{ \begin{array}{l}\chi_L \equiv F^* \chi (\cdot t^{\beta})F \qquad , \quad \chi_S = \1 - \chi_L \\ \\ {\check B}_1 = {\check B}_{1S} + {\check B}_{1L}\qquad , \quad  {\check B}_{1S}  = \chi_S  {\check B}_{1} \qquad , \quad  {\check B}_{1L}   = \chi_L  {\check B}_{1} \ . \end{array} \right . 
\eeq

\noi Corresponding to the fact that ${\check B}_0$ and ${\check B}_2$ are regarded as short range, we denote
\beq
\label{2.37e}
{\check B}_S = {\check B}_0 + {\check B}_{1S} + {\check B}_2\qquad , \quad {\check B}_L = {\check B}_{1L} \ .
\eeq

\noi As the equation for $\varphi$, we impose
\beq
\label{2.38e}
\partial_t \varphi = t^{-1} g(w) + {\check B}_{1L}
\eeq

\noi so that the remaining equation for $w$ becomes
\beq
\label{2.39e}
i \partial_t w = Hw
\eeq

\noi with
\beq
\label{2.40e}
H = - (1/2) \Delta_K + {\check B}_S\ .
\eeq

Using the fact that the dynamical equations (\ref{2.34e}) (\ref{2.39e}) contain $\varphi$ only through its gradient $s = \nabla \varphi$, we can replace (\ref{2.38e}) by its gradient, thereby obtaining a closed system of equations for $(w, s, B_2)$, namely 
\beq
\label{2.41e}
\left \{ \begin{array}{l} i\partial_t w = Hw\\ \\  \partial_t s = t^{-1} \nabla g (w) + \nabla    {\check B}_{1L}\\ \\ B_2 = {\cal B}_2 (w, w, K) \end{array} \right .
\eeq

\noi with $H$, ${\check B}_{1L}$ and ${\cal B}_2$ defined by (\ref{2.40e}) (\ref{2.36e}) (\ref{2.35e}). This is the final form of the auxiliary system that replaces the original system (\ref{2.6e}) (\ref{2.7e}).\\

\noi {\bf Remark 2.2.} In I, aside from the change $t \to 1/t$, we have used a slightly different equation for $\varphi$, whereby the term $s^2 = |\nabla \varphi |^2$ appearing in the expansion of the covariant Laplacian $\Delta_K$ was included in the $\varphi$ equation, thereby making it into a Hamilton-Jacobi equation. The present choice is simpler since (i) the equation for $\varphi$ or $s$ is now immediately solved by integration over time and (ii) a number of terms are treated together in the form of the $K$ covariant derivative, which simplifies the algebra and thereby the subsequent estimates. On the other hand the previous more complicated choice gives slightly better decay estimates, more precisely a gain by $(1 - \ell n t)^2$ in the convergence estimates at zero (or infinity).\\

For technical reasons, in addition to the system (\ref{2.41e}), it will be useful to consider also a partly linearized system for $(w, B_2)$, namely 
\beq
\label{2.42e}
\left \{ \begin{array}{l} i\partial_t w' = Hw'\\ \\  B'_2 = {\cal B}_2 (w, w, K) \end{array} \right .
\eeq

\noi for new variables $(w', B'_2)$, where $H$ and $K$ still correspond to $(w, B_2)$. There is no point at this stage to introduce a new variable $s'$, since the equation for $s$ makes it an explicit function of $w$ (up to some suitable initial condition). \par

The problem of constructing the wave operators, namely of solving the Cauchy problem at infinity for the original system (\ref{2.6e}) (\ref{2.7e}) is now replaced by the problem of solving the Cauchy problem at $t=0$ for the auxiliary system (\ref{2.41e}). Since that system is singular at $t=0$, that cannot be done directly, and we follow instead the procedure sketched at the beginning of this section on the example of the original system (\ref{2.6e}) (\ref{2.7e}). We choose a set of asymptotic functions $v = (w_a, s_a, B_{1a}, B_{2a})$ which are expected to be suitable asymptotic forms of $(w, s, B_1, B_2)$ at $t = 0$, and we try to construct solutions of the auxiliary system (\ref{2.41e}) that are asymptotic to $v$ in a suitable sense at $t=0$. Note at this point that although $B_1$ is an explicit function of $w$, we refrain from assuming that $B_{1a} = B_1(w_a)$ in order to allow for more flexibility. Actually the final choice of $B_{1a}$ will differ from $B_1 (w_a)$. A similar remark applies to $s_a$. We also define 
\beq
\label{2.43e}
B_a = B_0 + B_{1a} + B_{2a} \qquad , \quad K_a = s_a + B_a \ .
\eeq

In particular $B_0$ is its own asymptotic form.\par

In order to solve the auxiliary system (\ref{2.41e}) with the previous asymptotic behaviour at $t=0$, we define the difference variables
\beq
\label{2.44e}
\left ( q, \sigma , G_1, G_2\right ) \equiv \left ( w-w_a, s-s_a, B_1 - B_{1a}, B_2 - B_{2a} \right ) \ .
\eeq

\noi We also define
\beq
\label{2.45e}
G = G_1 + G_2 \quad , \qquad L = \sigma + G \ ,
\eeq

\noi so that
\beq
\label{2.46e}
B = B_a + G \quad , \qquad K = K_a + L\ .
\eeq

\noi We define in addition 
\beq
\label{2.47e}
g(w_1, w_2) = (4 \pi |x|)^{-1} \ * \  {\rm Re} \ \overline{w}_1 w_2
\eeq

\noi so that $g(w) = g(w, w)$, and
\beq
\label{2.48e}
Q_{K_1} (K_2 , \cdot ) = K_2 \cdot \nabla_{K_1} + (1/2) (\nabla \cdot K_2)
\eeq

\noi so that 
\beq
\label{2.49e}
\Delta_{K_1 + K_2} = \Delta_{K_1} - 2 i Q_{K_1}(K_2, \cdot ) - K_2^2 \ .
\eeq

\noi The separation of $B_a$ and of $G$ into short range and long range parts follows the same pattern as that of $B$, namely
\beq
\label{2.50e}
\left \{ \begin{array}{l} {\check B}_{1aS} = \chi_S   {\check B}_{1a} \qquad , \quad {\check B}_{1aL} = \chi_L {\check B}_{1a}\ , \\ \\ {\check B}_{aS}  =  {\check B}_{0} + {\check B}_{1aS}   +  {\check B}_{2a} \qquad , \quad {\check B}_{aL} = {\check B}_{1aL}\ , \end{array} \right . 
\eeq

\beq
\label{2.51e}
\left \{ \begin{array}{l} {\check G}_{1S} = \chi_S   {\check G}_{1} \qquad , \quad {\check G}_{1L} = \chi_L {\check G}_{1}\ , \\ \\ {\check G}_{S}  =  {\check G}_{1S} + {\check G}_{2}    \qquad , \quad {\check G}_{L} = {\check G}_{1L} \end{array} \right . 
\eeq

\noi with $\chi_L$ and $\chi_S$ defined in (\ref{2.36e}).

Using the definitions  (\ref{2.44e})-(\ref{2.48e}), we rewrite the auxiliary system  (\ref{2.41e}) in terms of the difference variables. We take $(q, G_2)$ as independent dynamical variables and we consider $G_1$ and $\sigma$ as functions of $q$ defined by 
\beq
\label{2.52e}
G_1 = B_1 (q, 2w_a + q) - R_3
\eeq

\noi and by the equation for $\sigma$ that follows from the equation for $s$, namely 
\beq
\label{2.53e}
\partial_t \sigma = t^{-1} \nabla g (q, 2w_a+q) + \nabla  {\check G}_{1L} - R_2
\eeq

\noi with initial conditon $\sigma (0) = 0$. The auxiliary system for $(q, G_2)$ then becomes
\beq
\label{2.54e}
\left \{ \begin{array}{l} i \partial_t q = Hq - \widetilde{R}_1 \\ \\ G_2 = {\cal B}_2 (q, 2w_a + q, K) - t F_2 (PL |w_a|^2) - R_4 \end{array} \right . 
\eeq

\noi where
\bea
\label{2.55e}
&&\widetilde{R}_1  = R_1 - H_1 w_a \ ,\\
&&H_1 = i Q_{K_a} (L, \cdot ) + (1/2)L^2 + {\check G}_{S} 
\label{2.56e}
\eea

\noi and the remainders $R_j$, $1 \leq j \leq 4$ are defined by
\beq
\label{2.57e}
\left \{ \begin{array}{l} R_1 = i \partial_t w_a + (1/2) \Delta_{K_a} w_a - {\check B}_{aS} w_a\\ \\ R_2 = \partial_t s_a - t^{-1} \nabla g (w_a) - \nabla {\check B}_{1aL} \\ \\ R_3 = B_{1a} - B_1 (w_a)\\ \\ R_4 = B_{2a} - {\cal B}_2 (w_a, K_a)\ .\end{array} \right . 
\eeq

\noi In the equation for $G_2$ in  (\ref{2.54e}), we have used the identity 
$${\cal B}_2(w,K) = {\cal B}_2(w_a, K_a) + {\cal B}_2 (q, 2w_a + q, K) - tF_2 (PL|w_a|^2)\ .$$

The remainders $R_j$, $1 \leq j \leq 4$, express the failure of $(w_a, s_a, B_{1a}, B_{2a})$ to satisfy the original system (\ref{2.41e}). In order to solve that system, an essential condition will be that they tend to zero in a suitable sense as $t \to 0$. Their rate of convergence to zero measures  the quality of $(w_a, s_a, B_{1a}, B_{2a})$ as an asymptotic form for a solution of the system (\ref{2.41e}).\par

Again for technical reasons, we shall need a partly linearized version of the system (\ref{2.54e}) for the independent dynamical variables $(q, G_2)$. With 
\beq
\label{2.58e}
w' = w_a + q'\qquad , \quad B'_2 = B_{2a} + G'_2 \ ,
\eeq

\noi the linearized version of (\ref{2.54e}) corresponding to (\ref{2.42e}) becomes
\beq
\label{2.59e}
\left \{ \begin{array}{l} i \partial_t q' = Hq' - \widetilde{R}_1 \\ \\ G'_2 = {\cal B}_2 (q, 2w_a + q, K) - t F_2 (PL |w_a|^2) - R_4\ . \end{array} \right . 
\eeq

\noi Again there is no point in introducing new variables $G'_1$ and $\sigma '$, since $G_1$ and $\sigma$ (with the initial condition $\sigma (0) = 0$) are explicit functions of $q$.\par

The construction of solutions $(w, s, B_2)$ of the system (\ref{2.41e}) with prescribed asymptotic behaviour at $t= 0$ is now performed in two steps. The first step consists in solving the system (\ref{2.54e}) with $(q, G_2)$ tending to zero at $t=0$ under assumptions on $(w_a, s_a, B_{1a}, B_{2a})$ of a general nature, the most important of which being decay assumptions on the remainders $R_j$, $1 \leq j \leq 4$. This is done as follows. One first considers the linearized system (\ref{2.59e}) and one solves that system for $(q', G'_2)$ for given $(q, G_2)$ tending to zero at $t= 0$. This requires no work for $G'_2$ which is given by an explicit formula. As regards $q'$, we first solve the Cauchy problem for the relevant equation with initial condition $q'(t_0) = 0$ for some $t_0 > 0$ and we take the limit of the solution thereby obtained when $t_0 \to 0$. This procedure defines a map $\Gamma : (q, G_2) \to (q', G'_2)$. One then proves by a contraction method that the map $\Gamma$ has a fixed point in a suitable function space.\par

The second step of the method consists in constructing asymptotic functions satisfying the assumptions needed for the first step, and in particular the decay properties of the remainders. This will be done as follows. We shall take for $(w_a, s_a, B_{1a}, B_{2a})$ the second approximation in an iterative solution of the system (\ref{2.41e}) with the contribution of $B_0$ omitted, starting from $w_a (0) = w_+ = Fu_+$ where $u_+$ is the Schr\"odinger asymptotic state. This will be sufficient to control the $B_0$ independent terms in the remainders. On the other hand, the $B_0$ dependent terms will be controlled with the help of a support condition on $w_+$, as mentioned in the Introduction. The formulas are too complicated to be given here and are deferred to Section 6 below.\par

With the solution of the auxiliary systems (\ref{2.41e}) or (\ref{2.54e}) available, it is an easy matter to construct the modified wave operator for the original MS system (\ref{2.6e})  (\ref{2.7e}). The starting point is the asymptotic state $(u_+, A_+, \dot{A}_+)$ for the Schr\"odinger and Maxwell fields. One constructs the solution of the auxiliary system (\ref{2.41e}) as just explained, in a neighorhood of $t=0$. From $s$ one reconstructs the phase $\varphi$ by using (\ref{2.38e}). One finally substitutes $(w, \varphi , B_2)$ into (\ref{2.16e})  (\ref{2.17e})  (\ref{2.31e})-(\ref{2.33e}), thereby obtaining a solution $(u, A)$ of the system (\ref{2.6e}) (\ref{2.7e})  defined for large time. The modified wave operator is the map $(u_+, A_+, \dot{A}_+) \to (u, A)$ thereby obtained.\par

The main result of this paper is the construction of $(u, A)$ from $(u_+, A_+, \dot{A}_+)$ as described above, together with the asymptotic properties of $(u, A)$ that follow from that construction. It will be stated in full mathematical detail in Propositions 7.1 and 7.2 below. We give here only a heuristic preview of that result, stripped from most technicalities. We take $\beta = 1/2$ in (\ref{2.36e}) for definiteness.\\

\noi {\bf Proposition 2.1.} {\it Let $\beta = 1/2$. Let $u_+$ be such that $w_+ \equiv Fu_+ \in H^5$, $xw_+ \in H^4$ and let $w_+$ satisfy the support condition (\ref{6.90e}). Let $(A_+, \dot{A}_+)$ be sufficiently regular and decaying at infinity. Define $(w_a, s_a, B_a)$ by (\ref{2.43e}) (\ref{2.17e})$_0$ (\ref{2.9e}) and (\ref{6.1e}) (\ref{6.2e}) (\ref{6.3e}). Then \par

(1) There exists $\tau = \tau (u_+ , A_+, \dot{A}_+)$, $0 < \tau \leq 1$, such that the auxiliary system (\ref{2.41e}) has a unique solution $(w, s, B_2)$ in a suitable space, defined for $0 < t \leq \tau$, and such that $(w- w_a, s- s_a, B_2 - B_{2a})$ tends to zero in suitable norms when $t \to 0$.\par

(2) There exist $\varphi$ and $\varphi_a$ such that $s = \nabla \varphi$, $s_a = \nabla \varphi_a$ and such that $\varphi - \varphi_a$ tends to zero in suitable norms when $t \to 0$. Define $(u, A)$ by (\ref{2.16e}) (\ref{2.17e}) (\ref{2.31e})-(\ref{2.33e}). Then $(u, A)$ solves the system (\ref{2.6e}) (\ref{2.7e}) for $t \geq \tau^{-1}$ and $(u, A)$ behaves asymptotically as 
$$\left ( M(t) D(t) \exp (i \varphi_a (1/t)\right ) \overline{w}_a (1/t), A_0 (t) - t^{-1} D_0 (t) \left ( B_{1a} + B_{2a})(1/t)\right )$$
\noi in the sense that the difference tends to zero in suitable norms (for which each term separately is 0(1)) when $t \to \infty$.}\\

\noi {\bf Remark 2.3.} The suitable space quoted in Part (1) includes in particular the fact that $w \in {\cal C} (I, H^3) \cap {\cal C}^1 (I, H^1)$ and $xw \in {\cal C}(I, H^2) \cap {\cal C}^1 (I, L^2)$ with $I = (0, \tau ]$ so that the construction involves a loss of two space derivatives from $w_+$ to $w$. The relevant space for $B_2$ is slightly more complicated.\\

We now describe the contents of the technical parts of this paper, namely Sections 3-7. In Section 3, we introduce some notation, define the relevant function spaces and collect a number of preliminary estimates. In Section 4, we study the Cauchy problem for the auxiliary system (\ref{2.41e}). We solve the Cauchy problem for finite (non zero) initial time for the equation for $w$ of that system (Proposition 4.1) and we prove a uniqueness result under some (weak) decay of the solutions at $t=0$. In Section 5, we study the Cauchy problem at $t=0$ for the auxiliary system (\ref{2.54e}) under suitable general boundedness and decay assumptions on $(w_a, s_a, B_{1a}, B_{2a})$ and on the remainders $R_j$, $1 \leq j \leq 4$. We prove in particular the existence of a unique solution defined in some interval $(0, \tau ]$ for $\tau$ sufficiently small, tending to zero in suitable norms as $t$ tends to zero, first for the linearized system (\ref{2.59e}) (Proposition 5.1) and then for the non linear system (\ref{2.54e}) (Proposition 5.2). In Section 6, we construct asymptotic $(w_a, s_a, B_{1a}, B_{2a})$ satisfying the conditions required in Section~5. We derive in particular suitable bounds for those quantities (Proposition 6.1) and for the remainders (Propositions 6.2 and 6.3). In Section 7, we first collect the results of Sections 5 and 6 to derive the main result on the Cauchy problem for the auxiliary system (\ref{2.41e}) at $t=0$ (Proposition 7.1). Finally we construct the modified wave operators for the system (\ref{2.6e}) (\ref{2.7e}) from the results previously obtained for the system (\ref{2.41e})  and we derive the asymptotic estimates for the solutions $(u, A)$ in their range that follow from the previous estimates (Proposition 7.2).

\mysection{Notation and preliminary estimates}
\hspace*{\parindent} In this section we introduce some notation and collect a number of estimates which will be used throughout this paper. We denote by $\parallel \cdot \parallel_r$ the norm in $L^r \equiv L^r ({I\hskip-1truemm R}^n)$, to be used mostly in ${I\hskip-1truemm R}^3$, and by $<\cdot , \cdot >$ the scalar product in $L^2$. For any non negative integer $k$ and for $1 \leq r \leq \infty$ we denote by $H_r^k \equiv H_r^k ({I\hskip-1truemm R}^n)$ the Sobolev spaces 
$$H_r^k = \left \{ u \in {\cal S}'({I\hskip-1truemm R}^n):\ \parallel u ; H_r^k \parallel \ = \sum_{\alpha : 0 \leq |\alpha | \leq k} \parallel \partial^{\alpha} u\parallel_r\ < \infty \right \}$$

\noi where $\alpha$ is a multiindex. As a shorthand notation we will use 
$$\nabla^k = \left \{ \partial^{\alpha} :\  | \alpha | = k \right \} \ .$$

\noi For $1 < r < \infty$ those spaces can be defined equivalently (with equivalent norms) by
$$H_r^k = \left \{ u \in {\cal S}'({I\hskip-1truemm R}^n):\ \parallel u ; H_r^k \parallel \ = \ \parallel <\omega >^k u\parallel_r\ < \infty \right \}$$

\noi where $\omega = (- \Delta )^{1/2}$ and $<\cdot > = (1 + |\cdot |^2)^{1/2}$. The latter definition extends immediately to any $k \in {I\hskip-1truemm R}$ and we occasionally use such spaces. The subscript $r$ in $H_r^k$ will be omitted in the case $r= 2$. Besides the standard Sobolev spaces, we will use the associated homogeneous spaces $\dot{H}_r^k$ with norm $\parallel u ; \dot{H}_r^k\parallel \ = \ \parallel \omega^k u\parallel_r$. In particular it will be understood that $\dot{H}^1({I\hskip-1truemm R}^3) \subset L^6 ({I\hskip-1truemm R}^3)$. In addition we shall use the notation 
$$\ddot{H}^k = \dot{H}^1 \cap \dot{H}^k$$

\noi for any $k \geq 1$. For any Banach space $X \subset {\cal S}' ({I\hskip-1truemm R}^n)$ we denote by $<x>X$ the space defined by 
$$<x> X = \left \{ u \in {\cal S}' ({I\hskip-1truemm R}^n): <x>^{-1} u \in X \right \} \ .$$

\noi For any interval I and for any Banach space $X$ we denote by ${\cal C}(I, X)$ (resp. ${\cal C}_w (I, X)$) the space of strongly (resp. weakly) continuous functions from $I$ to $X$. For any positive integer $k$, we denote by ${\cal C}^k (I, X)$ the space of $k$ times differentiable functions from $I$ to $X$. For any $r$, $1 \leq r\leq \infty$, we denote by $L^r (I, X)$ (resp. $L_{loc}^r(I, X)$) the space of $L^r$ integrable (resp. locally $L^r$ integrable) functions from $I$ to $X$ if $r < \infty$, and the space of measurable essentially bounded (resp. locally essentially bounded) functions from $I$ to $X$ if $r = \infty$. For $I$ an open interval we denote by $D' (I, X)$ the space of vector-valued distributions from $I$ to $X$. We say that an evolution equation has a solution in $I$ with values in $X$ if the equation is satisfied in $D'(I_0, X)$ where $I_0$ is the interior of $I$. For $I$ a given interval, we denote by $(X, f)$ the set 
\beq
\label{3.1e}
(X, f) = \left \{ u \in {\cal C} (I, X) : \ \parallel u(t) ; X \parallel \ \leq f(t) \quad \forall\ t \in I \right \}
\eeq 

\noi where $X$ is a Banach space and $f \in {\cal C}(I, {I\hskip-1truemm R}^+)$. For real numbers $a$ and $b$ we use the notation $a \vee b = {\rm Max} (a, b)$ and $a \wedge b = {\rm Min} (a,b)$.\par

We shall use extensively the following Sobolev inequalities, stated here in ${I\hskip-1truemm R}^n$, but used only in ${I\hskip-1truemm R}^3$, and the following Leibnitz and commutator estimates.\\

\noi {\bf Lemma 3.1.} {\it (1) Let $1 < r \leq \infty$, $1 < r_1, r_2 < \infty$ and $0 \leq j < k$. If $r = \infty$, assume in addition that \ $k-j > n/r_2$. Let $\sigma$ satisfy \ $j/k \leq \sigma \leq 1$ and
$$n/r - j = (1 - \sigma )n/r_1 + \sigma (n/r_2 - k)\ .$$

\noi Then the following estimate holds~:
\beq
\label{3.2e}
\parallel \omega^j v \parallel_r\ \leq \ c \parallel v \parallel_{r_1}^{1-\sigma} \ \parallel \omega^k v \parallel_{r_2}^{\sigma}\ .
\eeq 

(2) Let $1 < r, r_1, r_3 < \infty$ and
$$1/r = 1/r_1 + 1/r_2 = 1/r_3 + 1/r_4\ .$$

\noi Then the following estimates hold
\beq
\label{3.3e}
\parallel \omega^k (v_1 v_2) \parallel_r\ \leq \ c\left ( \parallel \omega^k v_1 \parallel_{r_1}\  \parallel v_2  \parallel_{r_2}\ + \  \parallel \omega^k v_2  \parallel_{r_3}\  \parallel v_1  \parallel_{r_4}\right ) 
\eeq 

\noi for $k \geq 0$, and 
\beq
\label{3.4e}
 \parallel [\omega^k , v_1]v_2  \parallel_r\ \leq \ c\left ( \parallel \omega^k v_1 \parallel_{r_1}\  \parallel v_2  \parallel_{r_2}\ + \  \parallel \omega^{k-1} v_2  \parallel_{r_3}\  \parallel \nabla v_1  \parallel_{r_4}\right ) 
\eeq 

\noi for $k \geq 1$, where $[\ ,\ ]$ denotes the commutator.}\\

The proof of Lemma 3.1, part (1) follows from the Hardy-Littlewood-Sobolev inequality (\cite{16r}, p. 117) (from the Young inequality if $r = \infty$), from Paley-Littlewood theory and interpolation. The proof of Lemma 3.1, part (2) is given in \cite{17r} \cite{18r} with $\omega$ replaced by $<\omega >$ and follows therefrom by a scaling argument.\par

Occasionally a special case of (\ref{3.2e}) will be used with the ordinary derivative $\nabla$ replaced by the covariant derivative $\nabla_A = \nabla - iA$, where $A$ is a real vector-valued function, namely
\beq
\label{3.5e}
 \parallel v  \parallel _r\ \leq \ c  \parallel v  \parallel _{r_1}^{1 - \sigma} \  \parallel  \nabla_A v  \parallel _{r_2}^{\sigma}
\eeq 

\noi which holds under the assumptions of Lemma 3.1, part (1) with $j = 0$, $k=1$. The proof of (\ref{3.5e}) is an immediate consequence of (\ref{3.2e}) with $j = 0$, $k=1$ applied to $|v|$ and of the inequality $|\nabla |v| |\leq |\nabla_A v| $.\par

We shall also make use of the following three lemmas, stated here in ${I\hskip-1truemm R}^n$, but used only in ${I\hskip-1truemm R}^3$. \\

\noi {\bf Lemma 3.2.} {\it Let $k > 1 + n/2$, $v_j \in H^k$, $xv_j \in H^{k-1}$, $j = 1,2$. Then $xv_1 v_2 \in H^k$.}\\

\noi {\bf Proof.} From 
$$\omega^k x v_1 v_2 = [\omega^k, v_1] xv_2 + v_1 [\omega^k , x ] v_2 + x v_1 \omega^k v_2$$

\noi and
$$[x , \omega^k] = k \omega^{k-2} \nabla \ ,$$

\noi using Lemma 3.1, part (2), it follows that 
$$ \parallel \omega^k xv_1v_2 \parallel _2\ \leq \ c \Big (  \parallel \omega^k v_1  \parallel _2\  \parallel xv_2  \parallel _{\infty} \ + \  \parallel  \nabla v_1  \parallel _{\infty} \  \parallel  \omega^{k-1} xv_2  \parallel_2$$
$$+\  \parallel  v_1  \parallel _{\infty}\  \parallel  \omega^{k-1} v_2  \parallel _2 \ +  \parallel xv_1  \parallel _{\infty}\  \parallel  \omega^k v_2  \parallel _2 \Big ) \ . $$
\par\nobreak \hfill $\sq$\par

\noi {\bf Lemma 3.3.} {\it Let $P = \1 - \nabla \Delta^{-1} \nabla$. Let $k > 1 + n/2$, $v_j \in H^k$, $j = 1,2$. Then $P \ {\rm Im} \ \overline{v}_1 \nabla v_2 \in H^k$.}\\

\noi {\bf Proof.} One checks easily that 
\beq
\label{3.6e}
P\ {\rm Im} \ \overline{v}_1 \nabla v_2 = P\ {\rm Im} \ \overline{v}_2 \nabla v_1 \ .
\eeq 

\noi From
\begin{eqnarray*}
\omega^k P\ {\rm Im} \ \overline{v}_1\nabla v_2 &=& P\ {\rm Im} \ [\omega^k, \overline{v}_1] \nabla v_2 + P\ {\rm Im} \ \overline{v}_1\omega^k \nabla v_2\\
&=& P\ {\rm Im} \ [\omega^k, \overline{v}_1] \nabla v_2 + P\ {\rm Im} \ (\omega^k\overline{v}_2) \nabla v_1
\end{eqnarray*}

\noi using Lemma 3.1, part (2), it follows that 
\beq
\label{3.7e}
 \parallel  \omega^k P\ {\rm Im} \ \overline{v}_1\nabla v_2 \parallel_2\ \leq\ c \left (  \parallel  \omega^{k} v_1  \parallel_2\ \parallel  \nabla v_2 \parallel  _{\infty}\ + \ \parallel   \omega^{k} v_2  \parallel_2\ \parallel  \nabla v_1 \parallel  _{\infty}\right ) \ .
\eeq 
\nobreak \hfill $\sq$ \par

\noi {\bf Lemma 3.4.} {\it Let $k > 2 + n/2$, $v_1 \in L^{\infty} \cap \dot{H}^k$, $x\nabla v_1 \in L^{\infty}$, $v_2 \in H^k$, $xv_2 \in H^{k-1}$. Then $P_{ij} \ {\rm Im} \ x \overline{v}_2 \nabla_j v_1 \in H^k$.}\\

\noi {\bf Proof.} Let $w_i = P_{ij} \ {\rm Im} \ x \overline{v}_2 \nabla_j v_1$. Using (\ref{3.6e}) we obtain 
$$\parallel w_i \parallel_2\ \leq \parallel v_1 \parallel_{\infty} \ \parallel \nabla x v_2 \parallel_2\ .$$

\noi We estimate the norm $\dot{H}^k$ of $w_i$. We compute 
$$\nabla P_{ij} \ {\rm Im}\ x \overline{v}_2 \nabla_j v_1 = w_{1i} + w_{2i}$$

\noi where
\begin{eqnarray*}
&&w_{1i} = P_{ij}\ {\rm Im} \ (\nabla x \overline{v}_2) \nabla_j v_1 \\
&&w_{2i} = P_{ij}\ {\rm Im}\  (\nabla  v_1) \nabla_j x v_2
\end{eqnarray*}

\noi and, in the expression for $w_{2i}$, we have used again (\ref{3.6e}). This implies
$$ \parallel  \omega^k w_i \parallel_2 \ \leq  \ \parallel  \omega^{k-1} w_{1i} \parallel_2 \ + \  \parallel  \omega^{k-1} w_{2i}  \parallel_2 \ .$$

\noi From
$$\omega^{k-1}  w_{1i}  = P_{ij}\ {\rm Im} \left (  [ \omega^{k-1}, (\nabla_j v_1 )] \nabla x \overline{v}_2 + (\nabla_j v_1) [\omega^{k-1}\nabla, x] \overline{v}_2 + x(\nabla_j v_1) \omega^{k-1} \nabla \overline{v}_2 \right ) \ ,$$

\noi using Lemma 3.1, part (2), it follows that
$$\parallel  \omega^{k-1}  w_1 \parallel_2 \ \leq\ c \Big (  \parallel  \omega^k v_1\parallel_2 \ \parallel \nabla x v_2 \parallel_{\infty}\ + \ \parallel \nabla^2 v_1 \parallel_{\infty}\  \parallel  \omega^{k-1}  x v_1 \parallel_2$$
$$+\ \parallel \nabla v_1 \parallel_{\infty}\  \parallel  \omega^{k-1}   v_2 \parallel_2\ + \ \parallel x \nabla v_1 \parallel_{\infty}\  \parallel  \omega^{k}   v_2 \parallel_2\Big )\ . $$

\noi From the computation (\ref{3.7e}) in the proof of Lemma 3.3 we obtain 
$$\parallel  \omega^{k-1}  w_2 \parallel_2 \ \leq\ c \left  (  \parallel  \omega^k v_1\parallel_2 \ \parallel \nabla x v_2 \parallel_{\infty}\ + \ \parallel  \omega^{k-1}  x v_2 \parallel_2\ \parallel \nabla^2 v_2 \parallel_{\infty}\right ) \ . $$

\noi The result now follows from 
$$\parallel \nabla^m v_1 \parallel_{\infty}\ \leq\ c \parallel v_1 \parallel_{\infty}^{1-\sigma} \ \parallel \omega^k v_1 \parallel_2^{\sigma}$$

\noi for $m = 1,2$, with $\sigma = m (k - n/2)^{-1}$, using the Sobolev inequalities with integer derivatives and Mikhlin's theorem. \par \nobreak \hfill $\sq$ \par

We shall also need some estimates of the Hartree function $g$ defined by (\ref{2.47e}).\\

\noi {\bf Lemma 3.5.} {\it Let $w_j \in H^k$, $j = 1,2$, $k > 3/2$. Then $g(w_1,w_2) \in \ddot{H}^{k+2}$.}\\

\noi {\bf Proof.} The estimate
$$\parallel  \omega^{k+2} g(w_1,  w_2) \parallel_2 \ \leq\ c   \parallel  \omega^k \overline{w}_1 w_2\parallel_2$$

\noi is obvious, while the estimate
$$\parallel  g(w_1,  w_2) \parallel_6 \ \leq\ c  \parallel  \omega g(w_1, w_2) \parallel_2  \ \leq\ c  \parallel \overline{w}_1 w_2\parallel_{6/5}$$

\noi follows from the Hardy-Littlewood-Sobolev inequality. \par \nobreak \hfill $\sq$ \par

We now define the spaces where to look for solutions of the auxiliary system. For any interval $I \subset (0, 1]$, we denote by $X_0 (I)$ the Banach space 
\bea
\label{3.8e}
&&X_0(I) = \Big \{ (w, B_2): w\in {\cal C} (I, H^3) \cap {\cal C}^1 (I, H^1)  \ , \nn \\
&&xw\in {\cal C}(I, H^2) \cap {\cal C}^1 (I, L^2) ,\ B_2, {\check B}_2 \in {\cal C}(I, \dot{H}^1 \cap \dot{H}^2) \cap {\cal C}^1 (I, \dot{H}^1) \Big \}
\eea

\noi where ${\check B}_2$ is defined by (\ref{2.20e}). In order to take into account the time decay of the norms of the variables $q$ and $G_2$ (see (\ref{2.44e})) as $t$ tends to zero, we introduce a function $h \in {\cal C} (I, {I\hskip-1truemm R}^+)$ where $I = (0, \tau_0]$ for some $0 < \tau_0 \leq 1$, such that the function $\overline{h}(t) = t^{-3/2} h(t)$ be non decreasing in $I$ and satisfy
 \beq
\label{3.9e}
\int_0^t dt'\ t{'}^{-1} \ \overline{h} (t') \leq c \ \overline{h}(t)
\eeq 

\noi for some $c > 0$ and for all $t \in I$. A typical example of such an $h$ is $h(t) = t^{-3/2-\lambda}$, with $\lambda > 0$, which satisfies (\ref{3.9e}) with $c = \lambda^{-1}$. We then define the Banach space
\bea
\label{3.10e}
&&X(I) = \Big \{ (q, G_2) \in X_0 (I): \ \parallel (q, G_2); X(I) \parallel\ = \ \mathrel{\mathop {\rm Sup}_{t\in I}}\ h(t)^{-1}  \nn \\
&&\Big ( \parallel <x> q(t) \parallel_2 \vee \ t \left ( \parallel <x> \partial_t q(t) \parallel_2 \vee \parallel <x> \Delta q (t) \parallel_2 \right ) \nn\\
&&\vee\  t^{3/2} \left ( \parallel \nabla \partial_t q(t) \parallel_2 \vee \parallel \nabla \Delta q(t) \parallel_2 \right )\nn \\
&&\vee \ t^{-1/2}  \parallel \nabla G_2(t) \parallel_2 \vee \ t^{1/2}  \big ( \parallel \nabla^2 G_2(t) \parallel_2 \vee \parallel\nabla \partial_t G_2(t) \parallel_2 \vee \parallel \nabla {\check G}_2(t) \parallel_2 \big ) \nn \\
&& \vee \ t^{3/2} \big( \parallel \nabla^2 {\check G}_2(t) \parallel_2 \vee \parallel \nabla \partial_t {\check G}_2(t) \parallel \big )  \Big ) \Big \} \  .
\eea 

We next give some estimates for various components of $B_1$ expressed by (\ref{2.36e}). It follows immediately therefrom that 
\beq
\label{3.11e}
\parallel  \omega^m {\check B}_{1S} \parallel_2 \ \leq \ t^{\beta (p-m)} \parallel \omega^p {\check B}_{1S} \parallel_2\ \leq \ t^{\beta (p-m)} \parallel \omega^p {\check B}_{1} \parallel_2
\eeq 

\noi for $m \leq p$, and similarly
\beq
\label{3.12e}
\parallel  \omega^m {\check B}_{1L} \parallel_2 \ \leq \ (2t^{-\beta})^{m-p} \parallel \omega^p {\check B}_{1L} \parallel_2\ \leq \ (2t^{-\beta})^{m-p} \parallel \omega^p {\check B}_{1} \parallel_2
\eeq

\noi for $m \geq p$.\par

We now estimate $F_j (M)$ defined by (\ref{2.30e}) (\ref{2.15e}). From (\ref{2.30e}) it follows that 
\bea
\label{3.13e}
&&\omega F_j (M) = F_{j+1}( \omega M)\\
\label{3.14e}
&&\partial_t F_j (M) = F_{j+1} (\partial_t M) \\
&&x\cdot F_j (PM) = F_{j-1} (x \cdot PM)\ .
\label{3.15e}
\eea 

\noi The first two identities are obvious while in (\ref{3.15e}) we have used the identity
$$[x, f(\omega )] \cdot P = 0$$

\noi which holds for any regular function $f$. In addition, a direct computation yields
$$x \cdot PM = P (x \otimes M) - 2 \omega^{-2} \nabla \cdot M$$

\noi from which (\ref{3.15e}) can be continued to
\beq
\label{3.16e}
x \cdot F_j (PM) = F_{j-1} \left ( P (x \otimes M) - 2 \omega^{-2} \nabla \cdot M \right ) \ .
\eeq 

\noi In order to estimate $F_j$ we define
\beq
\label{3.17e}
I_j(f) (t) = \int_1^{\infty} d\nu\ \nu^{-j-3/2}\ f(t/\nu )
\eeq 

\noi for any $j \in {I\hskip-1truemm R}$ and for any non negative function $f$ in ${I\hskip-1truemm R}^+$. The estimates on $F_j$ are summarized in the following lemma.\\

\noi {\bf Lemma 3.6.} {\it For any $m,j \in {I\hskip-1truemm R}$ the following estimates hold~:}\par
\bea
\label{3.18e}
&{\it (1)} &\parallel  \omega^m F_j (M) \parallel_2 \ \leq c\ I_{j+m-2} \left ( \parallel \omega^{m-1} M\parallel_2\ \wedge \ \parallel \omega^m M\parallel_2 \right )\qquad\qquad\qquad \qquad  \\
& &\nn \\
&{\it (2)}  & \parallel  \omega^m x\cdot F_j (PM) \parallel_2 \ \leq c\ I_{j+m-3} \left ( \parallel <x> \omega^{m-1} M\parallel_2 \right ) 
\label{3.19e}
\eea

 {\it (3) For any $r$, $2 \leq r \leq 4$, 
 \beq
\label{3.20e}
\parallel  F_j (M) \parallel_r \ \leq c \int_1^{\infty} d\nu (\nu - 1)^{-1 + 2/r}\ \nu^{-j+1/r} \parallel M(t/\nu ) \parallel_{r_1}
\eeq 

\noi with $3/r_1 = 2 + 1/r$.}\\

\noi {\bf Proof.} \underline{Part (1).} From the definition of $F_j$ and from (\ref{3.13e}), from the identity
$$\parallel  \omega^m D_0 (\nu ) v \parallel_2 \ = \ \nu^{-m+3/2} \parallel  \omega^m v\parallel_2$$

\noi and from the estimate
$$|\sin \omega (\nu - 1)| \leq 1 \wedge \nu \omega$$ 

\noi we obtain easily (\ref{3.18e}).\\

\noi \underline{Part (2).} It is an immediate consequence of (\ref{3.16e}) and of Part (1).\\

\noi \underline{Part (3).} From the pointwise estimate \cite{1r} \cite{24r}
$$\parallel \sin \omega (\nu - 1) v \parallel_r\ \leq \ c (\nu - 1)^{-1 + 2/r} \parallel \omega^{2-4/r} v \parallel_{\overline{r}}$$

\noi with $2 \leq r < \infty$ and $1/r + 1/\overline{r} = 1$, it follows that
$$\parallel  F_j (M) \parallel_r \ \leq c \int_1^{\infty} d\nu (\nu - 1)^{-1 + 2/r}\ \nu^{-j+1/r} \parallel \omega^{1-4/r} M(t/\nu ) \parallel_{\overline{r}}$$

\noi which implies (\ref{3.20e}) by Lemma 3.1, part (1). \par \nobreak \hfill $\sq$ \par

The use of Lemma 3.6 leads to integrals of the type $I_j (t^{-\alpha} h)$ which are estimated by an elementary calculation as
\beq
\label{3.21e}
I_j (t^{- \alpha} h) = t^{3/2 - \alpha} I_{j+ 3/2 - \alpha } (\overline{h}) \leq c\ t^{-\alpha} h
\eeq 

\noi provided $\alpha \leq j+2$.

We now collect some properties of divergence free vector valued solutions of the wave equation $\sq A_0 = 0$ with initial data $(A_+, \dot{A}_+)$, given by (\ref{2.9e}). The divergence free condition $\nabla \cdot A_0 = 0$, obviously equivalent to $\nabla \cdot A_+ = \nabla \cdot \dot{A}_+=0$, implies that $\sq x \cdot A_0 = 0$, so that $x \cdot A_0$ can be written as 
\beq
\label{3.22e}
x \cdot A_0 (t) = (\cos \omega t) x \cdot A_+ + \omega^{-1} (\sin \omega t) x \cdot \dot{A}_+\ .
\eeq 

\noi We shall need the dilation operator
\beq
\label{3.23e}
S = t \partial_t + x \cdot \nabla + \1 \ .
\eeq 

\noi The operator $S$ satisfies the following commutation relations
\bea
&&\left [ S, \exp (i \omega t) \right ] = 0 \nn \\
&&S\omega^{-j} = \omega^{-j} (S + j)\nn \\
\label{3.24e}
&&S \partial_t^j \nabla^k = \partial_t^j \nabla^k (S - j - k)\\
&&SD_0 (t) = D_0 (t) \ t \partial_t
\label{3.25e}
\eea

\noi for any non negative integers $j$, $k$. It follows from (\ref{3.24e}) that $\sq SA_0 = 0 $ so that $SA_0$ can be written as
\beq
\label{3.26e}
\left ( SA_0\right ) (t) = (\cos \omega t) (1 + x \cdot \nabla )A_+ + \omega^{-1} (\sin \omega t) (2 + x \cdot \nabla ) \dot{A}_+\ .
\eeq 

\noi When changing variables from $A_0$ to $B_0$ according to (\ref{2.17e})$_0$ we obtain 
\beq
\label{3.27e}
\left ( x \cdot A_0\right  )(t) = - t^{-1} D_0 (t) \ {\check B}_0 (1/t)
\eeq 

\noi and
\beq
\label{3.28e}
 \nabla^k S^j A_0(t) = (-1)^{j+1}\ t^{-1-k} D_0(t) \left ( \nabla^k (t \partial_t )^j B_0 \right ) (1/t)\ .
\eeq 

\noi We finally collect some estimates of divergence free vector solutions of the wave equation.\\

\noi {\bf Lemma 3.7.} {\it Let $j$, $k$ be non negative integers. Assume $(A_+, \dot{A}_+)$ to be divergence free and to satisfy the conditions
\beq
\label{3.29e}
\left \{ \begin{array}{ll} {\cal A} \in L^2 &\nabla^2 {\cal A} \in L^1\\ \\ \omega^{-1} \dot{\cal A} \in L^2\qquad\qquad\qquad &\nabla \dot{\cal A} \in L^1 \end{array} \right .
\eeq 

\noi for 
\beq
\label{3.30e}
\left \{ \begin{array}{ll} {\cal A} = (x \cdot \nabla )^{j'} \nabla^k A_+ &{\cal A} = (x \cdot \nabla )^{j'} \nabla^k (x \cdot A_+)\\ \\ \dot{\cal A}  = (x \cdot \nabla )^{j'} \nabla^k \dot{A}_+\qquad\qquad\qquad &\dot{\cal A}  = (x \cdot \nabla )^{j'} \nabla^k (x \cdot \dot{A}_+) \end{array} \right .
\eeq 

\noi for $0 \leq j' \leq j$. Then $A_0$ satisfies the following estimates~:
\beq
\label{3.31e}
\parallel (S^{j'} \nabla^k A_0) (t)\parallel_r \vee \parallel (S^{j'} \nabla^k x \cdot A_0 ) (t) \parallel_r\ \leq b_0 \ t^{-1+2/r}
\eeq 

\noi for $0 \leq j' \leq j$, for $2 \leq r \leq \infty$ and for all $t > 0$.\par

Let $B_0$ and ${\check B}_0$ be defined by (\ref{2.17e})$_0$ and (\ref{2.20e}). Then $B_0$ and ${\check B}_0$ satisfy the following estimates~:
\beq
\label{3.32e}
\parallel \partial_t^j\nabla^k B_0 (t)\parallel_r \vee \parallel \partial_t^j \nabla^k {\check B}_0  (t) \parallel_r\ \leq b_0 \ t^{-j-k+1/r} 
\eeq 
\noi for $2 \leq r \leq \infty$ and for all $t > 0$.}\\

\noi {\bf Proof.} The estimate (\ref{3.31e}) is standard for $j = 0$ \cite{27r}. For $j \not= 0$ it is a consequence of the case $j = 0$ and of the commutation relations satisfied by the operator $S$. The estimate (\ref{3.32e}) follows from (\ref{3.27e}), (\ref{3.28e}) and (\ref{3.31e}).\par \nobreak \hfill $\sq$ \par

In all the estimates in this paper we denote by $C$ a constant depending on the asymptotic functions $(w_a, K_a)$ through the available norms. Absolute constants will be in general omitted, except in special arguments where they are explicitly needed, in which case they are denoted by $c$. The letters $j$, $k$, $\ell$ will always denote non negative integers.

\mysection{Miscellaneous results on the auxiliary system}
\hspace*{\parindent} In this section we study the Cauchy problem for the auxiliary system (\ref{2.41e}) and its variants (\ref{2.42e}) (\ref{2.54e}) (\ref{2.59e}). We first solve the Cauchy problem with finite (non zero) initial time for the Schr\"odinger equation which occurs in (\ref{2.42e}) (\ref{2.59e}). We then derive a uniqueness result for the Cauchy problem at time zero for the non linear system (\ref{2.41e}). That result will eventually apply also to the system (\ref{2.54e}).

We rewrite the Schr\"odinger equation which occurs in (\ref{2.42e}) or (\ref{2.59e}) in the general form  
\beq
\label{4.1e}
i \partial_t v = - (1/2) \Delta_K v + V v + f_0 \equiv Hv + f_0
\eeq

\noi where $\Delta_K$ is the covariant Laplacian associated with a real vector field $K$ and $V$ is a real function. We first state the result on the Cauchy problem  at finite time for that equation. That result is a variant of Proposition 3.3 of \cite{9r}.\\

\noi {\bf Proposition 4.1.} {\it Let $0 < \tau \leq 1$, let $I = (0, \tau ]$ and let $t_0 \in I$. Let $K,V \in {\cal C}(I, \dot{H}^1 \cap \dot{H}^2) \cap {\cal C}^1 (I, \dot{H}^1)$ with $\partial_t \nabla \cdot K \in {\cal C}(I, \dot{H}^1)$. \par

(1) Let $f_0 \in L_{loc}^1 (I, H^{-1})$ and $v_0 \in H^1$. Then there exists at most one solution $v \in (L_{loc}^{\infty} \cap {\cal C}_w)(I, H^1)$ of (\ref{4.1e}) with $v(t_0) = v_0$. \par

(2) Let in addition $f_0 \in {\cal C}^1 (I, H^1)$ and $v_0 \in H^3$. Then there exists a unique solution $v \in {\cal C}(I, H^3) \cap {\cal C}^1 (I, H^1)$ of (\ref{4.1e}) with $v(t_0)= v_0$. That solution satisfies the conservation laws
\beq
\label{4.2e}
\parallel v(t_2) \parallel_2^2 \ - \ \parallel v(t_1) \parallel_2^2 \ = \int_{t_1}^{t_2} dt\ 2\  {\rm Im} <v,f_0>(t) \ ,
\eeq
\beq
\label{4.3e}
\parallel \partial_t v(t_2) \parallel_2^2 \ - \ \parallel \partial_t v(t_1) \parallel_2^2 \ = \int_{t_1}^{t_2} dt\ 2\  {\rm Im} <\partial_tv,f_2>(t) \ ,
\eeq
\beq
\label{4.4e}
\parallel \nabla_K \partial_t v(t_2) \parallel_2^2 \ - \ \parallel \nabla_K \partial_t v(t_1) \parallel_2^2 \ = \int_{t_1}^{t_2} dt\ 2\  {\rm Im} <\nabla_K\partial_tv,f_3>(t) 
\eeq

\noi for all $t_1, t_2 \in I$, where
\beq
\label{4.5e}
f_2 = \left ( \partial_t H \right ) v + \partial_t f_0\ ,
\eeq
\beq
\label{4.6e}
\partial_t  H = i \left ( \partial_t K \right )\cdot \nabla_K + (i/2) \left ( \partial_t \nabla \cdot K \right ) + \left ( \partial_t V \right )\ ,
\eeq
\beq
\label{4.7e}
f_3 = \left ( \partial_t K + \nabla V \right )  \partial_t v + \nabla_K f_2 \  .
\eeq

(3) Let in addition $<x> f_0 \in {\cal C}^1(I, L^2)$ and $<x> v_0 \in H^2$. Then the solution $v$ satisfies $<x>v\in {\cal C} (I, H^2) \cap {\cal C}^1 (I, L^2)$ and $v$ satisfies the conservation laws
\beq
\label{4.8e}
\parallel x v(t_2) \parallel_2^2 \ - \ \parallel x v(t_1) \parallel_2^2 \ = \int_{t_1}^{t_2} dt\ 2 \ {\rm Im} <xv,xf_0 + \nabla_K v >(t) \ ,
\eeq
\beq
\label{4.9e}
\parallel x \partial_t v(t_2) \parallel_2^2 \ - \ \parallel x\partial_t v(t_1) \parallel_2^2 \ = \int_{t_1}^{t_2} dt\ 2 \ {\rm Im} <x\partial_t v,xf_2 + \nabla_K \partial_t v>(t) 
\eeq

\noi for all $t_1, t_2 \in I$.}\\

\noi {\bf Proof.}

\noi \underline{Part (1).} It follows easily from the regularity assumptions on $v$ that the $L^2$ norm of the difference of two solutions is constant in time. This immediately implies uniqueness.\\

\noi \underline{Part (2).} The proof proceeds by a parabolic regularization and a limiting procedure. We consider separately the cases $t \geq t_0$ and $t \leq t_0$ and we begin with $t \geq t_0$. The regularization is parametrized by a parameter $\eta$ with $0 < \eta \leq 1$, and in addition we regularize $K$ in space by the use of a standard mollifier parametrized by $\eta$ so that the regularized $K_{\eta}$ belongs to ${\cal C}^1(I, \dot{H}^N)$ for any $N \geq 1$. The regularized equation is then 
\beq
\label{4.10e}
i \partial_t v = H_{\eta} v + f_0
\eeq

\noi where
\beq
\label{4.11e}
H_{\eta} = - (1/2) (1 - i \eta ) \Delta_{K_{\eta}} + V = H_{0 \eta} + A_{\eta} \ ,
\eeq
\beq
\label{4.12e}
H_{0\eta} = - (1/2) (1 - i \eta ) \Delta
\eeq
\beq
\label{4.13e}
A_{\eta} = (1 - i \eta ) \left ( i K_{\eta} \cdot \nabla + (i/2) (\nabla \cdot  K_{\eta} ) + (1/2) K_{\eta}^2\right ) + V \ .
\eeq

\noi The equation (\ref{4.10e}) can then be rewritten as 
\beq
\label{4.14e}
i \partial_t v = H_{0\eta} v + \widetilde{f}_{\eta}
\eeq

\noi where
\beq
\label{4.15e}
 \widetilde{f}_{\eta} = A_{\eta}v + f_0 \ .
\eeq

\noi Define
\beq
\label{4.16e}
U_{\eta}(t) = \exp \left ( - i t \ H_{0\eta}\right ) \ .
\eeq

The Cauchy problem for (\ref{4.10e}) with initial condition $v(t_0) = v_0$ can then be rewritten as a fixed point problem for the map $\phi$ defined by 
\beq
\label{4.17e}
(\phi (v)) (t) = U_{\eta} (t-t_0)v_0 - i \int_{t_0}^t dt'\ U_{\eta} (t-t')  \widetilde{f}_{\eta}(t')\ .
\eeq

\noi The map $\phi$ is easily seen to satisfy the following identities
\beq
\label{4.18e}
i \partial_t \phi (v) = H_{0\eta}  \ \phi (v) +  \widetilde{f}_{\eta} \ ,
\eeq
\beq
\label{4.19e}
\left ( i \partial_t \phi (v)\right ) (t) = U_{\eta} (t-t_0)\left ( H_{0\eta} v_0 + \widetilde{f}_{\eta}(t_0) \right ) + \int_{t_0}^t dt' U_{\eta} (t-t')\partial_t \widetilde{f}_{\eta}(t')\ .
\eeq

\noi We next show that the map $\phi$ is a contraction in ${\cal C} (J, H^3)\cap {\cal C}^1 (J, H^1)$ for $J = [t_0 , t_0 + T]$ for $T$ sufficiently small. For that purpose we estimate
\beq
\label{4.20e}
\parallel (\phi (v)) (t) \parallel _2\ \leq \ \parallel v_0 \parallel _2\ +  \int_{t_0}^t dt' \parallel \widetilde{f}_{\eta}(t')\parallel_2 \ ,
\eeq
\beq
\label{4.21e}
\parallel \left (  \partial_t \phi (v)\right ) (t) \parallel_2 \ \leq \ \parallel H_{0\eta} v_0\parallel_2\ + \ \parallel \widetilde{f}_{\eta}(t_0) \parallel_2\  + \int_{t_0}^t dt' \parallel \partial_t \widetilde{f}_{\eta}(t')\parallel_2 \ ,
\eeq
\beq
\label{4.22e}
\parallel \left (  \nabla \partial_t \phi (v)\right ) (t) \parallel_2 \ \leq \ \parallel \nabla H_{0\eta} v_0\parallel_2\ + \ \parallel \nabla \widetilde{f}_{\eta}(t_0) \parallel_2\  + \eta^{-1/2} \int_{t_0}^t dt'  (t-t')^{-1/2}\parallel \partial_t \widetilde{f}_{\eta}(t')\parallel_2 \ ,
\eeq
\beq
\label{4.23e}
\parallel \nabla \Delta \phi (v) \parallel _2\ \leq \ 2 \left ( \parallel \nabla \partial_t \phi (v) \parallel _2\ +  \ \parallel \nabla \widetilde{f}_{\eta}\parallel_2 \right )
\eeq

\noi where in (\ref{4.22e}) we have used the estimate
$$\parallel \nabla U_{\eta}(t) ; {\cal B}(L^2) \parallel \ \leq ( \eta t)^{-1/2}  \ .$$

\noi In order to estimate the relevant norms of $\widetilde{f}_{\eta}$, we need estimates of $A_{\eta}$, $\partial_t A_{\eta}$ and $\nabla A_{\eta}$ and for later purposes we also consider $\nabla \partial_t A_{\eta}$. From (\ref{4.13e}) we obtain 
\beq
\label{4.24e}
\partial_t A_{\eta} = (1 - i \eta )i\  Q_{K_{\eta}} \left ( \partial_t K_{\eta} , \cdot \right ) + \partial_t V\ ,
\eeq
\beq
\label{4.25e}
\nabla  A_{\eta} = (1 - i \eta )i\  Q_{K_{\eta}} \left ( \nabla K_{\eta} , \cdot \right ) + \nabla V\ ,
\eeq
\beq
\label{4.26e}
\nabla  \partial_t A_{\eta} = (1 - i \eta )i\  Q_{K_{\eta}} \left ( \nabla \partial_t K_{\eta} , \cdot \right ) + \left (  \partial_t K_{\eta} \right ) \left (  \nabla K_{\eta} \right ) + \nabla \partial_t V \ ,
\eeq

\noi where
$$Q_K (L, \cdot ) = L \cdot \nabla_K + (1/2) (\nabla \cdot K) \ .$$

\noi We estimate
\beq
\label{4.27e}
\parallel A_{\eta} v\parallel_2 \ \leq \ \left ( \parallel K \parallel_{\infty} \ + \ \parallel \nabla \cdot K \parallel_3 \right ) \parallel \nabla v \parallel_2 \ + \left (  \parallel K \parallel_{\infty}^2\ + \ \parallel V \parallel_{\infty} \right ) \parallel v \parallel_2 \ ,
\eeq
\bea
\label{4.28e}
\parallel \left ( \partial_t A_{\eta} \right ) v\parallel_2 &\leq& \parallel \partial_t  K \parallel_{6} \  \left (  \parallel \nabla v \parallel_3 + \parallel K \parallel_{\infty} \ \parallel v \parallel_{3}\right ) \nn \\
&&+ \left ( \parallel \partial_t \nabla \cdot K  \parallel_6 \ + \  \parallel \partial_t V \parallel_6\right )  \parallel v \parallel_3 \ ,
\eea
\bea
\label{4.29e}
\parallel \left ( \nabla A_{\eta} \right ) v\parallel_2 &\leq& \parallel \nabla  K \parallel_{3} \  \left (  \parallel \nabla v \parallel_6 + \parallel K \parallel_{\infty} \ \parallel v \parallel_{6}\right ) \nn \\
&&+ \left ( \parallel \nabla \nabla \cdot K  \parallel_2 \ + \  \parallel \nabla V \parallel_2\right )  \parallel v \parallel_{\infty} \ ,
\eea
\bea
\label{4.30e}
&&\parallel \left ( \nabla \partial_t A_{\eta} \right ) v\parallel_2 \ \leq \ \parallel \nabla \partial_t  K \parallel_{2} \  \left (  \parallel \nabla v \parallel_{\infty} + \parallel K \parallel_{\infty} \ \parallel v \parallel_{\infty}\right ) \nn \\
&&+ \left ( \parallel \nabla \partial_t \nabla \cdot K  \parallel_2 \ + \  \parallel \partial_t K \parallel_6\ \parallel \nabla K  \parallel_3\ + \ \parallel \nabla \partial_t v  \parallel_2\right )  \parallel v \parallel_{\infty} \ ,
\eea

\noi where we have used the fact that the norms of $K_{\eta}$ are bounded by those of $K$ whenever the latter are finite. Using (\ref{4.27e}) (\ref{4.28e}) (\ref{4.29e}), we estimate 
\beq
\label{4.31e}
\parallel \widetilde{f}_{\eta} \parallel_2\ \leq \ \parallel A_{\eta} v\parallel_2\ + \ \parallel  f_0 \parallel _2\ \leq C \left ( \parallel \nabla v \parallel _2 \ + \ \parallel  v \parallel _2 \ + 1 \right ) \ ,
\eeq
\bea
\label{4.32e}
&&\parallel \partial_t  \widetilde{f}_{\eta} \parallel_2\ \leq \ \parallel A_{\eta} \partial_t  v\parallel_2\ + \  \parallel\left ( \partial_t A_{\eta} \right ) v \parallel_2\ + \ \parallel  \partial_t  f_0 \parallel _2\nn \\
&&\ \leq C \Big ( \parallel \nabla \partial_t v \parallel _2 \ + \ \parallel  \partial_t  v \parallel _2 \ + \ \parallel  \nabla  v \parallel _3\ + \ \parallel  v \parallel _3\ + 1 \Big ) \ ,
\eea
\bea
\label{4.33e}
&&\parallel \nabla \widetilde{f}_{\eta} \parallel_2\ \leq \ \parallel A_{\eta} \nabla v\parallel_2\ + \  \parallel\left ( \nabla A_{\eta} \right ) v \parallel_2\ + \ \parallel  \nabla f_0 \parallel _2\nn \\
&&\ \leq C \Big ( \parallel \Delta v \parallel _2 \ + \ \parallel  \nabla  v \parallel _2 \ + \ \parallel  v \parallel _{\infty}\ +  1 \Big )
\eea

\noi where the constants $C$ depend on $K$, $V$ and $f_0$.\par

From (\ref{4.20e})-(\ref{4.23e}) and (\ref{4.31e})-(\ref{4.33e}), it follows that $\phi$ is a contraction in ${\cal C} (J, H^3) \cap {\cal C}^1 (J, H^1)$ for $T$ sufficiently small, so that (\ref{4.10e}) with initial condition $v(t_0) = v_0$ has a unique solution in that space. By standard arguments using the linearity of the equation, one can extend the solution to $I_+ = [t_0 , \tau ]$. Let $v_{\eta}$ be that solution.\par

We next take the limit $\eta \to 0$ and for that purpose we need estimates of $v_{\eta}$ uniform in $\eta$ in the relevant norms. Those estimates will follow from the conservation laws satisfied by $v_{\eta}$, which are the regularized version of (\ref{4.2e}) (\ref{4.3e}) (\ref{4.4e}), with regularized $K$, $f_2$ and $f_3$. In fact from (\ref{4.10e}) it follows formally that
\beq
\label{4.34e}
\partial_t \parallel v_{\eta} \parallel_2^2 \ = - \eta \parallel \nabla_{K_{\eta}}\  v_{\eta}\parallel_2^2 \ + 2\  {\rm Im} <v_{\eta}, f_0>
\eeq

\noi so that 
\beq
\label{4.35e}
\partial_t \parallel v_{\eta} \parallel_2\ \leq \ \parallel f_0 \parallel_2 \ .
\eeq

\noi Taking the time derivative of  (\ref{4.10e}) then yields
\beq
\label{4.36e}
i \partial_t \partial_t v_{\eta} = H_{\eta} \partial_t  v_{\eta}  + \left ( \partial_t A_{\eta} \right ) v_{\eta}  + \partial_t f_0
\eeq

\noi which implies
\beq
\label{4.37e}
\partial_t \parallel \partial_t v_{\eta}\parallel_2^2\  = - \eta \parallel \nabla_{K_{\eta}} \partial_t  v_{\eta}\parallel_2^2  + 2\ {\rm Im} <\partial_t v_{\eta}, f_{2\eta}> 
\eeq

\noi and therefore
\beq
\label{4.38e}
\partial_t \parallel \partial_t v_{\eta}\parallel_2\  \leq \  \parallel f_{2\eta} \parallel_2
\eeq

\noi where
\beq
\label{4.39e}
f_{2\eta} = \left ( \partial_t A_{\eta}\right ) v_{\eta} + \partial_t f_0 = i (1 - i \eta ) \left ( Q_{K_{\eta}} \left ( \partial_t K_{\eta}, v_{\eta}\right ) \right ) + \left ( \partial_t V \right ) v_{\eta} + \partial_t f_0 \ .
\eeq

\noi Taking the covariant gradient of (\ref{4.36e}) yields
\beq
\label{4.40e}
i \partial_t \nabla_{K_{\eta}} \partial_t v_{\eta} = - (1/2) (1 - i \eta ) \nabla_{K_{\eta}} \Delta_{K_{\eta}}\partial_t v_{\eta} + V  \nabla_{K_{\eta}} \partial_t v_{\eta} + f_{3\eta}
\eeq

\noi so that
\beq
\label{4.41e}
\partial_t \parallel \nabla_{K_{\eta}} \partial_t v_{\eta} \parallel_2^2\ = \ - \eta  \parallel \Delta_{K_{\eta}} \partial_t v_{\eta} \parallel_2^2\  + 2\ {\rm Im} < \nabla_{K_{\eta}} \partial_t v_{\eta}, f_{3\eta}> 
\eeq

\noi and therefore
\beq
\label{4.42e}
\partial_t \parallel \nabla_{K_{\eta}} \partial_t v_{\eta} \parallel_2\ \leq \ \parallel   f_{3\eta} \parallel _2
\eeq

\noi where
\beq
\label{4.43e}
 f_{3\eta} = \left ( \partial_t K_{\eta} + \nabla V \right ) \partial_t v_{\eta} +  \nabla_{K_{\eta}} f_{2\eta}\ .
\eeq

The conservation laws  (\ref{4.34e})   (\ref{4.37e})   (\ref{4.41e}) for $v_{\eta}$ are in fact derived in an integral form similar to  (\ref{4.2e})  (\ref{4.3e})  (\ref{4.4e}). The proof is immediate for  (\ref{4.34e})  (\ref{4.37e}) under the available regularity properties of $v_{\eta}$, while it requires a more delicate argument for  (\ref{4.41e}) involving an additional regularization and a limiting procedure. It is at this point that we need the $\eta$ regularization of $K$. We refer to \cite{2r} \cite{3r} for details on that problem.\par

We now estimate $v_{\eta}$. We define
\beq
\label{4.44e}
y_0 = \ \parallel v_{\eta} \parallel_2\quad , \qquad y_2 = \ \parallel \partial_t v_{\eta} \parallel_2\quad , \qquad y_3 = \ \parallel \nabla_{K_{\eta}}  \partial_t v_{\eta} \parallel_2\ .
\eeq

\noi From (\ref{4.35e}) we obtain
\beq
\label{4.45e}
y_0 \leq \ \parallel v_0 \parallel_2\ + \int_{t_0}^t dt'  \parallel f_0(t')  \parallel_2\ \leq C\ .
\eeq

\noi From (\ref{4.14e})  (\ref{4.15e})  (\ref{4.27e}) we obtain
\begin{eqnarray*}
\parallel \Delta v_{\eta} \parallel_2 &\leq& y_2 + C \left ( \parallel \nabla v_{\eta} \parallel_2\ + \ \parallel v_{\eta} \parallel_2\right ) + \ \parallel f_0 \parallel _2 \\
&\leq& y_2 + C \left ( y_0^{1/2} \parallel \Delta v_{\eta} \parallel_2^{1/2}\ + y_0 + 1 \right )
\end{eqnarray*}

\noi which implies
\beq
\label{4.46e}
\parallel \Delta v_{\eta} \parallel_2\ \leq y_2 + C \ .
\eeq

\noi From (\ref{4.38e})  (\ref{4.39e})  (\ref{4.28e}) we next obtain
\bea
\label{4.47e}
\partial_t y_2 \leq \ \parallel f_{2\eta} \parallel_2&\leq& C \left ( \parallel \nabla v_{\eta} \parallel_3\ + \ \parallel v_{\eta} \parallel_3\right ) + \ \parallel \partial_t  f_0\parallel_2 \nn \\
&\leq& C \left ( y_0^{1/4} \parallel \Delta v_{\eta} \parallel_2^{3/4}\ + \ y_0^{3/4} \parallel \Delta v_{\eta} \parallel_2^{1/4} \ + 1\right ) \nn \\
&\leq& C \left ( y_2^{3/4} + 1 \right )
\eea

\noi by (\ref{4.45e})  (\ref{4.46e}), which implies
\beq
\label{4.48e}
y_2 \leq C \ .
\eeq

\noi We next estimate $y_3$. From (\ref{4.42e})  (\ref{4.43e}) (\ref{4.39e}) we obtain
$$\partial_t y_3 \leq \ \parallel f_{3\eta} \parallel_2\ \leq \  \parallel \partial_t  K_{\eta} + \nabla V \parallel_6\ \parallel \partial_t v_{\eta} \parallel_3 \ +\  \parallel   \left ( \partial_t  A_{\eta}\right ) \nabla  v_{\eta}   \parallel_2$$
\beq
\label{4.49e}
+ \  \parallel \left ( \nabla \partial_t A_{\eta}\right )  v_{\eta}\parallel_2\ + \  \parallel K_{\eta}\parallel_{\infty}\ \parallel f_{2\eta} \parallel_2 \ + \ \parallel \nabla \partial_t f_0 \parallel_2 \nn\\
\eeq

\noi so that by (\ref{4.28e})  (\ref{4.30e}) and  (\ref{4.46e})  (\ref{4.47e})  (\ref{4.48e})
\bea
\label{4.50e}
&&\partial_t y_3 \leq \ \parallel f_{3\eta} \parallel_2\ \leq \  C\Big ( \parallel \partial_t  v_{\eta}  \parallel_3\ + \ \parallel \nabla^2 v_{\eta} \parallel_3 \ + \ \parallel \nabla v_{\eta}\parallel_3 \nn \\
&&+\  \parallel    \nabla  v_{\eta} \parallel_{\infty} \ + \  \parallel v_{\eta} \parallel_{\infty}\ + 1\Big ) \  +\ \parallel \nabla \partial_t f_0\parallel_2 \nn \\
&&\leq C \Big (  \parallel \partial_t v_{\eta} \parallel_3 \ + \ \parallel \Delta  v_{\eta}\parallel_3\ + \   \parallel \nabla v_{\eta} \parallel_{\infty} + 1 \Big ) \ . 
\eea

\noi Using a slight extension of (\ref{4.27e}), namely
$$ \parallel    A_{\eta}   v_{\eta} \parallel_r\ \leq \  \parallel   K  \parallel_{\infty}\  \parallel  \nabla  v_{\eta} \parallel_r\ + \  \parallel  \nabla \cdot K \parallel_6\ \parallel  v_{\eta} \parallel_{r_1} \ + \left ( \parallel K \parallel_{\infty}^2 \ + \ \parallel V \parallel_{\infty} \right ) \parallel  v_{\eta} \parallel_r$$

\noi with $2 \leq r \leq 6$ and $1/r_1 = 1/r - 1/6$, we estimate
\bea
\label{4.51e}
 \parallel  \Delta  v_{\eta}   \parallel_r &\leq&  \parallel   \partial_t v_{\eta}   \parallel_r \ + \ C\left (  \parallel   \nabla v_{\eta}   \parallel_r\ + \  \parallel   v_{\eta}   \parallel_{r_1} \ + \  \parallel   v_{\eta}   \parallel_r \right ) + \  \parallel   f_0   \parallel_r\nn \\
 &\leq&  \parallel   \partial_t v_{\eta}   \parallel_r\ + C \ \leq C \left ( y_3^{\delta} + 1 \right )
 \eea

\noi with $0 \leq \delta = 3/2 - 3/r \leq 1$, by (\ref{4.46e}) (\ref{4.48e}) and a covariant Sobolev inequality. From (\ref{4.50e}) (\ref{4.51e}) we obtain
\beq
\label{4.52e}
\partial_t y_3 \leq C  \left ( y_3^{1/2} + 1 \right )
\eeq

\noi which implies
\beq
\label{4.53e}
y_3 \leq C 
\eeq

\noi and therefore by (\ref{4.48e})
\beq
\label{4.54e}
\parallel   \nabla \partial_t v_{\eta}   \parallel_2\ \leq C \ .
\eeq

\noi Finally taking the gradient of (\ref{4.14e}), we estimate
\beq
\label{4.55e}
\parallel   \nabla \Delta v_{\eta}   \parallel_2\ \leq 2 \left ( \parallel   \nabla \partial_t v_{\eta}   \parallel_2\ + \ \parallel   \nabla ( A_{\eta}  v) \parallel_2\ + \ \parallel   \nabla f_0  \parallel_2\right ) \ \leq C
\eeq

\noi by (\ref{4.54e}) (\ref{4.27e}) (\ref{4.29e}) and (\ref{4.46e}) (\ref{4.48e}). The estimates (\ref{4.45e}) (\ref{4.48e}) (\ref{4.54e}) (\ref{4.55e}) show that $v_{\eta}$ is estimated in ${\cal C} (I_+ , H^3) \cap {\cal C}^1 (I_+ , H^1)$ uniformly in $\eta$. \par

We now take the limit $\eta \to 0$. Let $0 < \eta_2 \leq \eta_1$ and let $v_i = v_{\eta_i}$, $K_i = K_{\eta_i}$, $i = 1,2$. Taking the difference of (\ref{4.10e}) for $i = 1, 2$, we obtain
$$i \partial_t (v_1 - v_2) = (i/2) (\eta_1- \eta_2) \Delta_{K_1} v_1 - (1/2) (1 - i \eta_2) \Big( ( \Delta_{K_1} - \Delta_{K_2}) v_1$$
$$ + \Delta_{K_2} (v_1 - v_2)\Big ) + V(v_1 - v_2)$$

\noi so that
$$\partial_t \parallel v_1 - v_2\parallel_2^2 \ = (\eta_1 - \eta_2) {\rm Re} <v_1 - v_2, \Delta_{K_1} v_1>$$
$$-\  {\rm Im} <v_1 - v_2, (1 - i\eta_2) (\Delta_{K_1} - \Delta_{K_2}) v_1 > - \eta_2 \parallel \nabla_{K_2} (v_1 - v_2)\parallel_2^2$$

\noi and therefore
\beq
\label{4.56e}
\partial_t \parallel v_1 - v_2\parallel_2\ \leq  \eta_1\parallel \Delta_{K_1} v_1\parallel_2\ + \ \parallel (\Delta_{K_1} - \Delta_{K_2})v_1\parallel_2\ .
\eeq

\noi We then estimate
\bea
\label{4.57e}
\parallel (\Delta_{K_1} - \Delta_{K_2})v_1\parallel_2 &\leq&\parallel \nabla (K_1 - K_2) \parallel_2 \Big ( \parallel \nabla v_1 \parallel_3\ + \ \parallel v_1 \parallel_{\infty} \nn \\ 
&&+ \ \parallel K_1 + K_2 \parallel_6 \ \parallel v_1 \parallel_6 \Big )\ .
\eea

\noi One chooses the regularization of $K$ in such a way that
\beq
\label{4.58e}
\parallel K_{\eta} - K ; \dot{H}^1 \parallel\ \leq \ \eta \parallel K ; \dot{H}^2 \parallel
\eeq

\noi which is ensured for instance by regularizing by convolution with a function of $x/\eta$. From  (\ref{4.56e}) (\ref{4.57e}) (\ref{4.58e}) and from the uniform bounds  (\ref{4.45e}) (\ref{4.46e}) (\ref{4.48e}) it follows that 
$$\partial_t \parallel v_1 - v_2\parallel_2\ \leq C \ \eta_1\ .$$

 \noi Therefore $v_{\eta}$ has a limit $v \in {\cal C} (I_+, L^2)$ in $L^{\infty} (I_+, L^2)$ norm. From that convergence and from the uniform bound (\ref{4.55e}), it follows that $v \in {\cal C} (I_+, H^k) \cap (L^{\infty} \cap {\cal C}_w) (I_+ , H^3)$ for $0 \leq k < 3$ and that $v_{\eta}$ converges to $v$ in $L^{\infty} (I_+, H^k)$ norm, in the weak $*$ sense in $L^{\infty}(I_+, H^3)$ and weakly in $H^3$ pointwise in time. From (\ref{4.10e}) for $v_{\eta}$, (\ref{4.58e}) and the previous convergence, it follows that $\partial_t v_{\eta}$ converges to $\partial_t v$ in $L^{\infty} (I_+, L^2)$ norm. From that convergence and from the uniform bound (\ref{4.54e}), it follows that $\partial_t v \in {\cal C} (I_+, H^k) \cap (L^{\infty} \cap {\cal C}_w) (I_+ , H^1)$ for $0 \leq k < 1$ and that $\partial_t v_{\eta}$ converges to $\partial_t v$ in $L^{\infty} (I_+, H^k)$ norm, in the weak $*$ sense in $L^{\infty}(I_+, H^1)$ and weakly in $H^1$ pointwise in time. Clearly $v$ satisfies (\ref{4.1e}). Furthermore $v$ satisfies the conservations laws (\ref{4.2e}) (\ref{4.3e}) which can be proved directly from (\ref{4.1e}) under the available regularity properties of $v$. The conservation law (\ref{4.4e}) is more delicate and can be proved only as an inequality at the present stage. We start from the conservation law for $v_{\eta}$, namely
\bea
\label{4.59e}
&&\parallel   \nabla_{K \eta} \partial_t v_{\eta}(t)   \parallel_2^2 \ - \ \parallel   \nabla_{K \eta} (H_{0\eta}v_0 + \widetilde{f}_{\eta} (t_0))   \parallel_2^2 \ = - \eta \int_{t_0}^t dt' \parallel   \Delta_{K \eta} \partial_t v_{\eta}(t')   \parallel_2^2 \nn \\ 
&&+ \int_{t_0}^t dt' \ 2\ {\rm Im} <\nabla_{K \eta} \partial_t v_{\eta}, f_{3 \eta}> (t')\ .
\eea

\noi When $\eta \to 0$, $\nabla_{K \eta} \partial_t v_{\eta}(t)$ converges to $\nabla_{K} \partial_t v(t)$ weakly in $L^2$, so that 
$$\parallel   \nabla_{K} \partial_t v(t)   \parallel_2^2\ \leq \lim {\rm inf} \parallel   \nabla_{K \eta} \partial_t v_{\eta}(t)   \parallel_2^2\ .$$

\noi The second term in the left hand side of (\ref{4.59e}) converges to the obvious limit by estimates similar to (\ref{4.27e}) (\ref{4.29e}). The first term in the right hand side is negative. In the second term, $\nabla_{K \eta} \partial_t v_{\eta}$ is bounded in $L^{\infty} (I_+, L^2)$ and converges to $\nabla_K \partial_t v$ weakly in $L^2$ pointwise in time. On the other hand $f_{3\eta}$ is also bounded in $L^{\infty} (I_+, L^2)$ and converges to $f_3$ strongly in $L^2$ pointwise in time by  (\ref{4.49e}) (\ref{4.28e}) (\ref{4.30e}) (\ref{4.50e}). Therefore that term converges to the obvious limit by the Lebesgue dominated convergence theorem. Thus we obtain for all $t \in I_+$
\beq
\label{4.60e}
\parallel   \nabla_{K} \partial_t v(t)   \parallel_2^2 \ - \ \parallel   \nabla_{K} (Hv_0 + f_0 (t_0))   \parallel_2^2 \ \leq  \int_{t_0}^t dt' \ 2\ {\rm Im} <\nabla_{K} \partial_t v, f_{3}> (t')\ .
\eeq

At this stage we have proved Part (2) of the proposition in $I_+$, except for the strong continuity of $v$ and $\partial_t v$ at the upper levels of $\dot{H}^3$ and $\dot{H}^1$, and for the conservation law (\ref{4.4e}) for which we have only the one sided inequality (\ref{4.60e}).\par

We next consider the interval $I_- = (0, t_0]$. By the same method we obtain a solution $v \in {\cal C} (I_-, H^k) \cap (L_{loc}^{\infty} \cap {\cal C}_w)(I_-, H^3)$ for $0 \leq k < 3$ with $\partial_t v \in  {\cal C} (I_-, H^k) \cap (L_{loc}^{\infty} \cap {\cal C}_w)(I_-, H^1)$ for $0 \leq k < 1$. Let now $0 <t_1 <t_2 \leq \tau$. We apply the first construction with $t_0 = t_1$ in $[t_1, \tau ]$ and the second one with $t_0 = t_2$ in $(0, t_2]$. Using (\ref{4.60e}) with $t_0 = t_1$ and its reverse analogue with $t_0 = t_2$, and using the uniqueness of Part (1), we obtain (\ref{4.4e}) for any $t_1, t_2 \in I$. Finally (\ref{4.4e}) together with weak continuity implies strong continuity of $\partial_t v$ in $\dot{H}^1$, from which strong continuity of $v$ in $\dot{H}^3$ follows from the gradient of (\ref{4.1e}). This completes the proof of Part (2).\\

\noi \underline{Part (3).} Let $\psi \in {\cal C}_0^{\infty} ( {I\hskip-1truemm R}^3, {I\hskip-1truemm R}^+)$, $0 \leq \psi \leq 1$, $\psi (x) = 1$ for $|x| \leq 1$, $\psi (x) = 0$ for $|x| \geq 2$ and let $\rho_R = x \psi (x/R)$ for $R > 0$, so that  
\beq
\label{4.61e}
|\nabla \rho_R| \leq c \qquad , \qquad |\Delta \rho_R| \leq c\ R^{-1}
\eeq

\noi for all $x \in {I\hskip-1truemm R}^3$. Let $v$ be the solution of (\ref{4.1e}) obtained in Part (2). Then $\rho_Rv $ satisfies the equation
\beq
\label{4.62e}
i \partial_t \rho_Rv = H \rho_R v + ( \nabla \rho_R )\cdot \nabla_K v + (1/2) (\Delta \rho_R) v + \rho_R f_0\ . 
\eeq

\noi The function $\rho_R v$ satisfies the conclusions of Part (2) and in particular the conservation laws
\beq
\label{4.63e}
\parallel \rho_R v(t_2) \parallel_2^2\ -\ \parallel \rho_R v(t_1) \parallel_2^2\ = \int_{t_1}^{t_2} dt\ 2\  {\rm Im} <\rho_R v, f_{0R} >(t)\ ,
\eeq
\beq
\label{4.64e}
\parallel \rho_R \partial_t v(t_2) \parallel_2^2\ -\ \parallel \rho_R  \partial_t v(t_1) \parallel_2^2\ = \int_{t_1}^{t_2} dt\ 2\  {\rm Im} <\rho_R  \partial_t v, f_{2R} >(t)
\eeq

\noi for all $t_1, t_2 \in I$, with
\bea
\label{4.65e}
&&f_{0R} = (\nabla \rho_R) \cdot \nabla_K v + (1/2) (\Delta \rho_R) v + \rho_R f_0\ ,\\
&&f_{2R} = (\partial_t H) \rho_R v + \partial_t f_{0R}\nn \\
&&= (\partial_t H) \rho_R v + (\nabla \rho_R) \cdot \partial_t \nabla_K v + (1/2) (\Delta \rho_R) \partial_t v + \rho_R \partial_t f_0\ .
\label{4.66e}
\eea

\noi From (\ref{4.63e}) (\ref{4.65e}) (\ref{4.61e}) and the properties of $v$ obtained in Part (2), it follows that $\rho_Rv$ is estimated in $L_{loc}^{\infty}(I, L^2)$ uniformly in $R$, so that $xv \in (L_{loc}^{\infty} \cap {\cal C}_w)(I, L^2)$. Taking the limit $R \to \infty$ in (\ref{4.63e}) yields (\ref{4.8e}) which implies that $xv \in {\cal C}(I, L^2)$. From (\ref{4.64e}) (\ref{4.66e}) (\ref{4.61e}), by the same estimates as in the proof of Part (2), in particular (\ref{4.46e}) (\ref{4.47e}), it follows that $\rho_R \partial_t v$ and $\Delta \rho_R v$ are estimated in $L_{loc}^{\infty}(I,L^2)$ uniformly in $R$ so that $x\partial_t v \in (L_{loc}^{\infty} \cap {\cal C}_w)(I, L^2)$. Taking the limit $R \to \infty$ in (\ref{4.64e}) by using the convergence of $\rho_R v$ to $xv$ in $L^2$ and the uniform bound of $\Delta \rho_R v$, we obtain (\ref{4.9e}), which implies that $x \partial_t v\in {\cal C}(I,L^2)$. From that result and from (\ref{4.1e}) it follows that also $x \Delta v \in {\cal C}(I,L^2)$. This completes the proof of Part (3). \par \nobreak \hfill $\sq$ \par

We now turn to the nonlinear system (\ref{2.41e}) and its linearized version (\ref{2.42e}). We shall eventually derive a uniqueness result for the solutions of (\ref{2.41e}) under suitable assumptions on their behaviour at time zero. For that purpose, we need some estimates of the difference of two solutions of the system (\ref{2.42e}). Those estimates will be used again in the contraction proof of Section 5 below.

We shall use the following notation. Let $f_i$, $i = 1,2$ be two functions or operators associated with a pair $(w_i, s_i, B_{2i})$ or $(w'_i, B'_{2i})$, $i = 1,2$. We define $f_{\pm} = (1/2) (f_1 \pm f_2)$ so that $f_1 = f_+ + f_-$, $f_2 = f_+ - f_-$ and $(fg)_{\pm} = f_+ g_{\pm} + f_- g_{\mp}$. Let now $(w'_i, B'_{2i})$, $i = 1,2$, be a pair of solutions of the system  (\ref{2.42e}) associated with a pair $(w_i, s_i, B_{2i})$, $i = 1,2$. Then $(w'_-, B'_{2-})$ satisfies the equations
\beq
\label{4.67e}
\left \{ \begin{array}{l} i \partial_t w'_- = H_+ w'_- + H_- w'_+\\ \\ B'_{2-} = 2 {\cal B}_2 (w_+, w_-, K_+) - t F_2 (P K_- (|w_+|^2 + |w_-|^2))\end{array} \right .
\eeq

\noi where
\bea
\label{4.68e}
&&H_+ = - (1/2) \Delta_{K_+} + (1/2) K_-^2 + {\check B}_{S+}\ , \\
&&H_- = i K_- \cdot \nabla _{K_+} + (i/2) (\nabla \cdot K_-) + {\check B}_{S-}\ . 
\label{4.69e}
\eea

\noi We shall also need (see (\ref{2.31e}) (\ref{2.32e}))
\beq
\label{4.70e}
B_{1-} = 2B_1(w_+, w_-) \qquad , \quad B_{1S/L-} = 2 B_{1S/L} (w_+, w_-) \ .
\eeq

\noi If $s_i$, $i=1,2$, satisfy the second equation of (\ref{2.41e}), then 
\beq
\label{4.71e}
\partial_t s_- = 2 t^{-1} \nabla g (w_+, w_-) + \nabla {\check B}_{1L-}\ .
\eeq

\noi We can now state the difference estimates of two solutions of the system (\ref{2.42e}). The estimates of $w'_-$ are stated in differential form for brevity, but should be understood in integral form, in the same way as the conservation laws (\ref{4.2e})-(\ref{4.4e}).\\

\noi {\bf Lemma 4.1.} {\it Let $0 < \beta < 1$. Let $0 < \tau \leq 1$, let $I = (0, \tau ]$ and let $h_1 \in {\cal C}(I, {I\hskip-1truemm R}^+)$ satisfy
\beq
\label{4.72e}
\int_0^{\tau} dt \ t^{-3/2} h_1 (t) < \infty \ .
\eeq

\noi Let $w_i$, $i = 1,2$ satisfy $w_i \in L^{\infty} (I, H^3)$, $x w_i \in L^{\infty} (I, H^2)$ and 
\beq
\label{4.73e}
\parallel <x> w_- (t) \parallel_2 \ \leq C \ h_1(t)
\eeq 

\noi for all $t \in I$.\par

(1) Let $B_1 (w_i)$ be defined by (\ref{2.31e}) (\ref{2.32e}). Then $B_1(w_i) \in ({\cal C} \cap L^{\infty} ) (I , \ddot{H}^4)$, $t \nabla {\check B}_1 (w_i) \in ({\cal C} \cap L^{\infty} ) (I , \ddot{H}^2)$ and $B_{1-}$ satisfies the estimates
\bea
\label{4.74e}
&&\parallel \nabla B_{1-} \parallel_2 \ \leq C \ I_0 \left ( \parallel <x> w_-  \parallel_2 \right )\ , \\
&&\parallel \nabla {\check B}_{1-} \parallel_2 \ \leq C \ t^{-1} I_{-1}  \left ( \parallel <x> w_-  \parallel_2 \right )\ .
\label{4.75e}
\eea 

(2) Let $s_i$ satisfy the second equation of (\ref{2.41e}) with $w = w_i$, $i = 1,2$, with $s_i (t_0) \in \ddot{H}^2$ for some $t_0 \in I$. Then $s_i \in {\cal C} (I, \ddot{H}^2), s_- \in {\cal C} (I, H^2)$ and $s_-$ satisfies the estimates
\beq
\label{4.76e}
\parallel \nabla^k \partial_t s_{-} \parallel_2 \ \leq C \ t^{-1-k \beta} I_{-1} \left ( \parallel <x> w_-  \parallel_2 \right ) + C\ t^{-1} \left ( \parallel w_-  \parallel_2\ + \delta_{k2}   \parallel \nabla w_-  \parallel_2\right )
\eeq 
\noi for $k = 0, 1, 2$.\par

(3) Let $B_0$ satisfy (\ref{3.32e}) for $0 \leq j , k , j+k \leq 1$ and $r = \infty$. Let in addition $t \partial_t w_i \in L^{\infty}(I, H^1)$, $i = 1,2$. Let $B_{2i}$ satisfy $B_{2i} \in L^{\infty} (I, \ddot{H}^2)$, $t \partial_t B_{2i} \in L^{\infty}(I, \dot{H}^1)$, $t \nabla  {\check B}_{2i}  \in L^{\infty}(I, \dot{H}^1)$, $i = 1,2$. Let $(w'_i, B'_{2i})$, $i = 1,2$ be solutions of the linearized system (\ref{2.42e}) satisfying the same conditions as $(w_i, B_{2i})$. Then the following estimates hold~:
\bea
\label{4.77e}
&&\left | \partial_t \parallel w'_-  \parallel_2\right |\ \leq\ C \Big ( \parallel \nabla \cdot s_-  \parallel_2\  + (1 - \ell n\ t) \left ( \parallel s_-  \parallel_3\ + \ \parallel \nabla B_-  \parallel_2\right )\nn\\
&&\qquad\qquad +\  t^{\beta} \parallel \nabla   {\check B}_{1-} \parallel_2\ + \ \parallel \nabla  {\check B}_{2-}  \parallel_2\ \Big ) ,
\eea
\bea
\label{4.78e}
&&\left | \partial_t \parallel xw'_-  \parallel_2\right | \ \leq\ \parallel \nabla_{K_+} w'_-  \parallel_2\  + \ C \Big ( \parallel \nabla \cdot s_-  \parallel_2\ + (1 - \ell n\ t) \nn \\
&&\left  ( \parallel s_-  \parallel_3\ + \ \parallel \nabla B_-  \parallel_2\right ) \  +\  t^{\beta} \parallel \nabla   {\check B}_{1-} \parallel_2\ + \ \parallel \nabla  {\check B}_{2-}  \parallel_2\ \Big ) ,
\eea
\bea
\label{4.79e}
&&\left | \partial_t \parallel \nabla_{K_+} w'_-  \parallel_2\right |\ \leq\ C\Big ( t^{-1} \left (  \parallel w'_-  \parallel_2\ + \ \parallel w'_-  \parallel_3\right )  + \ t^{-1}  \parallel s_-  \parallel_2\nn \\
&& + \ (1 - \ell n\ t) \parallel \nabla s_-  \parallel_2\ +\ \parallel \nabla \nabla \cdot s_-\parallel_2\  + \  t^{-1} \parallel \nabla  B_-  \parallel_2\ + \  \parallel \nabla {\check B}_{1-} \parallel_2\nn\\
&& + \ (1 - \ell n\ t) \parallel \nabla  {\check B}_{2-}  \parallel_2\ \Big ) ,
\eea    
\beq
\label{4.80e}
\parallel \nabla  B'_{2-}  \parallel_2\ \leq\ C\ t \ I_1 \Big (  (1 - \ell n\ t) \parallel w_-  \parallel_2\ + \ \parallel s_-  \parallel_2\ + \ \parallel \nabla B_-  \parallel_2\Big )\ , 
\eeq
\beq
\label{4.81e}
\parallel \nabla  {\check B}'_{2-}  \parallel_2\ \leq\ C \ I_0 \Big (  (1 - \ell n\ t) \parallel <x>w_-  \parallel_2\ + \ \parallel s_-  \parallel_2\ + \ \parallel \nabla B_-  \parallel_2\Big )
\eeq  

\noi for all $t \in I$.}\\

\noi {\bf Proof.} We first estimate $B_1(w_i)$ and $s_i$. From (\ref{2.31e}) (\ref{2.32e}), from (\ref{3.18e}) and from Lemma 3.2, we estimate
\beq
\label{4.82e}
\parallel \omega^{m+1} B_1(w_i)\parallel_2\ \leq\ I_m \left ( \parallel  \omega^m x |w_i|^2\parallel _2\right ) \leq C
\eeq 

\noi for $0 \leq m \leq 3$. Similarly from (\ref{3.19e}), we estimate
\beq
\label{4.83e}
\parallel \omega^{m+1} {\check B}_1(w_i)\parallel_2\ \leq\ t^{-1} I_{m-1}  \left ( \parallel  <x> \omega^m x |w_i|^2\parallel _2\right ) \leq C\ t^{-1}
\eeq 

\noi for $1 \leq m \leq 2$. Similarly from (\ref{3.18e}), we estimate
\beq
\label{4.84e}
\parallel \omega^{m+1} \partial_t {B}_1(w_i)\parallel_2\ \leq\ I_{m+1}  \left ( \parallel \omega^m x\partial_t  |w_i|^2\parallel _2\right ) \leq C\ t^{-1}
\eeq 

\noi for $0 \leq m \leq 1$. From (\ref{4.83e}) and Lemma 3.5, we obtain
\beq
\label{4.85e}
\parallel \omega^{m+1} \partial_t s_i\parallel_2\ \leq\ C\ t^{-1}
\eeq 

\noi for $0 \leq m \leq 1$, and by integration over time
\beq
\label{4.86e}
\parallel \omega^{m+1} s_i\parallel_2\ \leq\ C(1 - \ell n \ t)
\eeq 

\noi for $0 \leq m \leq 1$ and for all $t\in I$. Lower norms of $B_1(w_i)$, ${\check B}_1(w_i)$ and $s_i$ can also be estimated. In particular one can show that $B_1(w_i)$, $\nabla {\check B}_1(w_i)$ and $s_i$ are estimated in $L^4$. The argument will be given on the example of $B_1(w_a)$, ${\check B}_1(w_a)$ and $s_a$ at the beginning of the proof of Lemma 6.1 below. From (\ref{4.82e})-(\ref{4.86e}) it follows in particular that
\beq
\label{4.87e}
\parallel B_{1+} \parallel_{\infty}\ \vee\ \parallel \nabla B_{1+} \parallel_{\infty} \ \vee\ t \parallel  \nabla  {\check B}_{1+}\parallel_{\infty}\ \vee\ t \parallel \partial_t B_{1+} \parallel _{\infty}\ \leq C \ ,
\eeq
\beq
\label{4.88e}
t \parallel\partial_t s_+ \parallel_{\infty}\ \leq C \qquad , \quad \parallel s_+ \parallel_{\infty} \ \vee\  \parallel  \nabla  s_+\parallel_6\ \leq C(1 - \ell n\ t )\ .\eeq  

\noi We now begin the proof of the lemma proper.\\

\noi \underline{Part (1).} The properties of $B_1(w_i)$ follow from (\ref{4.82e}) (\ref{4.83e}). The estimates (\ref{4.74e}) (\ref{4.75e}) follow immediately from (\ref{4.70e}) (\ref{3.18e}) (\ref{3.19e}). Note that at this point that the condition (\ref{4.72e}) ensures the convergence of the integral in the right hand side of (\ref{4.75e}).\\

\noi \underline{Part (2).} The properties of $s_i$ follow from (\ref{4.86e}). The estimate (\ref{4.76e}) follows from (\ref{4.75e}) (\ref{3.12e}) and Lemma 3.5.\\

\noi \underline{Part (3).} We first estimate $\parallel w'_-\parallel_2$. From (\ref{4.67e}), by using (\ref{3.11e}), we obtain 
$$|\partial_t \parallel w'_-\parallel_2| \ \leq \ \parallel H_-w'_+\parallel_2 \ \leq\ \parallel s_-\parallel_3 \ \parallel \nabla_{K_+} w'_+\parallel_6 \ + \ \parallel B_-\parallel_6 \ \parallel \nabla_{K_+} w'_+\parallel_3$$
$$+\ \parallel \nabla \cdot s_-\parallel_2\ \parallel w'_+\parallel_{\infty}\ + \ t^{\beta} \parallel\nabla {\check B}_{1-}\parallel_2 \ \parallel w'_+\parallel_{\infty} \ + \ \parallel {\check B}_{2-}\parallel_6 \ \parallel w'_+\parallel_3$$

\noi from which (\ref{4.77e}) follows by using the assumptions on $w_i$, $B_0$, $B_{2i}$ and the estimate (\ref{4.87e}). \par

We next estimate $\parallel xw'_-\parallel_2$. From (\ref{4.67e}) (\ref{3.11e}) and the commutation relation 
$$[x, H_+] = \nabla_{K_+}$$

\noi we obtain similarly
$$|\partial_t \parallel x w'_-\parallel_2| \ \leq \ \parallel \nabla_{K_+}w'_-\parallel_2 \ + \ \parallel x H_- w'_+\parallel_2$$
$$\leq  \ \parallel \nabla_{K_+} w'_-\parallel_2 \ + \ \parallel s_-\parallel_3 \ \parallel x \nabla_{K_+} w'_+\parallel_6\ + \ \parallel B_-\parallel_6 \ \parallel x \nabla_{K_+} w'_+\parallel_3$$
$$+\ \parallel \nabla \cdot s_-\parallel_2\ \parallel xw'_+\parallel_{\infty}\ + \ t^{\beta} \parallel\nabla {\check B}_{1-}\parallel_2 \\ \parallel xw'_+\parallel_{\infty}\ +\ \parallel {\check B}_{2-}\parallel_6 \ \parallel xw'_+\parallel_3$$

\noi from which (\ref{4.78e}) follows in the same way as before. \par

We next estimate $\parallel \nabla_{K_+} w'_-\parallel_2$. Taking the covariant gradient in (\ref{4.67e}) yields 
$$i \partial_t \nabla_{K_+} w'_- = - (1/2) \nabla_{K_+} \Delta_{K_+} w'_- + \left ( (1/2) K_-^2 + {\check B}_{S+}\right ) \nabla_{K_+} w'_-$$
$$+ \left ( \partial_t K_+ + K_- \nabla K_- + \nabla  {\check B}_{S+}\right ) w'_- + i K_- \cdot \nabla_{K_+}^2 w'_+ + i (\nabla K_-) \cdot \nabla_{K_+} w'_+$$
$$+ (i/2) (\nabla \cdot s_-) \nabla_{K_+} w'_+ + (i/2) (\nabla \nabla \cdot s_-) w'_+ +  {\check B}_{S-} \nabla_{K_+} w'_+ + (\nabla  {\check B}_{S-}) w'_+$$

\noi from which we estimate
\bea
\label{4.89e}
&&|\partial_t \parallel \nabla_{K_+} w'_-\parallel_2| \ \leq \ \parallel (\partial_t K_+ + \nabla  {\check B}_{S+})w'_-\parallel_2 \ + \ \parallel K_-\cdot  \nabla_{K_+}^2 w'_+\parallel_2\nn \\
 &&+\  \parallel \nabla K_- \parallel_2 \left ( \parallel \nabla_{K_+} w'_+\parallel_{\infty} \ + \  \parallel K_- w'_-\parallel_{\infty}\right ) + \ \parallel \nabla \nabla \cdot s_-\parallel_2\ \parallel w'_+\parallel_{\infty}\nn \\
&& + \ \parallel\nabla {\check B}_{1-}\parallel_2 \left ( t^{\beta}  \parallel \nabla_{K_+}  w'_+\parallel_{\infty} \ + \ \parallel w'_+\parallel_{\infty} \right ) \ \nn \\
&&+ \ \parallel \nabla  {\check B}_{2-}\parallel_2 \left ( \parallel \nabla_{K_+} w'_+\parallel_3\ + \ \parallel w'_+\parallel_{\infty} \right ) \ .
  \eea
  
\noi We next estimate the first two terms in the right hand side of (\ref{4.89e}).
\bea
\label{4.90e}
&&\parallel (\partial_t K_+ + \nabla {\check B}_{S+}) w'_-\parallel_2\ \leq \  \parallel \partial_t (s_+ + B_0 + B_{1+}) + \nabla ({\check B}_0 + {\check B}_{1+} ) \parallel_{\infty} \ \parallel w'_- \parallel_2 \nn \\
&&+\ \parallel \partial_t  B_{2+} + \nabla {\check B}_{2+} \parallel_6\ \parallel w'_-\parallel_3\ \leq \ C\ t^{-1} \left ( \parallel w'_-\parallel_2\ + \ \parallel w'_-\parallel_3 \right ) \ , 
\eea 
$$\parallel K_- \nabla_{K_+}^2 w'_+ \parallel_2\ \leq \ \parallel s_- \parallel_3 \Big ( \parallel \nabla^2 w'_+ \parallel_6 \ + \ \parallel \nabla (s_+ + B_{2+})\parallel_6\ \parallel w'_+ \parallel_{\infty} \Big )$$
$$+ \ \parallel s_- \parallel_2 \Big ( \parallel K_+ \parallel_{\infty} \ \parallel \nabla w'_+ \parallel_{\infty}\ + \left ( \parallel \nabla (B_0 + B_{1+}) \parallel_{\infty} \ + \ \parallel K_+ \parallel_{\infty}^2 \right ) \parallel w'_+ \parallel_{\infty} \Big )$$
$$+ \ \parallel B_- \parallel_6 \Big ( \parallel  \nabla^2 w'_+\parallel_3\ +\ \parallel  K_+ \parallel_{\infty}\ \parallel \nabla w'_+ \parallel_3\ + \ \parallel \nabla (B_0 + B_{1+}) \parallel_{\infty} \ \parallel w'_+ \parallel_3$$
$$+\ \parallel \nabla (s_+ + B_{2+} ) \parallel_6\ \parallel w'_+ \parallel_6\ + \ \parallel K_+\parallel_{\infty}^2 \ \parallel w'_+ \parallel_3 \Big )$$
\beq
\label{4.91e}
\leq C \left ( (1 - \ell n\ t) \parallel s_- \parallel_3 \ + \ t^{-1} \left ( \parallel s_- \parallel_2 \ + \ \parallel \nabla B_- \parallel_2 \right ) \right ) \ .
\eeq 

\noi Substituting  (\ref{4.90e})  (\ref{4.91e}) into  (\ref{4.89e}) and estimating the remaining terms in a similar way yields  (\ref{4.79e}).\par

We finally estimate $B'_{2-}$. From  (\ref{4.67e}) and (\ref{3.18e}) (\ref{3.19e}), we obtain 
$$\parallel  \nabla B'_{2-} \parallel_2\ \leq t\ I_1 \left ( \parallel  w_- \parallel _2\ \parallel  \nabla_{K_+} w_+ \parallel _{\infty} \ + \ \parallel  s_- \parallel _2 \ \parallel  w_+ \parallel _{\infty}^2\ + \ \parallel B_- \parallel _6 \ \parallel w_+ \parallel_6^2 \right ) \ ,$$
$$\parallel  \nabla {\check B}'_{2-} \parallel_2\ \leq \ I_0 \Big ( \parallel  <x> w_- \parallel _2\ \parallel  \nabla_{K_+} w_+ \parallel _{\infty} \ + \ \parallel  s_- \parallel _2 \ \parallel  w_+ \parallel _{\infty}\ \parallel  <x> w_+ \parallel _{\infty}$$
$$+ \ \parallel B_- \parallel _6 \ \parallel w_+ \parallel_6\ \parallel  <x> w_+ \parallel _6 \Big ) \ ,$$

\noi from which (\ref{4.80e}) and (\ref{4.81e}) follow.\par \nobreak \hfill $\sq$ \par

As a first application of Lemma 4.1, we give a uniqueness result for the nonlinear system (\ref{2.41e}) with initial condition at time zero. That result is a variant of Proposition 4.2, part (2) of I.\\

\noi {\bf Proposition 4.2.} {\it Let $0 < \beta < 1$. Let $0 < \tau \leq 1$, let $I = (0, \tau ]$ and let $h_1 \in {\cal C} (I, {I\hskip-1truemm R}^+)$ be such that $\overline{h}_1 (t) = (t^{-2\beta} \vee t^{-1/2}) h_1(t)$ be non decreasing and satisfy  
\beq
\label{4.92e}
\int_0^t dt'\ t{'}^{-1} \overline{h}_1 (t') \leq c\ \overline{h}_1(t)
\eeq

\noi for some $c > 0$ and for all $t \in I$. Let $B_0$ satisfy (\ref{3.32e}) for $r = \infty$ and $0 \leq j, k, j+k \leq 1$. Let $(w_i, s_i, B_{2i})$, $i = 1,2$, be two solutions of the system (\ref{2.41e}) such that $w_i \in L^{\infty} (I, H^3)$, $xw_i \in L^{\infty} (I, H^2)$, $t \partial_t w_i \in L^{\infty} (I, H^1)$, $B_{2i} \in L^{\infty} (I, \ddot{H}^2)$, $t \partial_t B_{2i} \in L^{\infty} (I, \dot{H}^1)$ and $t\nabla {\check B}_{2i} \in L^{\infty} (I, \dot{H}^1)$. Assume in addition that $s_- (0) = 0$ and that 
\beq
\label{4.93e}
\parallel <x> w_- (t) \parallel _2\ \leq C\ h_1(t)
\eeq

\noi for all $t \in I$. Then $(w_1, s_1, B_{21}) = (w_2, s_2, B_{22})$.}\\

\noi {\bf Proof.} Note first that (\ref{4.92e}) implies (\ref{4.72e}) so that Lemma 4.1 can be applied. From (\ref{4.76e}) with $k = 0$ and mild assumptions on $w_-$, it follows that $s_-(t)$ has an $L^2$ limit as $t \to 0$, thereby giving a meaning on the assumption $s_- (0) = 0$. Actually it follows from (\ref{4.76e}) (\ref{4.92e}) (\ref{4.93e}) that the limit exists in $H^2$. \par

We first prove the proposition for $\tau$ sufficiently small by using Lemma 4.1 with $(w'_i, B'_{2i}) = (w_i, B_{2i})$. We define
$$y_0 = \ \parallel <x> w_- \parallel_2 \qquad , \quad y_1 = \ \parallel \nabla_{K_+} w_- \parallel_2\ ,$$
$$Y_0 = \ \mathrel{\mathop {\rm Sup}_{t \in I}}\ h_1 (t)^{-1} y_0 (t)\ .$$

\noi From Lemma 4.1, especially  (\ref{4.74e})  (\ref{4.75e})  (\ref{4.76e}) and from  (\ref{4.92e}), we obtain
\beq
\label{4.94e}
\parallel \nabla B_{1-}\parallel_2\  \leq C\ Y_0\ h_1\ ,
\eeq
\beq
\label{4.95e}
\parallel \nabla {\check B}_{1-}\parallel_2\  \leq C\ Y_0\ t^{-1} h_1\ ,
\eeq
\beq
\label{4.96e}
\parallel \nabla^k  s_- \parallel_2\  \leq C\ Y_0\ t^{-k\beta } h_1 + C\ \delta_{k2} \int_0^t dt'\ t{'}^{-1} y_1(t')
\eeq

\noi for all $t \in I$ and for $k = 0, 1, 2$. The time integral of the last term in (\ref{4.76e}) converges because of the estimate
$$\parallel \nabla  w_- \parallel_2\  \leq\ \left ( \parallel w_- \parallel_2\ \parallel \Delta  w_- \parallel_2\right )^{1/2} \ \leq C(Y_0\ h_1)^{1/2}$$

\noi and we have replaced the ordinary derivative by the covariant derivative in that integral, which produces an innocuous term with $Y_0h_1$. On the other hand from (\ref{4.80e})  (\ref{4.94e})  (\ref{4.96e}) we obtain
\beq
\label{4.97e}
\parallel \nabla  B_{2-} \parallel_2\  \leq \ C\ Y_0 \ t(1 - \ell n\ t) h_1 + C\ t\ I_1 \left ( \parallel \nabla B_{2-} \parallel_2 \right ) \ .
\eeq

\noi From the assumptions on $B_{2i}$, it follows that $B_{2-} \in L^{\infty} (I, \dot{H}^1)$. Using that fact, one derives easily from (\ref{4.97e}) that
\beq
\label{4.98e}
\parallel \nabla  B_{2-} \parallel_2\  \leq \ C\ Y_0 \ t(1 - \ell n\ t) h_1 
\eeq

\noi for all $t \in I$ and for $\tau$ sufficiently small. Substituting that result into  (\ref{4.81e}) yields
\beq
\label{4.99e}
\parallel \nabla  {\check B}_{2-} \parallel_2\  \leq \ C\ Y_0 (1 - \ell n\ t) h_1 \ . 
\eeq

\noi Substituting (\ref{4.94e}) (\ref{4.95e}) (\ref{4.96e}) (\ref{4.98e}) (\ref{4.99e}) into (\ref{4.79e}) yields
$$|\partial_t y_1| \leq C \left ( Y_0 (t^{-1} + t^{-2 \beta} ) h_1 + \int_0^t dt'\ t{'}^{-1} y_1(t') + t^{-1} (Y_0 h_1 y_1)^{1/2}\right )$$

\noi which takes the form 
\beq
\label{4.100e}
|\partial_t y| \leq f + C \int_0^t dt'\ t{'}^{-1} y(t') + g y^{1/2}
\eeq

\noi with
$$y_1 = Y_0 \ y\quad , \quad f = C \left ( t^{-1} + t^{-2\beta}\right ) h_1 \quad , \quad g = C\ t^{-1} \ h_1^{1/2}\ .$$

\noi We define
$$z(t) = \int_0^t dt'\ t{'}^{-1} y(t')$$

\noi so that $y = t \partial_t z$, and
$$F(t) = \int_0^t dt' \ f(t')\ .$$

\noi Integrating (\ref{4.100e}) over time yields
$$y(t) \leq F(t) + C \int_0^t dt'\ t{'}^{-1}(t-t') y(t') + \int_0^t dt' \left ( g y^{1/2}\right ) (t')$$
$$\leq F(t) + C \ t\ z(t) + z(t)^{1/2} \left ( \int_0^t dt'\ t' g^2 (t') \right )^{1/2}$$

\noi so that
\beq
\label{4.101e}
\partial_t z \leq t^{-1} F + C z + g z^{1/2}
\eeq

\noi where we have used the estimate
$$\int_0^t dt'\ t' g^2 (t') = C \int_0^t dt'\ t{'}^{-1} h_1 (t') \leq C\ h_1 (t) \ .$$

\noi From (\ref{4.101e}) and Lemma 2.3 in \cite{9r}, we obtain 
\begin{eqnarray*}
z &\leq& e^{Ct} \left \{ \int_0^t dt'\ g(t') + \left ( \int_0^t dt'\ t{'}^{-1} F(t')\right )^{1/2} \right \}^2\\
&\leq& C \left ( 1 + t^{1-2\beta } \right ) h_1 (t)
\end{eqnarray*}

\noi by an elementary computation. Substituting that result into (\ref{4.101e}) yields 
\beq
\label{4.102e}
y_1 = C\ Y_0 \left ( 1 + t^{1-2\beta}\right ) h_1(t)\ .
\eeq

\noi Substituting (\ref{4.94e}) (\ref{4.95e}) (\ref{4.96e}) (\ref{4.98e}) (\ref{4.99e}) and (\ref{4.102e}) into (\ref{4.77e}) (\ref{4.78e}) yields
$$\partial_t y_0 \leq C\ Y_0 \left ( t^{-\beta} + t^{-1+\beta}\right ) h_1(t)$$

\noi and therefore by integration over time
$$Y_0 \leq C\ Y_0 \left ( \tau^{1-\beta} + \tau^{\beta} \right )$$

\noi which implies that $Y_0 = 0$ for $\tau$ sufficiently small. Together with (\ref{4.98e}), this concludes the proof in that case. The extension of the proof to the case of general $\tau$ proceeds by similar but more standard arguments.\par \nobreak \hfill $\sq$ \par

\mysection{Cauchy problem at time zero for the auxiliary system}
\hspace*{\parindent} In this section, we derive the main technical result of this paper, namely we construct solutions of the auxiliary system (\ref{2.41e}) with given asymptotic behaviour at time zero parametrized by asymptotic functions $(w_a, s_a, B_a)$ under assumptions of a general nature on those functions. This will be done by solving the auxiliary system (\ref{2.54e}) for the difference variables $(q, G_2)$ with $(q, G_2)$ tending to zero at time zero, with $G_1$ and $\sigma$ defined by (\ref{2.52e}) (\ref{2.53e}) with $\sigma (0) = 0$. The latter system will be solved in two steps. We first solve the linearized system (\ref{2.59e}) for $(q', G'_2)$ with given $(q, G_2)$, thereby defining a map $\Gamma : (q, G_2) \to (q', G'_2)$. We then prove that the map $\Gamma$ is a contraction in $X(I)$ for $I = (0, \tau ]$ and $\tau$ sufficiently small. The Schr\"odinger equation of the linearized system (\ref{2.59e}) with $q'$ tending to zero at time zero is solved by first constructing a solution $q'_{t_0}$ with initial condition $q'_{t_0}(t_0) = 0$ for some small $t_0 > 0$ by the use of Proposition 4.1 and then taking the limit of that solution when $t_0$ tends to zero. The system (\ref{2.59e})  yields $G'_2$ as an explicit function of $(q, G_2)$. The general assumptions made on $(w_a, s_a, B_a)$ consist of boundedness properties of $w_a$, boundedness properties of $(s_a, B_a)$ and decay properties of the remainders $R_j$, $1 \leq j \leq 4$, defined by (\ref{2.57e}), which we state below as assumptions (A1) (A2) and (A3) respectively. Those assumptions are stated in terms of an interval $I_0 = (0, \tau_0 ]$ with $0 < \tau_0 \leq 1$.\\

\noi (A1) $w_a$ satisfies the following properties
\beq
\label{5.1e}
w_a \in \left ( {\cal C} \cap L^{\infty}\right ) \left (I_0, H^3 \right ) \quad , \qquad x w_a \in \left ( {\cal C} \cap L^{\infty}\right ) \left (I_0, H^2 \right ) \ ,
\eeq
\beq
\label{5.2e}
t^{1/2} \partial_t w_a \in \left ( {\cal C} \cap L^{\infty}\right ) \left (I_0, H^2 \right )\quad , \quad  t^{1/2} x\partial_t w_a \in \left ( {\cal C} \cap L^{\infty}\right ) \left (I_0, H^1 \right )\ .
\eeq
\vskip 5 truemm

In order to state (A2), we recall that $B_a = B_0 + B_{1a} + B_{2a}$ and $K_a = s_a + B_a$, and that ${\check B}_{aS} = {\check B}_{0}+ {\check B}_{1aS} + {\check B}_{2a}$ (see (\ref{2.50e})).\\

\noi (A2) $s_a, B_a \in {\cal C}(I_0 , H_{\infty}^1)$ with sufficient additional regularity, and the following estimates hold for all $t\in I_0$~:
\beq
\label{5.3e}
\parallel K_a \parallel_{\infty} \ \leq C(1 - \ell n\ t)\ ,
\eeq
\beq
\label{5.4e}
\parallel \partial_t K_a \parallel_{\infty} \ \vee\ \parallel \nabla K_a \parallel_{\infty} \ \vee \ t \parallel \nabla \partial_t K_a \parallel_{\infty} \ \leq C \ t^{-1}\ ,\eeq
$$\parallel\nabla  s_a \parallel_{\infty} \ \vee\ \parallel \nabla \nabla \cdot s_a \parallel_{3} \ \vee \  \parallel \nabla (B_{1a} + B_{2a}) \parallel_{\infty} $$\beq
\label{5.5e}
\vee  \ t \Big ( \parallel \nabla \partial_t s_a \parallel_{\infty}\ \vee \parallel \nabla \partial_t \nabla \cdot s_a \parallel_{3} \ \vee\  \parallel \nabla \partial_t (B_{1a} + B_{2a}) \parallel_{\infty} \Big ) \leq C \ t^{-1/2}\ , \eeq
\beq
\label{5.6e}
\parallel \nabla {\check B}_a\parallel_{\infty}\ \vee \ t \parallel\nabla \partial_t {\check B}_a\parallel_{\infty}\ \leq C\ t^{-1}\ ,
\eeq 
\beq
\label{5.7e}
\parallel  {\check B}_{aS}\parallel_{\infty}\ \vee \ t \parallel  \partial_t {\check B}_{aS}\parallel_{\infty}\ \leq C\ t^{-1/2}\ .
\eeq 

Note that by Lemma 3.7, $B_0$ satisfies the assumptions made on $B_a$ under suitable assumptions on $(A_+, \dot{A}_+)$.\\

\noi (A3) The remainders $R_j$ satisfy the following estimates for all $t \in I_0$~:
\beq
\label{5.8e}
\left ( \parallel <x>R_1 \parallel _2\ \leq \right ) \ \parallel <x> \partial_t R_1; L^1((0, t], L^2 ) \parallel \ \leq r_1\ t^{-1}\ h(t)\ ,
\eeq
\beq
\label{5.9e}
\left ( \parallel \nabla R_1 \parallel _2\ \leq \right ) \ \parallel \nabla \partial_t R_1; L^1((0, t], L^2 ) \parallel \ \leq r_1\ t^{-3/2}\ h(t)\ ,
\eeq
\beq
\label{5.10e}
 \parallel \nabla^k R_2  \parallel_2\ \leq r_2\ t^{-1-k\beta}\ h(t)\qquad \hbox{for $k = 0, 1, 2$}\ ,
\eeq 
\bea
\label{5.11e}
&& \parallel \nabla R_3  \parallel_2\ \vee\ t^{1/2}  \parallel \nabla^2 R_3  \parallel_2\ \vee \  t\left ( \parallel \nabla \partial_t  R_3  \parallel_2\ \vee\ \parallel \nabla {\check R}_3  \parallel_2\right )\nn \\
&&\vee\ t^{3/2}  \parallel \nabla^2 {\check R}_3  \parallel_2\ \vee\ t^2  \parallel \nabla \partial_t {\check R}_3  \parallel_2\ \leq r_3 \ h(t)\ ,
\eea
\bea
\label{5.12e}
&& \parallel \nabla R_4  \parallel_2\ \vee\ t\left (  \parallel \nabla^2 R_4  \parallel_2\ \vee \  \parallel \nabla \partial_t  R_4  \parallel_2\ \vee\ \parallel \nabla {\check R}_4  \parallel_2\right ) \nn \\
&&\vee\ t^{2}  \left ( \parallel \nabla^2 {\check R}_4  \parallel_2\ \vee\  \parallel \nabla \partial_t {\check R}_4  \parallel_2\right ) \ \leq r_4\ t^{1/2}\  \ h(t)
\eea

\noi for some positive constants $r_j$, $1 \leq j \leq 4$, where $h$ is defined in Section 3.\\

The final result will require the full assumptions (A1) (A2) (A3),  but some intermediate ones will need only part of them. The final existence result for $(q, G_2)$ solution of (\ref{2.54e}) will be derived in an interval $(0, \tau ]$ for sufficiently small $\tau$. Smallness conditions on $\tau$ will occur at various stages of the proof and will be called asymptotic region conditions. They will eventually depend on $(w_a, s_a, B_a)$. Some of them will be imposed at early stages in order to eliminate terms higher than linear in the dynamical variables from the estimates. Others will be imposed at the final stage in order to ensure the contraction properties of the map $\Gamma$. \par

We first prepare the ground by deriving preliminary estimates of $G_1$, $\sigma$, $G'_2$, $H_1$ and $q'$. The assumptions will in general include the assumption (A1), parts of the assumption (A2) and some asymptotic region conditions. Whenever $(q, G_2)$ occurs, it will be assumed to be in $X_0(\cdot )$. A number of estimates will involve expressions of the type $I_j (\parallel q \parallel_*)$ for some norm of $q$ and we shall assume that norm to decay sufficiently fast to make the integral convergent. Such conditions will always be satisfied for $(q, 0) \in X(\cdot )$.\par

We first estimate $G_1$ and $\sigma$ defined by (\ref{2.52e}) (\ref{2.53e}).\\

\noi {\bf Lemma 5.1.} {\it Let $0 < \beta < 1$, $0 < \tau \leq \tau_0$ and $I = (0, \tau ]$. Let $w_a$ satisfy (A1) and let $(q, 0) \in X_0 (I)$ with 
\beq
\label{5.13e}
\left \{ \begin{array}{l}\parallel q ; L^{\infty} (I, H^3)\parallel \ \leq \ \parallel w_a ; L^{\infty} (I, H^3 ) \parallel \\ \\  \parallel xq ;  L^{\infty} (I, H^2)\parallel \ \leq \ \parallel xw_a ; L^{\infty} (I, H^2 )  \parallel \ . \end{array}\right .
\eeq 

\noi Then the following estimates hold for all $t \in I$~:
\beq
\label{5.14e}
\parallel \nabla G_1  \parallel_2\ \leq \ C\ I_0\left  ( \parallel q  \parallel_2\right ) \ + \ \parallel \nabla R_3  \parallel_2\ ,
\eeq
\beq
\label{5.15e}
\parallel \nabla {\check G}_1  \parallel_2\ \leq \ C\ t^{-1}\ I_{-1}\left  ( \parallel <x>q  \parallel_2\right ) \ + \ \parallel \nabla {\check R}_3  \parallel_2\ ,
\eeq
\beq
\label{5.16e}
\parallel \nabla^2 G_1  \parallel_2\ \leq \ C\ I_1\left  ( \parallel \nabla q  \parallel_2\right ) \ + \ \parallel \nabla^2 R_3  \parallel_2\ ,
\eeq
\beq
\label{5.17e}
\parallel \nabla^2 {\check G}_1  \parallel_2\ \leq \ C\ t^{-1}\ I_{0}\left  ( \parallel <x> \nabla q  \parallel_2\ + \ \parallel  q  \parallel_2\right ) \ + \ \parallel \nabla^2 {\check R}_3  \parallel_2\ ,
\eeq
\beq
\label{5.18e}
\parallel \nabla \partial_t G_1  \parallel_2\ \leq \ C\ I_1\left  ( \parallel \partial_t q  \parallel_2\ + \ t^{-1/2} \parallel  q  \parallel_3\right ) \ + \ \parallel \nabla \partial_t R_3  \parallel_2\ ,
\eeq
\bea
\label{5.19e}
&&\parallel \nabla \partial_t {\check G}_1  \parallel_2\ \leq\  C\ t^{-2}\ I_{-1}\left  ( \parallel <x> q  \parallel_2\right ) \nn  \\
&&+ \ C\ t^{-1}\ I_{0}\left  ( \parallel <x> \partial_t  q  \parallel_2\ +\ t^{-1/2} \parallel <x> q  \parallel_3\right ) \ + \ \parallel \nabla \partial_t {\check R}_3  \parallel_2\ ,
\eea
\beq
\label{5.20e}
\parallel \nabla^k \partial_t \sigma  \parallel_2\ \leq \ C\ t^{-1}\left  ( \parallel q  \parallel_2\ + \ \delta_{k2} \parallel  \nabla q  \parallel_2\right ) \ + \ t^{-k\beta} \parallel \nabla {\check G}_1\parallel_2\ + \ \parallel \nabla^k R_2  \parallel_2
\eeq

\noi for $k = 0, 1, 2$.}\\

\noi {\bf Proof.} From (\ref{2.32e}) (\ref{2.52e}) we obtain
$$G_1 = F_1 (PN_1 ) - R_3$$

\noi where
$$N_1 = - x \ {\rm Re}\ \overline{q} (2w_a + q)$$

\noi so that 
\begin{eqnarray*}
&&\partial_t G_1 = F_2 \left ( P \partial_t N_1 \right ) - \partial_t R_3\ ,\\
&&{\check G}_1 = t^{-1} F_0 \left (x\cdot PN_1 \right ) - {\check R}_3\ , \\
&&\partial_t  {\check G}_1 =-  t^{-2} F_0 \left (x\cdot PN_1 \right ) + t^{-1} F_1 \left (x\cdot P\partial_tN_1 \right ) - \partial_t {\check R}_3\ . 
\end{eqnarray*}

\noi Using (\ref{3.18e}) (\ref{3.19e}), we estimate successively
\begin{eqnarray*}
\parallel \nabla G_1  \parallel_2 &\leq& I_0\left  ( \parallel N_1  \parallel_2\right ) \ + \ \parallel \nabla R_3  \parallel_2\\
&\leq& I_0\left  ( \parallel q  \parallel_2\ \parallel x(2w_a + q)  \parallel_{\infty} \right ) \ + \ \parallel \nabla R_3  \parallel_2\ ,
\end{eqnarray*}
\begin{eqnarray*}
\parallel \nabla {\check G}_1  \parallel_2&\leq& t^{-1}\ I_{-1}\left  ( \parallel <x> N_1  \parallel_2\right )  + \ \parallel \nabla {\check R}_3  \parallel_2\\
 &\leq& t^{-1}\ I_{-1}\left  ( \parallel  x q  \parallel_2\  \parallel < x> (2w_a + q)  \parallel_{\infty} \right ) \ + \ \parallel \nabla {\check R}_3  \parallel_2\ ,
\end{eqnarray*}
\begin{eqnarray*}
&&\parallel \nabla^2 G_1  \parallel_2 \ \leq\  I_1\left  ( \parallel \nabla N_1  \parallel_2\right ) \ + \ \parallel \nabla^2 R_3  \parallel_2\\
&&\leq \ I_1\left  ( 2\parallel \nabla q  \parallel_2\ \parallel x(w_a + q)  \parallel_{\infty} \ +\ \parallel q\parallel_6 \left ( \parallel 2w_a + q  \parallel_3 \ + \ 2 \parallel x\nabla w_a  \parallel_3\right ) \right ) \\
&&+\ \parallel \nabla^2 R_3  \parallel_2\ ,
\end{eqnarray*}
$$\parallel \nabla^2 {\check G}_1  \parallel_2\ \leq\  t^{-1}\ I_{0}\left  ( \parallel <x> \nabla N_1  \parallel_2\right )  + \ \parallel \nabla^2 {\check R}_3  \parallel_2\ , $$
\begin{eqnarray*}
\parallel <x> \nabla N_1  \parallel_2&\leq& 2 \parallel  x \nabla q  \parallel_2\  \parallel <x>(w_a + q)  \parallel_{\infty}  \ + \ \parallel xq  \parallel_6\ \parallel <x>\nabla w_a  \parallel_3\\
&&+\ \parallel q \parallel_2\ \parallel <x> (2w_a + q)\parallel_2\ ,
\end{eqnarray*}
\begin{eqnarray*}
&&\parallel \nabla \partial_t G_1  \parallel_2 \ \leq\  I_1\left  ( \parallel \partial_t N_1  \parallel_2\right ) \ + \ \parallel \nabla \partial_t R_3  \parallel_2\\
&&\leq\ 2 I_1\left  ( \parallel \partial_t  q  \parallel_2\ \parallel x(w_a + q)  \parallel_{\infty} \ +\ \parallel q \parallel_3\ \parallel x \partial_t w_a \parallel_6\right ) \ + \ \parallel \nabla \partial_t R_3  \parallel_2\ ,
\end{eqnarray*}
$$\parallel \nabla  \partial_t  {\check G}_1  \parallel_2\ \leq\  t^{-2}\ I_{-1}\left  ( \parallel <x>  N_1  \parallel_2\right )  + \ t^{-1}\ I_0 \left (\parallel <x> \partial_t N_1  \parallel_2\right ) +  \ \parallel \nabla  \partial_t  {\check R}_3  \parallel_2\ , $$
\begin{eqnarray*}
\parallel <x> \partial_t N_1  \parallel_2&\leq& 2 \parallel  x \partial_t  q  \parallel_2\  \parallel <x>(w_a + q)  \parallel_{\infty}  \\
&&+\ 2\parallel xq \parallel_3\ \parallel <x> \partial_t w_a\parallel_6\ ,
\end{eqnarray*}

\noi from which (\ref{5.14e})-(\ref{5.19e}) follow by using (A1) and (\ref{5.13e}).\par

Finally (\ref{5.20e}) follows from (\ref{2.53e}) (\ref{3.12e}) (A1) (\ref{5.13e}) and Lemma 3.5. \par \nobreak \hfill $\sq$ \par

We next estimate $G'_2$ defined by (\ref{2.59e}).\\

\noi {\bf Lemma 5.2.} {\it Let $0 < \beta < 1$, $0 < \tau \leq \tau_0$ and $I = (0, \tau ]$. Let $w_a$ satisfy (A1). Let $K_a$ satisfy (\ref{5.3e}) and
\beq
\label{5.21e}
\parallel \nabla B_a \parallel_{\infty} \ \vee \ \parallel \partial_t K_a \parallel_{\infty} \ \leq \ C\ t^{-1}\ .
\eeq 

\noi Let $(q, G_2) \in X_0 (I)$ satisfy (\ref{5.13e}) and
\beq
\label{5.22e}
\parallel L \parallel_{\infty}\ \leq C(1 - \ell n\ t)\ .
\eeq

\noi Then the following estimates hold for all $t\in I$~:}
\beq
\label{5.23e}
\parallel \nabla G'_2 \parallel _2\ \leq \ C\ t\ I_1 \left ( \parallel q\parallel _2 (1 - \ell n\ t) \ + \ \parallel \sigma \parallel _2\ + \ \parallel \nabla G \parallel _2 \right ) \ + \ \parallel  \nabla R_4 \parallel _2\ ,
\eeq
\beq
\label{5.24e}
\parallel \nabla {\check G}'_2 \parallel _2\ \leq \ C\ I_0 \left ( \parallel <x> q\parallel _2 (1 - \ell n\ t) \ + \ \parallel \sigma \parallel _2\ + \ \parallel \nabla G \parallel _2 \right ) \ + \ \parallel  \nabla {\check R}_4 \parallel _2\ ,
\eeq
\beq
\label{5.25e}
\parallel \nabla^2 {G}'_2 \parallel _2\ \leq \ C\ t\ I_2 \left ( \parallel \nabla_K q\parallel _2 (1 - \ell n\ t) \ + \ t^{-1}\parallel q \parallel _2\ + \ \parallel \nabla L \parallel _2 \right ) \ + \ \parallel  \nabla^2 {R}_4 \parallel _2\ ,
\eeq
\bea
\label{5.26e}
\parallel \nabla^2 {\check G}'_2 \parallel _2&\leq& C\ I_1 \left ( \parallel <x> \nabla_K q\parallel _2 (1 - \ell n\ t) \ + \ t^{-1}\parallel q \parallel _2\ + \ \parallel \nabla L \parallel _2 \right ) \nn \\
&& + \ \parallel  \nabla^2 {\check R}_4 \parallel _2\ ,
\eea
\bea
\label{5.27e}
&&\parallel \nabla \partial_t G'_2 \parallel _2\ \leq \ C\ I_1 \left ( \parallel q\parallel _2 (1 - \ell n\ t) \ + \ \parallel \sigma \parallel _2\ + \ \parallel \nabla G \parallel _2 \right ) \nn \\
&&+ \ C\ t\ I_2 \Big ( \parallel \partial_t q \parallel _2 (1 - \ell n\ t) \ +\ \parallel \nabla_K q\parallel _2 \ t^{-1/2}\ + \ \ t^{-1}\parallel q \parallel _2\ + \ \parallel \partial_t \sigma \parallel _2\nn \\
&&+ \ \parallel  \nabla \partial_t G  \parallel _2\  + \left (  \parallel \sigma \parallel _2\ + \ \parallel \nabla G \parallel _2\right ) t^{-1/2} \Big ) \ + \  \parallel  \nabla \partial_t R_4 \parallel _2\ ,
\eea
$$\parallel \nabla \partial_t  {\check G}'_2 \parallel _2\ \leq \ C\ I_1 \Big ( \parallel <x> \partial_t q\parallel _2 (1 - \ell n\ t) \ + \ \parallel <x> \nabla_K q\parallel _2 \ t^{-1/2}$$
\beq
\label{5.28e}
+\ t^{-1}\parallel q \parallel _2\ + \ \parallel \partial_t \sigma \parallel _2\ + \ \parallel \nabla \partial_t G \parallel _2\ + \left (  \parallel \sigma \parallel _2\ +\ \parallel \nabla G \parallel _2 \right ) t^{-1/2} \Big ) \ + \ \parallel  \nabla \partial_t  {\check R}_4 \parallel _2\ .
\eeq

\vskip 5 truemm

\noi {\bf Proof.} From (\ref{2.35e}) (\ref{2.59e}) we obtain
\beq
\label{5.29e}
 G'_2 = t\ F_2(PN_2) - R_4
\eeq

\noi where
\beq
\label{5.30e}
N_2 = \ {\rm Im} \ \overline{q} \nabla_K (2w_a + q) - L |w_a|^2
\eeq

\noi so that 
\begin{eqnarray*}
&&{\check G}'_2 = F_1(x\cdot PN_2) - {\check R}_4\ , \\
&&\partial_t G'_2 = F_2 (PN_2) + t\ F_3 (P \partial_t N_2 ) - \partial_t  R_4\ ,\\
&&\partial_t {\check G}'_2 = F_2 \left (x \cdot P \partial_t N_2 \right ) - \partial_t {\check R}_4\ .
\end{eqnarray*}

\noi Using (\ref{3.18e}) (\ref{3.19e}) we first estimate 
$$\parallel \nabla G'_2 \parallel _2\ \leq \ t\ I_1 \left ( \parallel N_2 \parallel _2\right )  \ + \  \parallel  \nabla R_4 \parallel _2\ ,$$
$$\parallel  N_2 \parallel _2\ \leq \ \parallel  q\parallel _2 \ \parallel  \nabla_K (2w_a + q) \parallel _{\infty}\ +\ \parallel  \sigma \parallel _2 \ \parallel  w_a \parallel _{\infty}^2 \ + \ \parallel  G \parallel _6\ \parallel  w_a \parallel _6^2\ ,$$
$$\parallel \nabla {\check G}'_2 \parallel _2\ \leq \  I_0 \left ( \parallel <x> N_2 \parallel _2\right )  \ + \  \parallel  \nabla {\check R}_4 \parallel _2\ ,$$
$$\parallel  <x> N_2 \parallel _2\ \leq \ \parallel  <x> q\parallel _2 \ \parallel  \nabla_K (2w_a + q) \parallel _2\ +\ \parallel  \sigma \parallel _2 \ \parallel  <x> w_a \parallel _{\infty} \ \parallel  w_a \parallel _{\infty}$$
$$+ \ \parallel  G \parallel _6\ \parallel  <x> w_a \parallel _6 \  \parallel  w_a \parallel _6$$

\noi from which (\ref{5.23e}) (\ref{5.24e}) follow. In order to estimate $\nabla^2 G'_2$ and $\nabla^2 {\check G}'_2$ we use the identity $\nabla PN_2 = PN_3$ or more precisely
\beq
\label{5.31e}
\nabla_{\ell} \ P_{ij} \ N_{2j} = P_{ij}\ N_{3j\ell} 
\eeq 

\noi where
\bea
\label{5.32e}
N_{3j\ell} &=& {\rm Im} \left ( \overline{\nabla_{K\ell} q}\ \nabla_{Kj} (2 w_a + q) + \overline{\nabla_{K\ell} (2w_a + q)} \ \nabla_{Kj} q \right )\nn \\
&&+ F_{j\ell} \ {\rm Re}\ \overline{q} (2w_a + q) - \nabla_{\ell} (L_j |w_a|^2 )\ ,
\eea
$$F_{j\ell} = \nabla_j K_{\ell} - \nabla_{\ell} K_j = \nabla_j B_{\ell} - \nabla_{\ell} B_j \ .$$

\noi In fact, using the symmetry of the bilinear form $P\ {\rm Im}\ \overline{u} \nabla_k v$, we obtain
\begin{eqnarray*}
&&\nabla_{\ell} P_{ij}\ {\rm Im}\ \overline{u} \nabla_{Kj}v =  P_{ij}\ {\rm Im}\ \nabla_{\ell} \left ( \overline{u} \nabla_{Kj} v\right ) \\
&&= P_{ij}\ {\rm Im} \left ( \overline{\nabla_{K\ell} u}\ \nabla_{Kj} v + \overline{u} \nabla_{K\ell} \nabla_{Kj} v \right )\\
&&= P_{ij}\ {\rm Im} \left ( \overline{\nabla_{K\ell} u}\ \nabla_{Kj} v + i F_{j\ell} \overline{u} v + \overline{u} \nabla_{Kj} \nabla_{K\ell} v \right )\\
&&= P_{ij} \left \{ {\rm Im} \left ( \overline{\nabla_{K\ell} u}\ \nabla_{Kj} v +  \overline{\nabla_{K\ell} v }\ \nabla_{Kj} u\right ) + F_{j\ell} \ {\rm Re}\  \overline{u}v \right \}\\
\end{eqnarray*}
\vskip - 5 truemm
\noi from which (\ref{5.31e}) (\ref{5.32e}) follow. We can then rewrite $\nabla G'_2$ as
$$\nabla G'_2 = t\ F_3 (P \nabla N_2) - \nabla R_4 = t\ F_3 (PN_3) - \nabla R_4 $$

\noi so that by (\ref{3.18e}) (\ref{3.19e})
$$\parallel \nabla^2 G'_2 \parallel _2\ \leq \ t\ I_2 \left ( \parallel N_3 \parallel _2\right )  \ + \  \parallel  \nabla^2 R_4 \parallel _2\ ,$$
$$\parallel  N_3 \parallel _2\ \leq \ 2\parallel  \nabla_K q\parallel _2 \ \parallel  \nabla_K (2w_a + q) \parallel _{\infty}$$
$$+\ 2\left (  \parallel  \nabla B_a \parallel _{\infty} \ \parallel  q \parallel _2 \ + \ \parallel  \nabla L \parallel _2\ \parallel  q \parallel _{\infty} \right ) \parallel  2w_a + q \parallel _{\infty}$$
$$+ \ \left ( \parallel  \nabla L \parallel _2\ \parallel  w_a \parallel _{\infty} \ + \ 2 \parallel  L \parallel _6\ \parallel  \nabla w_a \parallel _3\right ) \parallel  w_a \parallel _{\infty}\ , $$
\vskip - 5 truemm

$$\parallel \nabla^2 {\check G}'_2 \parallel _2\ \leq \  I_1 \left ( \parallel <x> N_3 \parallel _2\right )  \ + \  \parallel  \nabla^2 {\check R}_4 \parallel _2\ ,$$
$$\parallel  <x> N_3 \parallel _2\ \leq \ 2 \parallel  <x> \nabla_K q\parallel _2 \ \parallel  \nabla_K (2w_a + q) \parallel _{\infty}$$
$$+\ 2\left ( \parallel  \nabla B_a \parallel _{\infty} \ \parallel  q \parallel _2\ + \ \parallel  \nabla L  \parallel _2\ \parallel  q \parallel _{\infty}\right )  \parallel  <x> (2w_a +q) \parallel _{\infty}$$
$$+\ \left ( \parallel  \nabla L \parallel _2 \ \parallel  w_a \parallel _{\infty}\ + \ 2\parallel  L  \parallel _6\ \parallel  \nabla w_a \parallel _3\right ) \ \parallel  <x> w_a \parallel _{\infty}\ ,$$

\noi from which (\ref{5.25e}) (\ref{5.26e}) follow by using (A1) (\ref{5.13e}) (\ref{5.3e}) (\ref{5.21e}) (\ref{5.22e}). We finally estimate $\nabla \partial_t G'_2$ and $\nabla \partial_t {\check G}'_2$. For that purpose, we use the identity $P \partial_tN_2 = PN_4$ with 
\begin{eqnarray*}
N_4 &=& {\rm Im} \left ( ( \partial_t \overline{q}) \nabla_K (2w_a + q) + (\partial_t (2 \overline{w}_a + \overline{q})) \nabla_K q \right )\\
&& - ( \partial_t  K) \ {\rm Re}\  \overline{q} (2w_a + q) - \partial_t (L|w_a|^2)\ .
\end{eqnarray*}

\noi Using again (\ref{3.18e}) (\ref{3.19e}), we estimate
$$\parallel \nabla \partial_t  G'_2 \parallel _2\ \leq \ I_1 \left ( \parallel N_2 \parallel _2\right )  \ + \ t\ I_2 \left ( \parallel N_4 \parallel _2\right )\ + \ \parallel  \nabla \partial_t  R_4 \parallel _2\ ,$$
$$\parallel  N_4 \parallel _2\ \leq \ 2  \parallel \partial_t q \parallel_2 \ \parallel  \nabla_K (w_a + q)\parallel _{\infty}\ + \ 2 \parallel  \nabla_K q \parallel _2\ \parallel \partial_t w_a \parallel_{\infty}$$
$$+\ \left (  \parallel  \partial_t K_a \parallel _{\infty} \ \parallel  q \parallel _2 \ + \ \parallel  \partial_t \sigma \parallel _2\ \parallel  q \parallel _{\infty}\ + \ \parallel \partial_t G \parallel _6 \ \parallel q \parallel_3  \right ) \parallel  2w_a + q \parallel _{\infty}$$
$$+ \  \parallel \partial_t \sigma \parallel _2\ \parallel  w_a \parallel _{\infty}^2\ + \ \parallel \partial_t G \parallel _6 \ \parallel w_a \parallel_6^2 $$
$$+ \ 2\left ( \parallel \sigma \parallel _2\ \parallel  w_a \parallel _{\infty}\ + \ \parallel G \parallel _6 \ \parallel w_a \parallel _3\right ) \parallel  \partial_t  w_a \parallel_{\infty}\ , $$
\vskip - 3 truemm
$$\parallel \nabla  \partial_t  {\check G}'_2 \parallel _2\ \leq \  I_1 \left ( \parallel <x> N_4 \parallel _2\right )  \ + \  \parallel  \nabla  \partial_t  {\check R}_4 \parallel _2\ ,$$
\vskip - 5 truemm
$$\parallel  <x> N_4 \parallel _2\ \leq \ 2 \parallel  <x>  \partial_t   q\parallel _2 \ \parallel  \nabla_K (w_a + q) \parallel _{\infty}\ + \ 2 \parallel  <x>  \nabla_K q\parallel _2\ \parallel  \partial_t   w_a\parallel _{\infty}$$
$$+\ \left ( \parallel  \partial_t  K_a \parallel _{\infty} \ \parallel  q \parallel _2\ + \ \parallel  \partial_t \sigma  \parallel _2\ \parallel  q \parallel _{\infty}\ +\ \parallel  \partial_t  G \parallel _6\ \parallel  q \parallel_3 \right )  \parallel  <x> (2w_a +q) \parallel _{\infty}$$
$$+ \ \parallel  \partial_t \sigma  \parallel _2\ \parallel  w_a \parallel _{\infty}\  \parallel  <x> w_a  \parallel _{\infty}\ +\ \parallel  \partial_t  G \parallel _6\ \parallel  w_a \parallel_6  \ \parallel  <x> w_a  \parallel _6$$
$$+\ 2\left ( \parallel  \sigma \parallel _2 \ \parallel  <x> w_a \parallel _{\infty}\ + \ \parallel  G  \parallel _6\ \parallel  <x> w_a \parallel _3\right ) \ \parallel  \partial_t w_a \parallel _{\infty}\ ,$$

\noi from which (\ref{5.27e}) (\ref{5.28e}) follow by using (A1) (\ref{5.13e}) (\ref{5.3e}) (\ref{5.21e}) (\ref{5.22e}).\par \nobreak \hfill $\sq$\par

We next estimate $H_1w_a$ with $H_1$ defined by (\ref{2.56e}).\\

\noi {\bf Lemma 5.3.} {\it Let $0 < \beta < 1$, $0 < \tau \leq \tau_0$ and $I = (0, \tau ]$. Let $w_a$ satisfy (A1). Let $K_a$ satisfy (\ref{5.3e}) (\ref{5.4e}). Let $(q, G_2) \in X_0(I)$ satisfy (\ref{5.13e}) and 
\beq
\label{5.33e}
\parallel L \parallel_{\infty} \ \vee \ \parallel  \nabla L \parallel_3 \ \leq \ C(1 - \ell n\ t) 
\eeq 

\noi for all $t \in I$. Then the following estimates hold for all $t\in I$~:}
\bea
\label{5.34e}
\parallel <x> H_1 w_a\parallel _2 &\leq& C\Big ( \left ( \parallel \sigma \parallel _3\ + \ \parallel \nabla G \parallel _2\right )   (1 - \ell n\ t)\nn \\
 &+& \parallel \nabla \cdot \sigma \parallel _2\ +  \ t^{\beta} \parallel \nabla {\check G}_1 \parallel _2 \ + \ \parallel  \nabla {\check G}_2  \parallel _2\Big ) \ ,    \eea
\bea
\label{5.35e}
\parallel <x> H_1 \partial_t w_a\parallel _2 &\leq& C \Big ( \parallel L \parallel _{\infty}\ + \ \parallel L \parallel _6   (1 - \ell n\ t)\nn \\
 &+&  \parallel \nabla \cdot \sigma \parallel _3\ +  \ t^{\beta /2} \parallel \nabla {\check G}_1 \parallel _2 \ + \ \parallel  \nabla {\check G}_2  \parallel _2\Big ) t^{-1/2} ,
\eea
$$\parallel <x> ( \partial_t H_1) w_a\parallel _2 \ \leq\  C \Big ( \left ( \parallel \partial_t \sigma  \parallel _3\ + \ \parallel \nabla \partial_t G \parallel _2 \right )  (1 - \ell n\ t)$$
 \beq
\label{5.36e}
 + \ \parallel \partial_t \nabla \cdot \sigma \parallel _2\ +  \ t^{\beta} \parallel \nabla \partial_t {\check G}_1 \parallel _2 \ + \ \parallel  \nabla \partial_t {\check G}_2  \parallel _2\ + \ t^{-1} \left ( \parallel \sigma  \parallel _2\ + \ \parallel \nabla G \parallel _2\right ) \Big ) \ ,
\eeq
$$\parallel \nabla_K H_1 w_a\parallel _2 \ \leq\  C \Big (\parallel \nabla L  \parallel _2 (1 - \ell n\ t)^2\ + \ \parallel \nabla \cdot \sigma  \parallel _3 (1 - \ell n\ t)\ + \  \parallel \nabla \nabla \cdot \sigma  \parallel_2$$
   \beq
\label{5.37e}
 + \ t^{-1}\left (  \parallel \ \sigma \parallel _2\ +  \  \parallel \nabla G \parallel_2 \right ) \ +\ \parallel \nabla {\check G}_1 \parallel _2 \ + \ \parallel  \nabla  {\check G}_2  \parallel _2\ (1 - \ell n\ t) \Big ) \ ,
\eeq
\bea
\label{5.38e}
&&\parallel \nabla_K H_1 \partial_t  w_a\parallel _2 \ \leq\  C \Big ( \left ( \parallel L  \parallel _{\infty}  \ + \ \parallel \nabla L  \parallel _3\right )  (1 - \ell n\ t)\ + \ \parallel \nabla L  \parallel _2 (1 - \ell n\ t)^2\nn \\
&&\ + \  \parallel \nabla  \cdot \sigma  \parallel_3 (1 - \ell n\ t)\ + \ \parallel \nabla \nabla \cdot \sigma  \parallel_2 \ + \ t^{-1}\left (  \parallel \sigma \parallel _2\ +  \  \parallel \nabla G \parallel_2 \right )\nn \\
&&+\ \parallel \nabla {\check G}_1 \parallel _2 \ + \ \parallel  \nabla  {\check G}_2  \parallel _2\ (1 - \ell n\ t) \Big ) t^{-1/2} \ ,
\eea
$$\parallel \nabla_K (\partial_t  H_1) w_a\parallel _2 \ \leq\  C \Big (\parallel \nabla \partial_t L  \parallel _2(1 - \ell n\ t)^2  \ + \ \parallel \nabla \partial_t \nabla \cdot \sigma   \parallel _2$$
$$ + \ t^{-1} \left  ( \parallel \partial_t \sigma   \parallel _2 \ + \  \parallel \nabla  \partial_t G  \parallel_2 \ + \ \parallel \nabla L  \parallel_2 (1 - \ell n\ t)\right )$$
\beq
\label{5.39e}
 + \ t^{-2}\left (  \parallel \sigma \parallel _2\ +  \  \parallel \nabla G \parallel_2 \right )\ +\ \parallel \nabla \partial_t {\check G}_1 \parallel _2 \ + \ \parallel  \nabla  \partial_t {\check G}_2  \parallel _2\ (1 - \ell n\ t)\Big ) \ .
\eeq
\vskip 5 truemm

\noi {\bf Proof.}

\noi (\ref{5.34e}) and (\ref{5.35e}). We rewrite $H_1$ defined by (\ref{2.56e}) as
\beq
\label{5.40e}
H_1 = i L \cdot \nabla_{K_a} + (i/2) (\nabla \cdot \sigma ) + (1/2) L^2 + {\check G}_S\ .
\eeq

\noi Using (\ref{3.11e}), we estimate
$$\parallel <x> H_1 w_a\parallel _2 \ \leq\ \parallel \sigma \parallel _3\ \parallel  <x> \nabla_{K_a} w_a \parallel _6 \ + \ \parallel G \parallel_6 \ \parallel <x> \nabla_{K_a} w_a\parallel _3$$
$$ + \ \left ( \parallel \nabla \cdot \sigma  \parallel _2\ + \ t^{\beta}  \parallel \nabla {\check G}_1 \parallel _2\right )  \parallel <x> w_a \parallel_{\infty} \ + \parallel {\check G}_2\parallel_6 \ \parallel <x>w_a\parallel_3 $$
$$+\  \parallel L \parallel _{\infty}\left ( \parallel \sigma \parallel_3\ \parallel <x> w_a \parallel_6\ +\ \parallel G  \parallel _6 \ \parallel  <x> w_a\parallel _3\right  ) \ , $$  
 $$\parallel <x> H_1 \partial_t w_a\parallel _2 \ \leq\  \parallel L \parallel _{\infty}\ \parallel <x> \nabla \partial_t w_a \parallel_2\ + \ \parallel L \parallel _6\    
 \parallel K_a \parallel_{\infty}\ \parallel <x> \partial_t w_a \parallel_3$$
 $$+\ \left ( \parallel \nabla \cdot \sigma \parallel _3\ +  \ t^{\beta /2} \parallel \nabla {\check G}_1 \parallel _2 \right ) \parallel <x> \partial_t w_a \parallel_6\  + \ \parallel  {\check G}_2  \parallel _6\ \parallel <x> \partial_t w_a \parallel_3$$
$$+\  \parallel L \parallel _{\infty}\ \parallel L \parallel _6\ \parallel <x> \partial_t w_a \parallel_3$$

\noi from which (\ref{5.34e}) (\ref{5.35e}) follow by using (A1) (\ref{5.3e}) (\ref{5.33e}).\\

\noi (\ref{5.36e}). Taking the time derivative of (\ref{5.40e}), we obtain
\beq
\label{5.41e}
\partial_t H_1 = i (\partial_t L) \cdot \nabla_{K_a} + (i/2) (\partial_t \nabla \cdot \sigma ) + L\cdot (\partial_t K) + (\partial_t {\check G}_S) \ .
\eeq

\noi Using (\ref{3.11e}), we estimate
$$\parallel <x> ( \partial_t H_1) w_a\parallel _2 \ \leq\  \parallel \partial_t \sigma  \parallel _3\ \parallel <x> \nabla_{K_a} w_a \parallel_6\ + \ \parallel \partial_t G \parallel _6\ \parallel <x> \nabla_{K_a} w_a \parallel_3$$
$$ + \ \left ( \parallel \partial_t \nabla \cdot \sigma \parallel _2\ +  \ t^{\beta} \parallel \nabla \partial_t {\check G}_1 \parallel _2 \right ) \parallel <x> w_a \parallel_{\infty} \ +\  \parallel  \partial_t {\check G}_2  \parallel _6\ \parallel <x>  w_a \parallel_3$$
$$+\  \parallel \partial_t K_a  \parallel _{\infty} \left ( \parallel \sigma  \parallel _2\  \parallel <x>w_a  \parallel _{\infty}\ + \ \parallel G \parallel _6\ \parallel <x>w_a  \parallel _3 \right ) $$
$$+\  \parallel L \parallel _{\infty} \left ( \parallel \partial_t \sigma  \parallel _3\  \parallel <x>w_a  \parallel_6\ + \ \parallel \partial_t G \parallel _6\ \parallel <x>w_a  \parallel _3 \right ) $$

\noi from which (\ref{5.36e}) follows by using (A1) (\ref{5.3e}) (\ref{5.4e}) (\ref{5.33e}).\\

\noi (\ref{5.37e}) and (\ref{5.38e}). Taking the covariant gradient of $H_1v$, using (\ref{5.40e}) and the identity
\beq
\label{5.42e}
\nabla_K \nabla_{K_a} = \nabla_{K_a} \nabla - i K \nabla_{K_a} - i (\nabla K_a)
\eeq

\noi we obtain
$$\nabla_K H_1 v = iL \cdot \nabla_{K_a} \nabla v + (i (\nabla L) + KL) \cdot \nabla_{K_a} v$$
\beq
\label{5.43e}
+ \left ( (i/2) (\nabla \cdot \sigma ) + (1/2) L^2 + {\check G}_S \right ) \nabla_K v + \left ( L\cdot (\nabla K) + (i/2) (\nabla \nabla \cdot \sigma ) + (\nabla {\check G}_s)\right ) v\ .
\eeq 

\noi Using (\ref{3.11e}), we estimate
$$\parallel \nabla_K H_1 v\parallel _2 \ \leq\  \parallel \nabla L  \parallel _2 \left (   \parallel \nabla_{K_a} \nabla v   \parallel _3\ + \ \parallel \nabla_{K_a}  v   \parallel _{\infty} \ + \ \parallel K  \parallel _{\infty} \ \parallel \nabla_{K_a} v   \parallel _3\right )$$
$$+ \  \parallel \nabla  \cdot \sigma  \parallel_3 \ \parallel \nabla_{K} v   \parallel _6\ + \ \parallel \nabla \nabla \cdot \sigma  \parallel_2 \   \parallel v\parallel_{\infty}$$
$$+  \  \parallel \nabla L \parallel_2 \  \parallel L \parallel_{\infty}\left ( \  \parallel \nabla_K v \parallel_3 \ + \ \  \parallel v \parallel_{\infty} \right )$$
$$+\ \  \parallel \nabla K_a \parallel_{\infty} \left ( \  \parallel \sigma \parallel_2\ \  \parallel v \parallel_{\infty} \ + \ \  \parallel G \parallel_6\ \  \parallel v\parallel_3 \right )$$
\beq
\label{5.44e}
+\ \parallel \nabla {\check G}_1 \parallel _2 \left ( t^{\beta /2} \  \parallel \nabla_K v  \parallel_6\  +  \  \parallel v \parallel_{\infty} \right )\ +\  \parallel  \nabla  {\check G}_2  \parallel _2\ \left ( \  \parallel \nabla_K v  \parallel_3\ +\ \  \parallel v \parallel_{\infty}\right ) \ .
\eeq

\noi Applying (\ref{5.44e}) with $v = w_a$ and using (A1)  (\ref{5.3e}) (\ref{5.4e}) (\ref{5.33e}) yields (\ref{5.37e}), while (\ref{5.38e}) follows from (\ref{5.44e}) with $v = \partial_t w_a$, except for a slightly different estimate of the contribution of the first two terms in the right hand side of (\ref{5.43e}), namely  
$$\parallel L \cdot \nabla_{K_a} \nabla \partial_t w_a  \parallel_2 \ + \ \parallel (\nabla L) \nabla_{K_a} \partial_t w_a \parallel_2$$
$$\leq \ \parallel L \parallel_{\infty} \ \parallel \nabla_{K_a}  \nabla \partial_t w_a \parallel_2\ + \   \parallel \nabla L \parallel_3 \ \ \parallel \nabla_{K_a} \partial_t w_a \parallel_6\ .$$
\vskip 5 truemm

\noi (\ref{5.39e}). Taking the covariant gradient of $(\partial_t H_1)v$, using (\ref{5.41e}) (\ref{5.42e}), we obtain
$$ \nabla_K (\partial_t  H_1) v = i ( \partial_t L)  \cdot  \nabla_{K_a} \nabla v +\Big ( i (\nabla  \partial_t L) + K(\partial_t L) \Big ) \cdot  \nabla_{K_a}v$$   
$$ +  \left  ((i/2) ( \partial_t \nabla \cdot \sigma)  + L \cdot (  \partial_t K) +   (\partial_t {\check G}_S)\right )  \nabla_K v$$
$$+ \left (  (i/2) (\nabla \partial_t \nabla \cdot \sigma ) + (\partial_t L)\cdot (\nabla K_a) + (\nabla L) \cdot (\partial_t K) + L (\nabla \partial_t K) + ( \nabla \partial_t {\check G}_S) \right ) v\ .$$

\noi Using (\ref{3.11e}), we estimate
$$\parallel \nabla_K (\partial_t  H_1) v\parallel _2 \ \leq\  \parallel \nabla \partial_t L  \parallel _2\ \left ( \parallel \nabla_{K_a} \nabla v \parallel_3\ + \  \parallel \nabla_{K_a} v \parallel _{\infty}\ + \ \parallel K \parallel _{\infty} \ \parallel \ \nabla_{K_a} v \parallel _3\right )$$
$$ + \  \parallel \partial_t \nabla \cdot \sigma   \parallel _2 \ \parallel  \nabla_K v \parallel _{\infty}\ + \  \parallel \nabla  \partial_t \nabla \cdot \sigma   \parallel_2 \ \parallel v  \parallel_{\infty}$$
$$+\ \parallel \nabla K_{a} \parallel _{\infty}\left ( \parallel  \partial_t \sigma \parallel _2\ \parallel v \parallel_{\infty} \ +\ \parallel \partial_t G \parallel _6 \ \parallel  v  \parallel _3\right )$$
$$+\ \parallel \partial_t K_a\parallel _{\infty}\  \parallel \nabla L \parallel _2\ \left ( \parallel \nabla_K v \parallel_3 \ + \parallel  v  \parallel _{\infty} \right) $$
$$ + \  \parallel \nabla \partial_t K_a   \parallel _{\infty} \left (  \ \parallel  \sigma \parallel _2 \  \parallel v  \parallel_{\infty}\ + \  \parallel G   \parallel_6 \ \parallel v  \parallel_3\right )$$
$$+\ \parallel \nabla \partial_t L  \parallel _2\Big ( \parallel L  \parallel_{\infty} \ \parallel \nabla_K v  \parallel_3\ +\ \left ( \parallel \nabla L  \parallel_3\ + \ \parallel L  \parallel_{\infty}\right )\parallel v  \parallel_{\infty} \Big )$$
$$ + \  \parallel \nabla \partial_t {\check G}_1 \parallel _2 \left ( t^{\beta} \parallel \nabla_Kv \parallel _{\infty}\ +  \  \parallel v \parallel_{\infty} \right )\ +\ \parallel \nabla \partial_t {\check G}_2 \parallel _2 \left (  \parallel  \nabla_K  v  \parallel _3\ + \ \parallel v \parallel_{\infty} \right )$$

\noi from which (\ref{5.39e}) follows by substituting $w_a$ for $v$ and using (A1) (\ref{5.3e}) (\ref{5.4e}) (\ref{5.33e}). \par \nobreak \hfill $\sq$\par

We next estimate the solutions $q'$ of the Schr\"odinger equation in (\ref{2.59e}). That equation is of the type of (\ref{4.1e}) with $v = q'$, $V = {\check B}_S$ and $f_0 = - \widetilde{R}_1$. The new estimates are extensions and refinements of those derived for $v_{\eta}$ in the proof of Proposition 4.1, starting from the conservation laws in differential form (\ref{4.35e}) (\ref{4.38e}) (\ref{4.42e}), but we shall now take into account the asymptotic behaviour at time zero, which was not considered in the latter proof. In addition, we shall also estimate $xq'$ and $x \partial_tq'$ from the analogues of the conservation laws (\ref{4.8e}) (\ref{4.9e}). The estimates are written in differential form for brevity, but should be understood in integral form like the conservation laws  (\ref{4.2e}) (\ref{4.3e}) (\ref{4.4e}) and  (\ref{4.8e}) (\ref{4.9e}).\\

\noi {\bf Lemma 5.4.} {\it Let $0 < \beta < 1$, $0 < \tau \leq \tau_0$ and $I = (0, \tau ]$. Let $(w_a,K_a)$ satisfy (A1) (A2) and let $B_0$ satisfy (\ref{3.32e}) for $0 \leq j, k \leq 1$ and $r = \infty$. Let $(q, G_2) \in X_0(I)$ satisfy (\ref{5.13e}) and (\ref{5.33e}) and in addition
\beq
\label{5.45e}
\parallel \nabla {\check G}_S \parallel _2\ \leq C\ t^{-1/4} \ ,
\eeq
\beq
\label{5.46e}
\parallel \nabla \partial_t L \parallel_2\ \vee \ \parallel \nabla {\check G}_S \parallel _6\ \leq C\ t^{-3/4} \ ,
\eeq  
\beq
\label{5.47e}
\parallel \nabla \partial_t \nabla \cdot \sigma  \parallel_2\ \vee \ \parallel \nabla \partial_t {\check G}_S \parallel_2\  \leq C\ t^{-5/4} 
\eeq  

\noi for all $t \in I$. Let $q'$ with $(q', 0) \in X_0 (I)$ be a solution of the Schr\"odinger equation in (\ref{2.59e}), satisfying the conservation laws (\ref{4.2e}) (\ref{4.3e}) (\ref{4.4e}) (\ref{4.8e}) (\ref{4.9e}). Then the following estimates hold for all $t\in I$~:}
\beq
\label{5.48e}
\left | \partial_t \parallel q' \parallel_2 \right | \ \leq \ \parallel  \widetilde{R}_1 \parallel_2\ ,
\eeq
\beq
\label{5.49e}
\left | \partial_t \parallel xq' \parallel_2 \right | \ \leq \ \parallel  \nabla_K q' \parallel _2\ +\ \parallel x \widetilde{R}_1 \parallel_2\ ,
\eeq
\beq
\label{5.50e}
\left | \partial_t \parallel \partial_t q' \parallel_2 \right | \ \leq \ C \left ( t^{-1} \parallel  \nabla_K q' \parallel _2\ +\ t^{-3/4}\parallel  \nabla_K q' \parallel _3\ +\ t^{-3/2} \parallel  q' \parallel _2\right ) \ +\  \parallel \partial_t \widetilde{R}_1 \parallel \ ,
\eeq
\bea
\label{5.51e}
&&\left | \partial_t \parallel x \partial_t q' \parallel_2 \right | \ \leq \  \parallel  \nabla_K \partial_t q' \parallel _2\ +\ C \Big ( t^{-1} \parallel  x\nabla_K q' \parallel _2\ +\ t^{-3/4}\parallel  x\nabla_K q' \parallel _3\nn \\
&&+\ t^{-3/2} \parallel  xq' \parallel _2\Big ) \ +\  \parallel x\partial_t \widetilde{R}_1 \parallel_2 \ ,
\eea
$$\left | \partial_t \parallel \nabla_K \partial_t q' \parallel_2 \right | \ \leq \  C\Big ( t^{-1}\left ( \parallel   \partial_t q' \parallel _2\ + \ \parallel  \nabla_K^2 q' \parallel _2\right )$$
$$ t^{-3/4} \left (  \parallel   \partial_t q' \parallel _3\ + \ \parallel  \nabla_K^2 q' \parallel _3\ +\ \parallel  \nabla_K q' \parallel _{\infty}\right ) \ +\ t^{-3/2}
 \parallel  \nabla_K q' \parallel _2$$
\beq
\label{5.52e}
+\ t^{-5/4}\parallel  q' \parallel _{\infty}\ +\ t^{-2}\parallel  q' \parallel _2\Big ) \ +\ b_0\ t^{-2} \parallel  \nabla_K q'\parallel _2 \ +\  \parallel \nabla_K \partial_t \widetilde{R}_1 \parallel_2 \ ,
\eeq
\beq
\label{5.53e}
 \parallel  <x> \Delta_K q' \parallel _2\ \leq \  \parallel  <x> \partial_tq' \parallel _2\ + \ C\ t^{-1/2}  \parallel  <x>  q' \parallel _2\ + \ \parallel <x>\widetilde{R}_1 \parallel_2 \ ,
\eeq
\beq
\label{5.54e}
\parallel \nabla_K \Delta_K q' \parallel_2 \  \leq \ \parallel \nabla_K \partial_t q' \parallel_2\ + \ C\left ( t^{-1/2} \parallel  \nabla_K q'\parallel _2\ + \ t^{-1} \parallel  q'\parallel _2\right ) \ + \ \parallel \nabla_K  \widetilde{R}_1 \parallel_2 \ .
\eeq
\vskip 5 truemm

\noi {\bf Proof.} Before starting the proof, we remark that the conditions (\ref{5.45e})-(\ref{5.47e}) will eventually become asymptotic region conditions, namely upper bounds on $\tau$ (see (\ref{5.76e}) below). On the other hand under the regularity assumption $(q', 0) \in X_0 (I)$, all the conservation laws are satisfied except possibly (\ref{4.4e}), thereby making the assumption on $q'$ partly redundant.\par

We now begin the proof. The estimate (\ref{5.48e}) follows immediately from (\ref{2.59e}). From the commutation relation
\beq
\label{5.55e}
[x, H] = \nabla_K
\eeq

\noi it follows that $q'$ satisfies the equation
\beq
\label{5.56e}
i \partial_t xq' = \nabla_K q' + Hxq' - x  \widetilde{R}_1
\eeq 
 
 \noi from which (\ref{5.49e}) follows immediately.
 
 We next estimate $\partial_t q'$. Taking the time derivative of the first equation of (\ref{2.59e}) yields 
 \beq
\label{5.57e}
i \partial_t \partial_t q' = H \partial_t q' +  (\partial_tH) q' - \partial_t   \widetilde{R}_1
\eeq 

\noi with
\beq
\label{5.58e}
\partial_t H = i (\partial_t K) \cdot \nabla_K + (i/2) (\partial_t \nabla \cdot s) + \partial_t {\check B}_S
\eeq 
 
 \noi so that
$$ \left | \partial_t \parallel \partial_t q' \parallel_2 \right | \ \leq \ \parallel (\partial_t H) q' \parallel _2\ +\ \parallel \partial_t \widetilde{R}_1 \parallel_2 \ ,$$

$$\parallel (\partial_t H)q' \parallel_2 \ \leq \ \parallel \partial_t K_a \parallel_{\infty}\ \parallel \nabla_K q' \parallel_2 \ + \ \parallel \partial_t L \parallel_6\ \parallel \nabla_K q' \parallel_3$$
$$+ \left ( \parallel \partial_t \nabla \cdot s_a\parallel_{\infty} \ +\ \parallel \partial_t {\check B}_{aS} \parallel_{\infty} \right ) \parallel q'\parallel_2 \ + \ \left ( \parallel \partial_t \nabla \cdot \sigma \parallel_6 \ + \ \parallel \partial_t {\check G}_S \parallel_6 \right ) \parallel q' \parallel_3\ ,$$
 
\noi from which (\ref{5.50e}) follows by the use of (A2) (\ref{5.46e}) (\ref{5.47e}) and a covariant Sobolev inequality.

We next estimate $x\partial_t q'$. From (\ref{5.55e}) (\ref{5.57e}) we obtain 
 \beq
\label{5.59e}
i \partial_t x \partial_t q' = \nabla_K \partial_t q' + Hx \partial_t q' + x (\partial_t H) q' - x \partial_t \widetilde{R}_1
\eeq
 
 \noi so that 
$$\partial_t \parallel x \partial_t q'  \parallel_2\ \leq\  \parallel \nabla_K \partial_t q'  \parallel_2\ +\  \parallel x (\partial_t H)q'  \parallel_2\ +\  \parallel  x \partial_t  \widetilde{R}_1 \parallel_2$$

\noi from which (\ref{5.51e}) follows by the same estimates as before, with $\nabla_K q'$ and $q'$ replaced by $x \nabla_Kq'$ and $xq'$ respectively.\par

We next estimate $\nabla_K \partial_t q'$. Taking the covariant gradient of (\ref{5.57e}) yields 
\bea
\label{5.60e}
i \partial_t \nabla_K  \partial_t q' &=& \left (  \partial_t K + \nabla {\check B}_S \right )  \partial_t q' - (1/2) \nabla_K \Delta_K  \partial_t q' + {\check B}_S \nabla_K  \partial_t q'\nn \\
&&+ \ \nabla_K ( \partial_t H) q' - \nabla_K  \partial_t \widetilde{R}_1
\eea

\noi so that 
\beq
\label{5.61e}
 \partial_t \parallel \nabla_K  \partial_t q'\parallel_2 \ \leq \  \parallel \left (  \partial_t K + \nabla {\check B}_S \right )  \partial_t q' \parallel_2\ +\ \parallel \nabla_K (  \partial_t H)q'\parallel _2\ + \ \parallel  \nabla_K \partial_t \widetilde{R}_1 \parallel _2
\eeq

\noi with
\begin{eqnarray*}
\nabla_K (\partial_t H) q' &= &i ( \partial_t K)\cdot \nabla_K^2 q' + \left ( i(\nabla \partial_t K) + (i/2) (\partial_t \nabla \cdot s) + (\partial_t  {\check B}_S)\right ) \cdot  \nabla_K  q' \\
&&+ \left ( (i/2) (\nabla \partial_t \nabla \cdot s) + (\nabla \partial_t  {\check B}_S)\right )q' \  .  \end{eqnarray*}

\noi We estimate the various terms of (\ref{5.61e}) successively.
$$\parallel ( \partial_t K + \nabla {\check B}_S) \partial_t q' \parallel _2\ \leq \parallel \partial_t K_a + \nabla {\check B}_{aS} \parallel _{\infty} \ \parallel  \partial_t q' \parallel _2$$
$$+\ \parallel \partial_t L + \nabla {\check G}_S \parallel _6\ \parallel  \partial_t q' \parallel _3 \ \leq C \left ( t^{-1} \parallel  \partial_t q' \parallel _2 \ +\ t^{-3/4} \parallel  \partial_t q'\parallel _3\right )\ ,$$

$$\parallel  (\partial_t K) \nabla_K^2 q' \parallel _2\ \leq C \left ( t^{-1} \parallel \nabla_K^2 q' \parallel _2\ + \ t^{-3/4} \parallel \nabla_K^2 q' \parallel _3 \right )\ ,$$

$$\parallel (\nabla \partial_t K) \nabla_K q' \parallel _2\ +\ \parallel (\partial_t \nabla \cdot s ) \nabla_K q'\parallel_2 \ +\ \parallel (\partial_t {\check B}_S) \nabla_K q' \parallel _2$$
$$\leq \ \left ( \parallel \nabla \partial_t K_a \parallel _{\infty} \ + \ \parallel  \partial_t \nabla \cdot s_a \parallel _{\infty} \ + \ \parallel  \partial_t {\check B}_{aS} \parallel _{\infty} \right ) \parallel  \nabla_K q'\parallel _2$$
$$+ \ \parallel  \nabla \partial_t L \parallel _2\ \parallel  \nabla_K q' \parallel _{\infty} \ + \ \left ( \parallel  \partial_t \nabla \cdot \sigma \parallel _6 \ +\ \parallel \partial_t {\check G}_S \parallel _6 \right ) \parallel  \nabla_K q'\parallel _3$$
$$\leq \ \left ( b_0 \ t^{-2} + C \ t^{-3/2}\right ) \parallel  \nabla_K q' \parallel _2 \ + \ C\ t^{-3/4} \parallel \nabla_K q' \parallel_{\infty}\ ,$$

$$\parallel (\nabla \partial_t \nabla \cdot s ) q'\parallel_2 \ +\ \parallel (\nabla \partial_t {\check B}_S )q' \parallel_2$$
$$\leq\ \parallel \nabla \partial_t \nabla \cdot s_a \parallel_3  \  \parallel q'  \parallel _6\ +\ \parallel \nabla \partial_t {\check B}_{aS} \parallel_{\infty}  \ \parallel  q'  \parallel _2$$
$$+ \ \left ( \parallel \nabla \partial_t \nabla \cdot \sigma \parallel_2 \ +\ \parallel \nabla \partial_t {\check G}_S \parallel_2\right ) \parallel  q'  \parallel _{\infty}$$
$$\leq\ C \left ( t^{-5/4} \parallel  q'  \parallel _{\infty}\ +\ t^{-2} \parallel  q'  \parallel _2\right )$$

\noi by (A2) and (\ref{5.46e}) (\ref{5.47e}), and by eliminating some intermediate norms of $\nabla_K q'$ and $q'$ by the H\"older inequality. Note also that in the contribution of $\nabla \partial_t K_a \cdot \nabla_K q'$, the $t^{-2}$ decay comes only from $B_0$, so that the corresponding term appears with a coefficient $b_0$ instead of a more general constant $C$. Substituting the previous estimates into (\ref{5.61e}) yields (\ref{5.52e}).\par

We next estimate $<x> \Delta_K q'$ and $\nabla_K \Delta_K q'$. From (\ref{2.59e}) we obtain immediately 
$$\parallel  <x> \Delta_K q'  \parallel _2\ \leq \  \parallel  <x> \partial_t q'  \parallel _2\ + \ \parallel    {\check B}_S \parallel _{\infty}\  \parallel <x>  q'  \parallel _2 \ + \parallel   <x> \widetilde{R}_1 \parallel _2\ , $$
$$\parallel  \nabla_K  \Delta_K q'  \parallel _2\ \leq \ \parallel  \nabla_K  \partial_t q'  \parallel _2\ + \ \parallel    {\check B}_S \parallel _{\infty}\  \parallel \nabla_K   q'  \parallel _2 \ + \parallel  (\nabla {\check B}_S) q' \parallel _2\ + \ \parallel \nabla_K \widetilde{R}_1 \parallel_2\ . $$

\noi Now
\begin{eqnarray*}
 \parallel  (\nabla {\check B}_S) q' \parallel _2 &\leq&  \parallel  \nabla {\check B}_{aS}  \parallel _{\infty} \ \parallel  q'  \parallel _2 \ +  \parallel  \nabla {\check G}_S)  \parallel _3\ \parallel  q'  \parallel _6\\
&\leq& C \left ( t^{-1} \parallel  q'  \parallel _2\ +\ t^{-1/2} \parallel  \nabla_K q'  \parallel _2\right )\ .
\end{eqnarray*}

\noi by (A2) (\ref{5.45e}) (\ref{5.46e}). This implies (\ref{5.53e}) (\ref{5.54e}). \par \nobreak \hfill $\sq$ \par

In order to estimate $q'$, we shall need in addition some estimates relating covariant and noncovariant derivatives. The following lemma holds pointwise in time and does not require other assumption than finiteness of the norms occurring in the estimates.\\

\noi {\bf Lemma 5.5.} {\it We define
\bea
\label{5.62e}
&&m = \ \parallel K \parallel_{\infty}^2 \ + \  \parallel \nabla \cdot \sigma \parallel_3^2\ + \ \parallel \nabla \cdot s_a \parallel_{\infty} \ +\ \parallel \nabla \cdot \sigma \parallel_3\ , \\
&&\overline{m} = \ \parallel K \parallel_{\infty}^2 \ + \  \parallel \nabla L \parallel_3^2\ + \ \parallel \nabla K_a \parallel_{\infty} \ .
\label{5.63e}
\eea

\noi Then the following estimates hold~:
\bea
\label{5.64e}
&&\parallel <x> \nabla_K v \parallel_3 \ \leq\  \parallel <x> \nabla v \parallel_2^{1/2} \left ( \parallel <x> \Delta v \parallel_2\ +\ \parallel \nabla v \parallel_2 \right )^{1/2}\nn \\
&&+\parallel K \parallel_{\infty} \ \parallel <x> v \parallel_2^{1/2} \left ( \parallel <x> \nabla v \parallel_2 \ +\ \parallel v \parallel_2 \right )^{1/2}\ ,
\eea
\beq
\label{5.65e}
\parallel  \nabla_K v \parallel_{\infty} \ \leq\ \parallel  \nabla \Delta v \parallel_2^{1/2} \left ( \parallel \Delta v \parallel_2\ +\ \parallel K \parallel_{\infty} \ \parallel  v \parallel_2 \right )^{1/2}\ ,
\eeq
\beq
\label{5.66e}
 \parallel <x> \Delta_K v\parallel_{r} \ \leq \  \parallel <x> \Delta v  \parallel_r \ + \ m \parallel <x>  v  \parallel_r \ ,
\eeq
\beq
\label{5.67e}
 \parallel <x> \Delta v\parallel_{r} \ \leq \  \parallel <x> \Delta_K v  \parallel_r \ + \ m \parallel <x>  v  \parallel_r \ ,
\eeq
\beq
\label{5.68e}
 \parallel \nabla_K^2 v\parallel_{r} \ \leq \  \parallel \Delta v  \parallel_r \ + \ \overline{m} \parallel   v  \parallel_r 
\eeq

\noi for $2 \leq r \leq 3$,}
$$\parallel \nabla \Delta  v\parallel_2 \ \leq \  \parallel \nabla_K \Delta_K v\parallel_2 \ + \left (  \parallel \nabla\nabla \cdot s_a \parallel_3\ +\  \parallel \nabla K_a \parallel_{\infty} \right )  \parallel \nabla v\parallel_2$$
\beq
\label{5.69e}
+ \left (  \parallel K\parallel_{\infty}^3 \ +\  \parallel \nabla L\parallel_3^3 \ +\ \parallel K\parallel_{\infty}\ \parallel \nabla K_a \parallel_{\infty}\ +\ \parallel \nabla \nabla \cdot \sigma\parallel_2^2\right ) \parallel v\parallel_2\ .
\eeq
\vskip 5 truemm

\noi {\bf Proof.} The proof uses H\"older and Sobolev inequalities.

\noi (\ref{5.64e}). We estimate
$$\parallel <x> \nabla_K v \parallel_3\ \leq\ \parallel <x> \nabla v \parallel_3\ + \ \parallel  K \parallel_{\infty} \ \parallel <x> v \parallel_3\ ,$$
$$\parallel <x> \nabla v \parallel_3\ \leq\ \parallel <x> \nabla v \parallel_2^{1/2}\left ( \parallel  <x> \Delta v \parallel_2 \ +\ \parallel \nabla v \parallel_2\right )^{1/2} \ ,$$
$$\parallel <x>  v \parallel_3\ \leq\ \parallel <x>  v \parallel_2^{1/2}\left ( \parallel  <x> \nabla v \parallel_2 \ +\ \parallel  v \parallel_2\right )^{1/2} \ ,$$

\noi from which (\ref{5.64e}) follows.\\

\noi (\ref{5.65e}). We estimate similarly
$$\parallel \nabla_K v \parallel_{\infty}\ \leq\ \parallel \nabla v \parallel_{\infty}\ + \ \parallel  K \parallel_{\infty} \ \parallel v \parallel_{\infty}\ ,$$
$$\parallel  \nabla v \parallel_{\infty}\ \leq\ \parallel  \nabla \Delta v \parallel_2^{1/2}\ \parallel \Delta v \parallel_2^{1/2}\ ,$$
$$\parallel  v \parallel_{\infty}\ \leq\ \parallel  \nabla \Delta v \parallel_2^{1/2}\ \parallel v \parallel_2^{1/2}\ ,$$

\noi from which (\ref{5.65e}) follows.\\

\noi (\ref{5.66e}) (\ref{5.67e}). We expand and estimate
$$\Delta_K v = \Delta v - 2i K \cdot \nabla v - i (\nabla \cdot K) v - K^2 v\ ,$$
$$\left | \ \parallel <x> \Delta_K v \parallel_r\ -\ \parallel <x> \Delta v \parallel_r \right | \ \leq \ \parallel K \parallel_{\infty}\ \parallel <x> \nabla v \parallel_r$$
\beq
\label{5.70e}
+\ \parallel \nabla \cdot L \parallel_3\ \parallel <x> v \parallel_{r_1}\ + \ \left ( \parallel \nabla \cdot K_a \parallel_{\infty}\ +\ \parallel K \parallel_{\infty}^2 \right ) \parallel <x> v \parallel_r
\eeq

\noi with $2 \leq r \leq 3$ and $1/r_1 = 1/r - 1/3$. We next estimate
$$\parallel <x> v \parallel_{r_1}\ \vee \ \parallel <x> \nabla v \parallel_r\ \leq \ \parallel <x> v \parallel_r\ + \ \parallel<x> v\parallel_r^{1/2} \ \parallel <x> \Delta v \parallel_r^{1/2}$$

\noi so that (\ref{5.70e}) can be continued as 
$$\cdots \leq\ \left ( \parallel K \parallel_{\infty} \ +\ \parallel \nabla \cdot L \parallel_3 \right ) \parallel <x> v \parallel_r^{1/2} \ \parallel <x> \Delta v \parallel_r^{1/2}$$
$$+\ \left ( \parallel \nabla \cdot L \parallel_3 \ +\ \parallel \nabla \cdot K_a \parallel_{\infty} \ +\ \parallel K \parallel_{\infty}^2 \right ) \parallel <x> v \parallel_r$$

\noi from which (\ref{5.66e}) (\ref{5.67e}) follow by elementary algebraic manipulations. \\

\noi (\ref{5.68e}). We expand and estimate similarly
$$\nabla_K^2 v = \nabla^2 v - 2 i K \nabla v - i (\nabla K) v - K^2 v\ ,$$
$$\parallel \nabla_K^2 v \parallel_r\ \leq \ \parallel  \Delta v \parallel _r \ +\ \parallel K\parallel_{\infty} \ \parallel  \nabla v \parallel_r \ +\ \parallel  \nabla L \parallel_3\ \parallel  v \parallel_{r_1}$$
$$+ \left ( \parallel \nabla K_a \parallel_{\infty}\ +\ \parallel K \parallel_{\infty}^2\right ) \parallel v \parallel _r$$

\noi from which (\ref{5.68e}) follows by the same computation as before.\\

\noi (\ref{5.69e}). We expand and estimate similarly
$$\nabla_K \Delta_K v = (\nabla - i K ) \left (\Delta - 2 i K \cdot \nabla - i(\nabla \cdot K) - K^2 \right ) v$$
$$= \nabla \Delta v - i K \Delta v - 2i K \cdot \nabla^2 v - 2 i (\nabla K) \cdot \nabla v - 2 K^2\cdot \nabla v - i(\nabla \cdot K) \nabla v$$
$$- i (\nabla \nabla \cdot K) v - K(\nabla \cdot K) v - K^2 \nabla v - 2 K \cdot ( \nabla K) v + i K^3 v\ ,$$
$$\parallel \nabla_K \Delta_K v - \nabla \Delta v \parallel_2\ \leq \ 3 \parallel K \parallel_{\infty}\ \parallel \Delta v \parallel_2\ + \ 3 \parallel (\nabla K) \nabla v \parallel_2$$
$$+\ 3 \parallel K\parallel_{\infty}^2\ \parallel \nabla v \parallel_2\ +\ \parallel (\nabla \nabla \cdot K)v\parallel_2\ +\ 3 \parallel K \parallel_{\infty} \ \parallel (\nabla K)v \parallel_2\ +\ \parallel K \parallel_{\infty}^3 \ \parallel v \parallel_2$$
$$\leq \ \left ( \parallel K \parallel_{\infty}\ + \ \parallel \nabla L \parallel_3 \right ) \parallel \Delta v \parallel_2 \ +\ \parallel K \parallel_{\infty} \left ( \parallel K \parallel_{\infty} \ +\ \parallel \nabla L \parallel_3 \right ) \parallel \nabla v \parallel_2$$
$$+ \ \parallel \nabla K_a \parallel_{\infty}\ \parallel \nabla v \parallel_2÷ +\ \parallel \nabla \nabla \cdot K_a \parallel_3\ \parallel \nabla v \parallel_2 \ +\ \parallel \nabla \nabla \cdot L \parallel_2\ \parallel v \parallel_{\infty}$$
$$+\ \parallel K \parallel_{\infty}^3\ \parallel v \parallel_2\ +\ \parallel K \parallel_{\infty}\ \parallel \nabla K_a \parallel_{\infty} \ \parallel v \parallel_2\ .$$

\noi Now,
$$\parallel \Delta v \parallel_2\ \leq \ \parallel \nabla \Delta v \parallel_2^{2/3} \ \parallel v \parallel_2^{1/3}\ ,$$
$$\parallel \nabla v \parallel_2\ \leq \ \parallel \nabla \Delta v \parallel_2^{1/3} \ \parallel v \parallel_2^{2/3}\ ,$$
$$\parallel v \parallel_{\infty}\ \leq \ \parallel \nabla \Delta v \parallel_2^{1/2} \ \parallel v \parallel_2^{1/2}\ .$$

\noi The previous estimates imply (\ref{5.69e}) by elementary algebraic manipulations. \par \nobreak \hfill $\sq$ \par

From now on we shall work with $(q, G_2)$ in a bounded set of $X(I)$ for $I = (0, \tau ]$, so that $(q, G_2)$ satisfies
\beq
\label{5.71e}
\parallel <x> q(t) \parallel_2 \ \leq Y_0\ h(t) \ ,
\eeq
\beq
\label{5.72e}
\parallel <x> \partial_t q(t) \parallel_2\ \vee \ \parallel <x> \Delta q (t) \parallel_2\ \leq Y_2\ t^{-1}\ h(t) \ ,
\eeq
\beq
\label{5.73e}
\parallel \nabla  \partial_t q(t) \parallel_2\ \vee \ \parallel \nabla \Delta q (t) \parallel_2\ \leq Y_3\ t^{-3/2}\ h(t) \ ,
\eeq
\bea
\label{5.74e}
&&\parallel \nabla  G_2(t) \parallel_2\ \vee \ t\left ( \parallel \nabla^2 G_2(t)  \parallel_2\ \vee \ \parallel \nabla \partial_t G_2(t) \parallel_2\ \vee \ \parallel \nabla {\check G}_2(t) \parallel_2\right )\nn\\
&&\vee \ t^2\left ( \parallel \nabla^2 {\check G}_2(t)  \parallel_2\ \vee \ \parallel \nabla \partial_t {\check G}_2(t) \parallel_2\right )\ 
\leq Z\ t^{1/2}\ h(t) 
\eea

\noi for some constants $Y_0$, $Y_2$, $Y_3$, $Z$ and for all $t\in I$, with $h$ defined in Section 3. Note that from
$$\parallel <x> q \parallel_2\ \leq \ \int_0^t dt' \parallel <x> \partial_t q \parallel_2$$

\noi it follows that $Y_0 \leq Y_2$. It follows also from (\ref{5.71e})-(\ref{5.73e}) that 
\beq
\label{5.75e}
\parallel <x> \nabla q(t) \parallel_2\ \vee \ t^{-1/2} \parallel <x> q (t) \parallel_2\ \leq Y_1\ t^{-1/2}\ h(t) 
\eeq

\noi for some constant $Y_1 \leq Y_2$. \par

For $(q, G_2)$ satisfying (\ref{5.71e})-(\ref{5.74e}), the estimates of Lemmas 5.1-3 imply estimates of $G_1$, $\sigma$, $G'_2$ and $H_1$ in terms of $Y_j$, $0 \leq j \leq 3$, and $Z$.\par

Furthermore, the asymptotic region conditions (\ref{5.33e}) (\ref{5.45e})-(\ref{5.47e}) can be expressed in terms of $Y_j$, $0 \leq j \leq 3$, and $Z$, and the integrals $I_j$ occurring in Lemmas 5.1-3 are convergent. We collect the estimates thereby obtained in the following lemma.\\

\noi {\bf Lemma 5.6.} {\it Let $1/4 \leq  \beta < 3/4$, let $0 < \tau \leq \tau_0$ and $I = (0, \tau ]$. Let $w_a$, $K_a$ and the remainders satisfy the assumptions (A1) (A2) (A.3). Let $(q, G_2) \in X(I)$ satisfy the conditions (\ref{5.71e})-(\ref{5.74e}), and let $\tau$ be sufficiently small so that (\ref{5.13e}) holds and that 
\beq
\label{5.76e}
(Y_2 + r_5)h \leq t^{3/4} \wedge t^{-1/4 + 2 \beta }\qquad , \quad Z\ h \leq t^{3/4}
\eeq

\noi for all $t \in I$, where $r_5 = r_2 + r_3$. Then the following estimates hold for all $t\in I$~:}
\beq
\label{5.77e}
\parallel \nabla G_1 \parallel_2\ \vee \ t \parallel \nabla {\check G}_1\parallel_2\ \leq C \left ( Y_0 + r_3 \right ) h \ ,
\eeq
\beq
\label{5.78e}
\parallel \nabla^2 G_1 \parallel_2\ \vee \ t \parallel \nabla^2 {\check G}_1\parallel_2\ \leq C \left ( Y_1 + r_3 \right ) t^{-1/2}\ h \ ,
\eeq
\beq
\label{5.79e}
\parallel \nabla \partial_t  G_1 \parallel_2\ \vee \ t \parallel \nabla \partial_t  {\check G}_1\parallel_2\ \leq C \left ( Y_2 + r_3 \right ) t^{-1}\ h \ ,
\eeq
\beq
\label{5.80e}
\parallel \nabla^k \sigma  \parallel_2\ \vee \ t \parallel \nabla^k \partial_t \sigma \parallel_2\ \leq C \left ( Y_0 + r_5 \right ) t^{-k\beta }\ h \quad \hbox{{\it for k = 0, 1}}\ ,
\eeq
\beq
\label{5.81e}
\parallel \nabla^2 \sigma  \parallel_2\ \vee \ t \parallel \nabla^2 \partial_t \sigma \parallel_2\ \leq C \left ( Y_1 + r_5 \right ) t^{-2\beta }\ h \ ,
\eeq
\beq
\label{5.82e}
\parallel \nabla G'_2 \parallel_2\ \vee \ t \parallel \nabla {\check G}'_2\parallel_2\ \leq \ \left \{ C \left ( (Y_0 + r_5)t^{1/2}(1 - \ell n\ t) + Z t \right )+ r_4 \right \}  t^{1/2}\ h \ ,
\eeq
\beq
\label{5.83e}
\parallel \nabla^2 G'_2 \parallel_2\ \vee \ t \parallel \nabla^2 {\check G}'_2\parallel_2\ \leq \ \left \{ C \left ( (Y_1 + r_5)t^{1/2}+ Zt ^2\right )+ r_4 \right \}  t^{-1/2}\ h \ ,
\eeq
\beq
\label{5.84e}
\parallel \nabla \partial_t G'_2 \parallel_2\ \vee \ t \parallel \nabla\partial_t  {\check G}'_2\parallel_2\ \leq \ \left \{ C \left ( (Y_2 + r_5)t^{1/2}(1 - \ell n\ t) + Zt \right )+ r_4 \right \}  t^{-1/2}\ h \ ,
\eeq
\beq
\label{5.85e}
\parallel <x> H_1 w_a  \parallel_2\  \leq C \left \{ ( Y_0 + r_5) (t^{-\beta} + t^{-1 + \beta }) + Z\ t^{-1/2} \right \}  h \ ,
\eeq
\beq
\label{5.86e}
\parallel <x> H_1 \partial_t w_a  \parallel_2\  \leq C \left \{ ( Y_1 + r_5) (t^{-3\beta /2} + t^{-1 + \beta /2 }) + Z\ t^{-1/2} \right \}  t^{-1/2}\ h \ ,
\eeq
\beq
\label{5.87e}
\parallel <x> ( \partial_t H_1)w_a  \parallel_2\  \leq C \left \{ ( Y_2 + r_5) (t^{-\beta} + t^{-1 + \beta }) + Z\ t^{-1/2} \right \}  t^{-1}\ h \ ,
\eeq
\beq
\label{5.88e}
\parallel \nabla_K H_1 w_a  \parallel_2\  \leq C \left \{ ( Y_1 + r_5) (t^{-2\beta} + t^{-1 }) + Z\ t^{-1/2} (1 - \ell n\ t)\right \}   h \ ,
\eeq
\beq
\label{5.89e}
\parallel \nabla_K H_1 \partial_t w_a  \parallel_2\  \leq C \left \{ ( Y_1 + r_5) (t^{-2\beta} + t^{-1 }) + Z\ t^{-1/2} (1 - \ell n\ t)\right \}  t^{-1/2}\  h \ ,
\eeq
\beq
\label{5.90e}
\parallel \nabla_K (\partial_t H_1) w_a  \parallel_2\  \leq C \left \{ ( Y_2 + r_5) (t^{-2\beta} + t^{-1 }) + Z\ t^{-1/2} (1 - \ell n\ t)\right \}  t^{-1}\  h \ .
\eeq
\vskip 5 truemm

\noi {\bf Proof.} The estimates (\ref{5.77e})-(\ref{5.90e}) are obtained by substituting the bounds (\ref{5.71e})-(\ref{5.74e}) on $(q, G_2)$ into the estimates of Lemmas 5.1-3. Substituting (\ref{5.74e}) and (\ref{5.77e})-(\ref{5.81e}) into (\ref{5.33e}) (\ref{5.45e})-(\ref{5.47e}) then shows that the latter conditions follow from (\ref{5.76e}). \par \nobreak \hfill $\sq$ \par

We now turn to the construction of solutions $(q', G'_2)$ of the linearized system (\ref{2.59e}). We consider $(q, G_2)$ belonging to a bounded set of $X((0, \tau ])$, defined by (\ref{5.71e})-(\ref{5.74e}) for some $\tau$, $0 < \tau \leq \tau_0$. We shall deal with solutions $(q', G'_2)$ of the system (\ref{2.59e}) defined in an interval $I = [t_0, \tau ] \cap (0, \tau ]$ for some $t_0$ with $0 \leq t_0 < \tau$. We shall need to estimate $(q', G'_2)$ in $X(I)$ and for that purpose we define the relevant norms
\beq
\label{5.91e}
Y'_0 = \ \mathrel{\mathop {\rm Sup}_{t\in I}}\ h(t)^{-1} \parallel <x> q'(t) \parallel_2\ ,
\eeq
\beq
\label{5.92e}
Y'_2 = \ \mathrel{\mathop {\rm Sup}_{t\in I}}\ t\ h(t)^{-1} \left ( \parallel <x> \partial_t q'(t) \parallel_2\ \vee \ \parallel <x> \Delta q'(t) \parallel_2\right ) \  ,
\eeq
\beq
\label{5.93e}
Y'_3 = \ \mathrel{\mathop {\rm Sup}_{t\in I}}\ t^{3/2} \ h(t)^{-1} \left ( \parallel \nabla  \partial_t q'(t) \parallel_2\ \vee \ \parallel \nabla  \Delta q'(t) \parallel_2\right ) \  ,
\eeq
\bea
\label{5.94e}
&&Z' = \ \mathrel{\mathop {\rm Sup}_{t\in I}}\ t^{-1/2} \ h(t)^{-1} \Big \{  \parallel \nabla  G'_2(t) \parallel_2\ \vee \ t \Big ( \parallel \nabla^2  G'_2(t) \parallel_2\ \vee \ \parallel \nabla \partial_t G'_2(t)\parallel_2\nn \\
&&\vee \ \parallel \nabla {\check G'}_2(t)\parallel_2\Big ) \ \vee \ t^2 \left ( \parallel \nabla^2 {\check G'}_2(t)\parallel_2\ \vee \ \parallel \nabla \partial_t {\check G'}_2(t)\parallel_2\right ) \Big \}\ .
\eea

\noi For technical reasons, we shall also need the following auxiliary norms~:
\beq
\label{5.95e}
Y'_1 = \ \mathrel{\mathop {\rm Sup}_{t\in I}}\ t^{1/2}\ h(t)^{-1} \parallel <x> \nabla q'(t) \parallel_2\ ,
\eeq
\beq
\label{5.96e}
\widetilde{Y}'_1 = \ \mathrel{\mathop {\rm Sup}_{t\in I}}\ t^{1/2}\ h(t)^{-1} \parallel <x> \nabla_K  q'(t) \parallel_2\ ,
\eeq
\beq
\label{5.97e}
\widetilde{Y}'_{3/2} = \ \mathrel{\mathop {\rm Sup}_{t\in I}}\ t^{3/4}\ h(t)^{-1} \parallel <x> \nabla_K  q'(t) \parallel_3\ ,
\eeq
\beq
\label{5.98e}
Y'_{3/2} = \ \mathrel{\mathop {\rm Sup}_{t\in I}}\ t^{3/4}\ h(t)^{-1} \parallel q'(t) \parallel_{\infty}\ ,
\eeq
\beq
\label{5.99e}
Y'_{2,t} = \ \mathrel{\mathop {\rm Sup}_{t\in I}}\ t\ h(t)^{-1} \parallel <x> \partial_t q'(t) \parallel_2\ ,
\eeq
\beq
\label{5.100e}
\widetilde{Y}'_{2,x} = \ \mathrel{\mathop {\rm Sup}_{t\in I}}\ t\ h(t)^{-1} \parallel <x> \Delta_K  q'(t) \parallel_2\ ,
\eeq
\beq
\label{5.101e}
\widetilde{Y}'_{2} = Y'_{2,t}\vee \widetilde{Y}'_{2,x} \ ,
\eeq
\beq
\label{5.102e}
\widetilde{Y}''_{2} = Y'_{2, t} \vee \ \mathrel{\mathop {\rm Sup}_{t\in I}}\ t\ h(t)^{-1} \parallel  \nabla_K^2  q'(t) \parallel_2\ ,
\eeq
\beq
\label{5.103e}
\widetilde{Y}'_{5/2} =  \ \mathrel{\mathop {\rm Sup}_{t\in I}}\ t^{5/4} \ h(t)^{-1} \left ( \parallel  \partial_t q'  \parallel_3\ \vee \  \parallel\nabla_K  q' \parallel_{\infty}\ \vee\  \parallel \nabla_K^2 q'  \parallel_3\right )\ ,
\eeq
\beq
\label{5.104e}
\widetilde{Y}'_{3,t} = \ \mathrel{\mathop {\rm Sup}_{t\in I}}\ t^{3/2}\ h(t)^{-1} \parallel \nabla_K \partial_t q'(t)  \parallel_2\  ,
\eeq
\beq
\label{5.105e}
\widetilde{Y}'_{3,x} = \ \mathrel{\mathop {\rm Sup}_{t\in I}}\ t^{3/2}\ h(t)^{-1} \parallel \nabla_K \Delta_Kq'(t)  \parallel_2\  ,
\eeq
\beq
\label{5.106e}
\widetilde{Y}'_{3} = \widetilde{Y}'_{3,t}\vee \widetilde{Y}'_{3,x} \ .
\eeq

\noi We shall also need the following norms of $\widetilde{R}_1 = R_1 - H_1 w_a$~:
\beq
\label{5.107e}
N_0 = \ \mathrel{\mathop {\rm Sup}_{t\in (0, \tau ]}}\ h(t)^{-1} \int_0^t dt' \parallel <x> \widetilde{R}_1(t') \parallel_2\ ,
\eeq
\beq
\label{5.108e}
N_2 = \ \mathrel{\mathop {\rm Sup}_{t\in (0, \tau ]}}\ t\ h(t)^{-1} \int_0^t dt' \parallel <x> \partial_t \widetilde{R}_1(t') \parallel_2\ ,
\eeq
\beq
\label{5.109e}
N_3 = \ \mathrel{\mathop {\rm Sup}_{t\in (0, \tau ]}}\ t^{3/2}\ h(t)^{-1} \left \{ \parallel  \nabla_K \widetilde{R}_1(t)   \parallel_2\ \vee \ \int_0^t dt'  \parallel\nabla_K  \partial_t \widetilde{R}_1(t')  \parallel_2\right \}  
\eeq

\noi The $N_i$, $i = 0, 2, 3$ are finite and are estimated in the following lemma.\\

\noi {\bf Lemma 5.7.} {\it Let the assumptions of Lemma 5.6 be satisfied. Then the following estimates hold~:}
\beq
\label{5.110e}
N_0 \leq N_2 \leq r_1 + C \left \{ (Y_2 + r_5) (\tau^{1-\beta} + \tau^{\beta} ) + Z\ \tau^{1/2}\right \}\ ,
\eeq
\beq
\label{5.111e}
N_3 \leq r_1 \left ( 1 + C\ \tau^{1/2}(1 - \ell n \ \tau)\right )  + C \left \{ (Y_2 + r_5) (\tau^{3/2-2\beta} + \tau^{1/2} ) + Z\ \tau (1 - \ell n\ \tau )\right \}\ .
\eeq

\noi {\bf Proof.} The contribution of $H_1 w_a$ is estimated by Lemma 5.6, especially (\ref{5.85e})-(\ref{5.90e}), under the condition $\beta \leq 3/4$. The contribution of $R_1$ is estimated by the assumption (A3), especially (\ref{5.8e}) (\ref{5.9e}), except for the contribution of $K$ to $\nabla_K \partial_t R_1$ in $N_3$. That contribution is estimated by
$$\int_0^t dt' \parallel K \partial_t R_1(t') \parallel_2\ \leq C \int_0^t dt' (1 - \ell n\ t') \parallel \partial_t R_1(t')\parallel_2$$
$$\leq C \sum_{j=0}^{\infty} \int_{Ij} dt' (1 - \ell n \ t') \parallel \partial_t R_1(t') \parallel_2$$

\noi where $I_j = [t 2^{-(j+1)} , t2^{-j}]$,
\begin{eqnarray*}
&&\leq C\ r_1 \sum_j \left (1 - \ell n\ t + (j+1) \ell n\ 2\right )  t^{-1} \ 2^j \ h(t\ 2^{-j})\\
&&\leq C\ r_1 \sum_j \left (1 - \ell n\ t + (j+1) \ell n\ 2\right )  t^{1/2} \ 2^{-j/2} \ \overline{h}(t\ 2^{-j})\\
&&\leq C\ r_1 \ t^{-1} \ h(t) \sum_j \left (1 - \ell n\ t + (j+1) \ell n\ 2\right ) 2^{-j /2}\\
&&\leq C\ r_1(1 - \ell n\ t) t^{-1} \ h(t)
\end{eqnarray*}

\noi which completes the proof of (\ref{5.111e}).\par \nobreak \hfill $\sq$\par

We can now state the existence result of solutions of the linearized system (\ref{2.59e}).\\

\noi {\bf Proposition 5.1.} {\it Let $1/4 \leq \beta < 3/4$, let $0 < \tau \leq \tau_0$ and $I = (0, \tau ]$. Let $w_a$, $K_a$ and the remainders satisfy the assumptions (A1) (A2) (A3) and let $B_0$ satisfy (\ref{3.32e}) for $0 \leq j, k \leq 1$ and $r = \infty$. Let $(q, G_2) \in X(I)$, satisfying the bounds (\ref{5.71e})-(\ref{5.74e}). Then, for $\tau$ sufficiently small, there exists a unique solution $(q', G'_2)$ of the system (\ref{2.59e}) in $X(I)$, and that solution is estimated in the norms $Y'_0$, $Y'_2$, $Y'_3$, $Z'$ defined by (\ref{5.91e})-(\ref{5.94e}) by }
\bea
\label{5.112e}
&&Y'_0 \leq N_0 + N_2\tau + N_3 \tau^{3/2}\ , \\
\label{5.113e}
&&Y'_2 \leq N_2 + N_3\tau^{1/2}  + C\  N_0 \tau^{1/2}\ , \\
\label{5.114e}
&&Y'_3 \leq N_3 + b_0 (N_0 N_2)^{1/2}  + C\  N_0 \tau^{1/4} + C\ N_2 \tau^{1/2} (1 - \ell n\ \tau)\ , \\
&&Z' \leq C \left ( (Y_2 + r_5) \tau^{1/2} (1 - \ell n\ \tau ) + Z \tau \right ) + r_4\ .
\label{5.115e}
\eea
\vskip 5 truemm

\noi {\bf Proof.} The existence of $G'_2$ defined in the same interval as $(q, G_2)$ and the estimate (\ref{5.115e}) follow from the explicit expression in  (\ref{2.59e}), from the estimates (\ref{5.82e})-(\ref{5.84e}) of Lemma 5.6 and from the definition (\ref{5.94e}) of $Z'$.\par

We now consider the case of $q'$. We first take $\tau$ sufficiently small to satisfy the conditions (\ref{5.13e}) and (\ref{5.76e}) of Lemma 5.6. Let $0 < t_0 < \tau$ and let $q'_{t_0}$ be the solution of the Schr\"odinger equation in (\ref{2.59e}) with initial condition $q'_{t_0}(t_0) = 0$ obtained by Proposition 4.1. We shall construct $q'$ as the limit of $q'_{t_0}$ as $t_0$ tends to zero and for that purpose we need estimates of $q'_{t_0}$ in $X(I)$ with $I = [t_0, \tau ]$ that are uniform in $t_0$. Those estimates make an essential use of Lemmas 5.4 and 5.5. In all the computation we omit the subscript $t_0$ for brevity and we use the definitions (\ref{5.91e})-(\ref{5.106e}) with $I = [t_0, \tau ]$.\par

We integrate successively (\ref{5.48e})-(\ref{5.52e}) in $[t_0, t]$. Integrating (\ref{5.48e}) (\ref{5.49e}) yields 
$$\parallel <x> q'(t) \parallel_2 \ \leq \widetilde{Y}'_1 \int_{t_0}^t dt'\ t{'}^{-1/2} \ h(t') + \int_{t_0}^t dt' \parallel <x>  \widetilde{R}_1(t') \parallel_2$$

\noi so that
\beq
\label{5.116e}
Y'_0 \leq \widetilde{Y}'_1\ \tau^{1/2} + N_0 \ .
\eeq

\noi Integrating (\ref{5.50e}) (\ref{5.51e}) with $\partial_tq'(t_0) = i \widetilde{R}_1(t_0)$ yields
\begin{eqnarray*}
\parallel <x> \partial_t q'(t) \parallel_2 &\leq& \left ( \widetilde{Y}'_{3,t} + C (\widetilde{Y}'_1 + \widetilde{Y}'_{3/2} + Y'_0) \right ) \int_{t_0}^t dt'\ t{'}^{-3/2} \ h(t')\\ 
&+& \parallel <x>  \widetilde{R}_1(t_0) \parallel_2\ + \ \int_{t_0}^t dt' \parallel <x>  \partial_t \widetilde{R}_1(t') \parallel_2
\end{eqnarray*} 

\noi so that
$$Y'_{2,t} \leq \widetilde{Y}'_{3,t} \ \tau^{1/2} + C \left ( \widetilde{Y}'_1 + \widetilde{Y}'_{3/2} + Y'_0 \right )  \tau^{1/2} + N_2\ .$$

\noi Finally integrating (\ref{5.52e}) yields 
\begin{eqnarray*}
&&\parallel \nabla_K  \partial_t q'(t) \parallel_2 \ \leq b_0\ \widetilde{Y}'_1 \int_{t_0}^t dt' \ t{'}^{-5/2}\ h(t')\\
&& + C \left ( \widetilde{Y}''_2 + \widetilde{Y}'_{5/2} +  \widetilde{Y}'_1 + Y'_{3/2} + Y'_0 \right ) \int_{t_0}^t dt' \ t{'}^{-2}\ h(t')\\
&&+\ \parallel \nabla_K   \widetilde{R}_1(t_0) \parallel_2\ + \ \int_{t_0}^t dt' \ \parallel \nabla_K   \partial_t \widetilde{R}_1(t') \parallel_2
\end{eqnarray*} 

\noi so that 
$$\widetilde{Y}'_{3,t} \leq b_0 \ \widetilde{Y}'_1 + C \left ( \widetilde{Y}''_2 + \widetilde{Y}'_{5/2} + \widetilde{Y}'_1 + Y'_{3/2} + Y'_0 \right ) \tau^{1/2} + N_3\ .$$ 

\noi From (\ref{5.53e}) (\ref{5.54e}), we then obtain
\begin{eqnarray*}
&&\widetilde{Y}'_{2,x} \leq Y'_{2,t} + C\ Y'_0 \tau^{1/2} + N_0\ , \\
&&\widetilde{Y}'_{3,x} \leq \widetilde{Y}'_{3,t} + C \left ( \widetilde{Y}'_1 + Y'_0 \right ) \tau^{1/2} + N_3\ ,
\end{eqnarray*} 

\noi so that $\widetilde{Y}'_2$ and $\widetilde{Y}'_3$ satisfy the same estimates as $Y'_{2,t}$, $\widetilde{Y}'_{3,t}$, namely 
\bea
\label{5.117e}
&&\widetilde{Y}'_2 \leq \widetilde{Y}'_{3,t}\  \tau^{1/2} + C \left ( \widetilde{Y}'_1 + \widetilde{Y}'_{3/2} + Y'_0\right ) \tau^{1/2} + N_2\ , \\
&&\widetilde{Y}'_3 = b_0\ \widetilde{Y}'_1 + C  \left ( \widetilde{Y}''_2 + \widetilde{Y}'_{5/2} +\widetilde{Y}'_1 +  Y'_{3/2} + Y'_0 \right ) \tau^{1/2} + N_3\ .
\label{5.118e}
\eea

 The next step consists in eliminating the covariant derivatives in (\ref{5.116e}) (\ref{5.117e}) (\ref{5.118e}) by using the assumption (A2), Sobolev inequalities and Lemma 5.5. We first estimate
 \beq
 \label{5.119e}
 \widetilde{Y}'_1  \leq Y'_1 + C\ Y'_0 \ \tau^{1/2} (1 - \ell n\ \tau )
 \eeq

\noi so that (\ref{5.116e}) implies
 \beq
 \label{5.120e}
Y'_0 \leq Y'_1 \ \tau^{1/2}+ C\ Y'_0 \ \tau (1 - \ell n\ \tau ) + N_0\ .
 \eeq

 \noi We next consider (\ref{5.117e}). We estimate
 \beq
 \label{5.121e}
 \widetilde{Y}'_{3,t}   \leq Y'_3 + C\ Y'_2 \ \tau^{1/2} (1 - \ell n\ \tau ) \ ,
 \eeq
 \bea
  \label{5.122e}
 \widetilde{Y}'_{3/2}  & \leq& \left ( Y'_1 \left ( Y'_2 + Y'_1 \ \tau^{1/2}\right ) \right )^{1/2} + C \left ( Y'_0 \left ( Y'_1 + Y'_0\ \tau^{1/2}\right ) \right )^{1/2} \ \tau^{1/2} (1 - \ell n\ \tau)\nn \\
 &\leq& \left ( Y'_1 Y'_2\right )^{1/2} + C \left ( Y'_1 + Y'_0 \ \tau^{1/2}\right ) \tau^{1/4} (1 - \ell n\ \tau ) \ ,
 \eea
 
 \noi while
  \beq
 \label{5.123e}
Y'_2 \leq   \widetilde{Y}'_2 +  C\ Y'_0 \ \tau (1 - \ell n\ \tau )^2
 \eeq
 
\noi by (\ref{5.67e}) with $r= 2$ and $m \leq C(1 - \ell n\ \tau)^2$. Substituting (\ref{5.119e}) (\ref{5.121e}) (\ref{5.122e}) into (\ref{5.117e}) and substituting the result into (\ref{5.123e}) yields
 \beq
 \label{5.124e}
Y'_2 \leq   Y'_3  \ \tau^{1/2} +  C\left (  Y'_2  \ \tau (1 - \ell n\ \tau ) + \left ( \left ( Y'_1 Y'_2 \right )^{1/2} + Y'_1 + Y'_0 \right ) \tau^{1/2} \right ) + N_2 \ .
 \eeq

\noi We next consider (\ref{5.118e}). We estimate successively
 \beq
 \label{5.125e}
\widetilde{Y}''_2 \leq   Y'_2  +  C\ Y'_0 
 \eeq

\noi by (\ref{5.68e}) with
\beq
\label{5.126e}
\overline{m} \leq b_0\ t^{-1} + C\ t^{-1/2} \leq C\ t^{-1}\ ,
\eeq
\begin{eqnarray*}
&&\parallel \partial_t q' \parallel_3 \ \leq \left ( Y'_2 Y'_3\right )^{1/2}\ t^{-5/4}\ h(t)\ ,\\
&&\parallel \nabla_K q' \parallel_{\infty}  \ \leq \left ( \left ( Y'_2 Y'_3\right )^{1/2}+ C\left ( Y'_1 + Y'_2 \right ) t^{1/2}(1 - \ell n\ t) \right ) t^{-5/4} \ h(t)\ ,\\
&&\parallel \nabla_K^2 q' \parallel_3  \ \leq \left ( \left ( Y'_2 Y'_3\right )^{1/2}+ C\left ( Y'_1 + Y'_0 \right )\right )  t^{-5/4} \ h(t)\ ,
\end{eqnarray*}

\noi by (\ref{5.68e}) with $r= 3$ and (\ref{5.126e}), so that 
\beq
 \label{5.127e}
\widetilde{Y}'_{5/2}  \leq   \left ( Y'_2    Y'_3\right )^{1/2} + C \left ( Y'_2 \  \tau^{1/2}(1 - \ell n\ \tau) + Y'_1 + Y'_0\right ) 
 \eeq

\noi and finally 
\beq
 \label{5.128e}
Y'_{3/2}  \leq   \left ( Y'_1    Y'_2\right )^{1/2} \leq Y'_1 + Y'_2\ .
 \eeq

\noi Substituting (\ref{5.119e}) (\ref{5.125e}) (\ref{5.127e}) (\ref{5.128e}) into (\ref{5.118e}) yields
\beq
 \label{5.129e}
\widetilde{Y}'_3  \leq   b_0\ Y'_1 + C \left ( \left ( Y'_2    Y'_3\right )^{1/2} + Y'_2 + Y'_1 + Y'_0(1 - \ell n\ \tau )\right )  \tau^{1/2} + N_3\ .
 \eeq

\noi Using (\ref{5.69e}) (\ref{5.126e}) (A2) (\ref{5.76e}) we next estimate
$$Y'_3 \leq \widetilde{Y}'_3 + C\left ( Y'_2  +  Y'_0\right ) \tau^{1/2}(1 - \ell n\ \tau) + \left ( b_0 + C\ \tau^{1/2}\right ) Y'_1$$

\noi which together with (\ref{5.129e}) yields
\beq
 \label{5.130e}
Y'_3  \leq   b_0\ Y'_1 + C \left ( \left ( Y'_2    Y'_3\right )^{1/2} + \left ( Y'_2 + Y'_0\right )  (1 - \ell n\ \tau ) + Y'_1 \right )  \tau^{1/2} + N_3\ .
 \eeq

\noi We next simplify the resulting inequalities (\ref{5.120e}) (\ref{5.124e}) (\ref{5.130e}) by using the inequality $Y'_1 \leq (Y'_0 Y'_2)^{1/2}$ to eliminate $Y'_1$, the obvious inequality $Y'_0 \leq Y'_2$ at some harmless places, smallness conditions of the type $C\tau (1 - \ell n\ \tau)\leq 1$ to eliminate the diagonal terms in $Y'_0$ and $Y'_2$ in  (\ref{5.120e}) (\ref{5.124e}), and some elementary algebraic manipulations. We obtain 
\bea
 \label{5.131e}
 &&Y'_0 \leq Y'_2\ \tau + N_0\ ,\\
 \label{5.132e}
 &&Y'_2 \leq Y'_3 \ \tau^{1/2} + C Y'_0 \ \tau^{1/2} + N_2\ , \\
\label{5.133e}
&&Y'_3 \leq b_0(Y'_0 Y'_2)^{1/2} + C Y'_2  \ \tau^{1/2} (1 - \ell n\ \tau ) + N_3\ .
\eea

\noi Substituting (\ref{5.131e}) (\ref{5.133e}) into (\ref{5.132e}) yields
\begin{eqnarray*}
Y'_2 &\leq & b_0 \left ( N_0 Y'_2\right )^{1/2} \tau^{1/2} + C Y'_2 \tau  (1 - \ell n\ \tau )  + N_3 \ \tau^{1/2} + C\ N_0\ \tau^{1/2}\\
&&+ C\ Y'_2\ \tau^{3/2} + N_2
\end{eqnarray*}

\noi which yields (\ref{5.113e}) under an additional smallness condition on $\tau$. Substituting (\ref{5.113e}) into (\ref{5.131e}) (\ref{5.133e}) yields (\ref{5.112e}) (\ref{5.114e}).\par

We have derived the estimates (\ref{5.112e}) (\ref{5.113e}) (\ref{5.114e}) for the solution $q'_{t_0}$ in the interval $[t_0, \tau ]$. We now take the limit of that solution when $t_0$ tends to zero. Let $0 < t_0 \leq t_1 \leq \tau$. From the conservation of the $L^2$ norm for the difference of two solutions, it follows that 
\beq
\label{5.134e}
\parallel q'_{t_0}(t) - q'_{t_1}(t) \parallel _2\ = \ \parallel q'_{t_0}(t_1)\parallel _2\ \leq C\ h(t_1)
\eeq

\noi where $C$ is the right hand side of (\ref{5.112e}). Let now $I = (0 , \tau ]$. It follows from (\ref{5.134e}) that $q'_{t_0}$ converges in $L_{loc}^{\infty} (I, L^2)$ to a limit $q' \in {\cal C} (I, L^2)$. From that convergence and from the uniform estimates (\ref{5.112e})-(\ref{5.114e}) it follows that $q' \in {\cal C}(I, H^k) \cap (L^{\infty} \cap {\cal C}_w) (I, H^3)$ for $0 \leq k < 3$, and that $q'_{t_0}$ converges to $q'$ in $L_{loc}^{\infty} (I, H^k)$ norm and weakly in $H^3$ pointwise in time. Similar properties hold for $xq'$. From the Schr\"odinger equation in (\ref{2.59e}) and the previous convergence, it follows that $\partial_t q'_{t_0}$ converges to $\partial_t q'$ in $L_{loc}^{\infty} (I, L^2)$. From that convergence and from the uniform bounds  (\ref{5.113e}) (\ref{5.114e})  it follows that $\partial_t q' \in {\cal C} (I, H^k) \cap (L^{\infty} \cap {\cal C}_w)(I, H^1)$ for $0 \leq k' < 1$ and that $\partial_t q'_{t_0}$ converges to $\partial_t q'$ in $L_{loc}^{\infty} (I,  H^k)$ norm and weakly in $H^1$ pointwise in time. Similar properties hold for $x \partial_t q'$. From the previous convergences and from the uniform estimates (\ref{5.112e})-(\ref{5.114e}), it follows that $q'$ satisfies the same estimates in I. Clearly $q'$ is a solution of the Schr\"odinger equation in (\ref{2.59e}). Finally by the uniqueness part of Proposition 4.1, $q'$ coincides with the solution obtained from that proposition with initial data $q'(\overline{t})$ for some (any) $\overline{t} \in I$, and therefore enjoys all the properties obtained in that proposition, so that $(q', 0) \in X(I)$.\par \nobreak \hfill $\sq$ \par

We can now derive the main result of this section, namely the existence of solutions of the non linear system (\ref{2.54e}).\\

\noi {\bf Proposition 5.2.} {\it Let $1/4 \leq \beta < 3/4$. Let $w_a$, $K_a$ and the remainders satisfy the assumptions (A1) (A2) (A3) and let $B_0$ satisfy (\ref{3.32e}) for $0 \leq j, k \leq 1$ and $r = \infty$. Then there exists $\tau$, $0 < \tau \leq \tau_0$, and there exists a unique solution $(q, G_2) \in X((0, \tau ])$ of the system (\ref{2.54e}). In particular $(q, G_2)$ satisfies the estimates (\ref{5.71e})-(\ref{5.74e}) for some constants $Y_0$, $Y_2$, $Y_3$ and $Z$ depending on $\beta$, on $w_a$, $K_a$ and on the remainders through the norms occuring in the assumptions (A1) (A2) (A3). The solution $(q, G_2)$ is unique under the assumption that $(q, G_2) \in X_0(I)$ where $I = (0, \tau ]$, and that $(q, G_2)$ satisfy the following conditions~:
$$q \in L^{\infty} (I, H^3)\ ,\ xq \in L^{\infty} (I, H^2)\ ,\ t \partial_t q \in L^{\infty} (I, H^1)\ ,$$
$$\nabla G_2 \in L^{\infty} (I, H^1)\ ,\ t \nabla \partial_t G_2 \in L^{\infty} (I, L^2)\ ,\ t \nabla^2 {\check G}_2 \in L^{\infty} (I, L^2)\ ,$$
$$\parallel <x> q(t) \parallel_2\ \leq C\ h_1(t)$$ 

\noi for all $t \in I$, for some $h_1$ satisfying the conditions of Proposition 4.2.}\\

\noi {\bf Proof.} Let $0 < \tau \leq \tau_0$ and $I = (0 , \tau ]$. For $\tau$ sufficiently small, Proposition 5.1 defines a map $\Gamma : (q, G_2) \to (q', G'_2)$ from $X(I)$ into itself. We now show that for $\tau$ sufficiently small, the map $\Gamma$ is a contraction on the subset ${\cal R}$ of $X(I)$ defined by (\ref{5.71e})-(\ref{5.74e}) for a suitable choice of $Y_0$, $Y_2$, $Y_3$ and $Z$, in the norms associated with Lemma 4.1. \par

We first ensure that ${\cal R}$ is stable under $\Gamma$. From Proposition 5.1, especially (\ref{5.112e})-(\ref{5.115e}) and from Lemma 5.7, especially  (\ref{5.110e}) (\ref{5.111e}), it follows that $Y'_0$, $Y'_2$, $Y'_3$ and $Z'$ defined by (\ref{5.91e})-(\ref{5.94e}) satisfy 
\bea
\label{5.135e}
&&Y'_0 \leq c\ r_1 + o(\tau ; Y_2, Z)\\
\label{5.136e}
&&Y'_2 \leq c\ r_1 + o(\tau ; Y_2, Z)\\
\label{5.137e}
&&Y'_3 \leq c\ r_1 (1 + b_0) + o(\tau ; Y_2, Z)\\
&&Z' \leq r_4 + o(\tau ; Y_2, Z) 
\label{5.138e}
\eea 

\noi where $c$ is an absolute constant which was omitted in the previous estimates and which we reintroduce for the sake of the present argument, and where $o(\tau ; Y_2, Z)$ denotes a function of $\tau$, $Y_2$, $Z$ increasing in $Y_2$, $Z$ and tending to zero as a power of $\tau$ when $\tau$ tends to zero. We now choose
\beq
\label{5.139e}
Y_0 = Y_2 = 2c\ r_1 \quad , \quad Y_3 = 2c\ r_1 (1 + b_0)\quad , \quad Z = 2r_4
\eeq

\noi and we take $\tau$ sufficiently small so that $o(\tau ; 2cr_1, 2r_4) \leq cr_1$ in (\ref{5.135e})-(\ref{5.137e}) and $o(\tau ; 2cr_1, 2r_4) \leq r_4$ in (\ref{5.138e}). This ensures that ${\cal R}$ is stable under $\Gamma$.

We next show that $\Gamma$ is a contraction on ${\cal R}$ for the norms occurring in Lemma 4.1, namely the $L^2$ norms of $<x>q$, $\nabla q$, $\nabla G_2$ and $\nabla {\check G}_2$. Let $(q_i, G_{2i}) \in {\cal R}$, $i = 1,2$, and let $(q'_i, G'_{2i})$, $i = 1,2$, be their images under $\Gamma$. We define $(q_{\pm}, G_{2\pm})$ and $(q'_{\pm}, G'_{2\pm})$ by $f_{\pm} = (1/2) (f_1 \pm f_2)$, so that in particular all those quantities belong to ${\cal R}$. We define the norms
\bea
\label{5.140e}
&&Y_{0-} = \ \mathrel{\mathop {\rm Sup}_{t\in I}}\ h(t)^{-1} \parallel <x> q_- (t) \parallel_2\ , \\
\label{5.141e}
&&Y_{1-} = \ \mathrel{\mathop {\rm Sup}_{t\in I}}\ t^{1/2}\ h(t)^{-1} \parallel \nabla q_- (t) \parallel_2\ , \\ 
&&Z_- = \ \mathrel{\mathop {\rm Sup}_{t\in I}}\ t^{-1/2}\ h(t)^{-1} \left ( \parallel \nabla G_{2-} (t) \parallel_2\ \vee\ t \parallel \nabla {\check G}_{2-} (t) \parallel_2\right )\ ,
\label{5.142e}
\eea 

\noi and similarly for the primed quantities, and we estimate $Y'_{0-}$, $Y'_{1-}$, $Z'_-$ in terms of \break\noindent $Y_{0-}$, $Y_{1-}$, $Z_-$ by Lemma 4.1 with $(w_-, B_{1-}, s_-, B_{2-}) = (q_-, G_{1-}, \sigma_-, G_{2-})$,\break\noindent $(w_+, B_{1+}, s_+, B_{2+}) = (w_a, B_{1a}, s_a, B_{2a}) + (q_+, G_{1+}, \sigma_+, G_{2+})$ and similarly for the primed quantities. Note that the proof of Lemma 4.1 has to be slightly modified for the present application since $B_{2a}$ does not satisfy the assumptions on $B_2$ made in that Lemma. It suffices to separate $B_{2+} = B_{2a} + G_{2+}$ in (\ref{4.90e}) (\ref{4.91e}) and to estimate $B_{2a}$ in the same way as $B_{1+}$. Omitting the subscript $-$ in $Y_{0-}$, $Y_{1-}$, $Z_-$ and in $Y'_{0-}$, $Y'_{1-}$, $Z'_-$ for brevity, we obtain from Lemma 4.1 
\begin{eqnarray*}
&&\parallel \nabla G_{1-}\parallel_2\ \leq C\ Y_0\ h\ , \\
&&\parallel \nabla {\check  G}_{1-}\parallel_2\ \leq C\ t^{-1}\ Y_0\ h\ , \\
&&\parallel \nabla^k \sigma_- \parallel_2\ \leq C\left (  Y_0\ t^{-k \beta} + \delta_{k,2}\ Y_1 \ t^{-1/2}\right )  h \quad \hbox{for $k = 0, 1, 2$} \ , 
\end{eqnarray*}
\beq
\label{5.143e}
\left | \partial_t \parallel q'_- \parallel_2\right | \ \leq C\left (  Y_0\left (  t^{-\beta} + t^{-1 + \beta }\right ) + Z\ t^{-1/2}\right )  h \quad  \ , \\
\eeq
\bea
\label{5.144e}
\left | \partial_t \parallel xq'_- \parallel_2\right | &\leq& \left (  Y'_1\ t^{-1/2} + C\ Y'_0 (1 - \ell n\ t) \right ) h\nn\\
&&+ C\left ( Y_0 \left (t^{-\beta } + t^{-1 + \beta } \right ) + Z\ t^{-1/2}\right  )  h \ , 
\eea
\bea
\label{5.145e}
\left | \partial_t \parallel \nabla_{K_+} q'_- \parallel_2\right | &\leq& C\Big  (  Y'_0\ t^{-1} + \left (Y'_0 Y'_1\right )^{1/2} t^{-5/4} + Y_0\left  (t^{-1} + t^{-2\beta}\right ) \nn\\
&&+ Y_1\  t^{-1/2 } + Z\  t^{-1/2} (1 - \ell n\ t) \Big )  h \ . 
\eea

\noi Integrating (\ref{5.143e})-(\ref{5.145e}) over time and estimating $G'_{2-}$ by Lemma 4.1, especially (\ref{4.80e}) (\ref{4.81e}), we obtain
\beq
\label{5.146e}
Y'_0 \leq Y'_1 \ \tau^{1/2} + C \left ( Y'_0\ \tau (1 - \ell n\ \tau) + Y_0 \left ( \tau^{1-\beta} + \tau^{\beta}\right ) + Z\ \tau^{1/2}\right )\ ,
\eeq
\bea
\label{5.147e}
Y'_1 &\leq& C \Big( Y'_0 \ \tau^{1/2} (1 - \ell n\ \tau )+  \left ( Y'_0 Y'_1 \right )^{1/2} \ \tau^{1/4} + Y_0\  \tau^{1/2} \left ( 1 +  \tau^{1 - 2\beta}\right ) \nn \\
&&+ Y_1 \ \tau + Z\ \tau  (1 - \ell n\ \tau ) \Big )\ ,
\eea
\beq
\label{5.148e}
Z' \leq C \left ( Y_0 \ \tau^{1/2} (1 - \ell n\ \tau ) + Z\  \tau \right )\ .
\eeq

\noi Eliminating $(Y'_0 Y'_1)^{1/2}$ in (\ref{5.147e}) by an elementary algebraic manipulation and taking the sum, we obtain 
$$Y'_0 + Y'_1 + Z' \leq C \left ( Y'_0 + Y'_1 \right ) \tau^{1/2} (1 - \ell n\ \tau )$$
$$+ C \left ( Y_0 \left (\tau^{1-\beta} + \tau^{\beta} + \tau^{1/2} (1 - \ell n\ \tau ) +  \tau^{3/2 - 2\beta}\right ) + Y_1\ \tau + Z\  \tau^{1/2}\right )$$

\noi which implies
$$Y'_0 + Y'_1 + Z' \leq (1/2) \left ( Y_0 + Y_1 + Z \right )$$

\noi for $\tau$ sufficiently small. This proves that $\Gamma$ is a contraction in the norms (\ref{5.140e})-(\ref{5.142e}). The existence result now follows from the fact that ${\cal R}$ is closed in those norms.\par

Finally, the uniqueness result follows from Proposition 4.2.\par \nobreak \hfill $\sq$ \par

\mysection{Asymptotic functions and remainder estimates}
\hspace*{\parindent} In this section we construct approximate solutions $(w_a, s_a, B_a)$ of the auxiliary system (\ref{2.41e}) satisfying the assumptions made in Section 5 and in particular the remainder estimates needed for the Cauchy problem at $t=0$ for that system. We construct $(w_a, s_a, B_a)$ by solving the system (\ref{2.41e}) by iteration with the contribution of $B_0$ omitted, so that $B_0$ will appear only in the remainders $R_1$ and $R_4$. Accordingly we define 
$$B'_a = B_{1a} + B_{2a}\ .$$

\noi The contribution of the $j$-th order in the iteration for $(w_a, s_a, B'_a)$ will behave as $t^j$ modulo logarithms. Since we need an accuracy of order $t^{3/2 + \lambda}$ with $\lambda > 0$ (see the definition of $h$ in Section 3), it will suffice to iterate once in order to reach that accuracy with $\lambda < 1/2$. We recall also that as regards $B_a$, we choose separately asymptotic forms $B_{1a}$ and $B_{2a}$ for $B_1$ and $B_2$, in spite of the fact that $B_1$ is an explicit function of $w$. Thus we define
\beq
\label{6.1e}
\left \{ \begin{array}{l}w_a = w_{a0} + w_{a1} \qquad , \quad s_a = s_{a0} + s_{a1}\ ,   \\ \\ B_{1a} = B_{1a0} + B_{1a1} \qquad , \quad B_{2a} = B_{2a1} + B_{2a2}\ . \end{array} \right . \eeq

\noi The shift by 1 in the last index of $B_{2a}$ reflects the fact that $B_2$ itself is already of order $t$, so that in all cases the last index indicates the power of $t$ in the asymptotic behaviour. The term $B_{2a2}$ is apparently of order $t^2$ and could probably be omitted at the cost of serious technical complications. We define the lowest order quantities by 
\beq
\label{6.2e}
\left \{ \begin{array}{l} i \partial_t w_{a0} + (1/2) \Delta w_{a0} = 0 \qquad , \quad w_{a0}(0) = w_+\ ,  \\ \\ B_{1a0} = B_1 (w_{a0}) \\ \\  \partial_t s_{a0} = t^{-1} \nabla g (w_{a0}) + \nabla {\check B}_{1L} (w_{a0}) \qquad , \quad s_{a0} (1) = 0 \ , \\ \\  B_{2a1} = {\cal B}_2(w_{a0}, w_{a0} , s_{a0} + B_{1a0} )\end{array} \right . 
\eeq

\noi and the next order quantities by 
\bea
\label{6.3e}
\left \{ \begin{array}{l}
i\partial_t w_{a1} = i\left ( s_{a0} + B_{1a0}\right ) \cdot \nabla w_{a0} + (i/2) (\nabla \cdot s_{a0}) w_{a0} +\\
\qquad   \quad (1/2) \left ( s_{a0} + B_{1a0}\right )^2 w_{a0} + \left ( {\check B}_{1a0S} +  {\check B}_{2a1} \right ) w_{a0}\qquad , \quad  w_{a1} (0) = 0\\ \\
B_{1a1} = 2B_1 (w_{a0}, w_{a1})\\ \\
\partial_t s_{a1} = 2t^{-1} \nabla g (w_{a0}, w_{a1}) + 2 \nabla {\check B}_{1L} (w_{a0} , w_{a1}) \qquad , \quad s_{a1}(0) = 0\\ \\
B_{2a2} = 2 {\cal B}_2(w_{a0}, w_{a1}, s_{a0} + B_{1a0}) - t F_2 \left ( P (s_{a1} + B_{1a1} + B_{2a1}) |w_{a0}|^2 \right ) \ . 
\end{array} \right . 
\eea

\noi For any polynomial function $f(w_a, s_a, B_{1a}, B_{2a})$ and any nonnegative integer $p$, we define
\beq
\label{6.4e}
f(w_a,s_a,B_{a1}, B_{a2})_{\geq p} = \sum_{j+k+\ell + m \geq p} f(w_{aj}, s_{ak}, B_{1a\ell }, B_{2am})\ .
\eeq

\noi The remainders $R_j$, $1 \leq j \leq 4$, defined by (\ref{2.57e}) then become
\beq
\label{6.5e}
R_1 = R_{10} + R_{11}
\eeq

\noi where
\bea
\label{6.6e}
&&R_{10} = - i B_0 \cdot \nabla w_a - \left ( B_0 (s_a + B_{1a} + B_{2a}) + (1/2) B_0^2 + {\check B}_0 \right ) w_a \ ,\\
\label{6.7e}
&&R_{11} = (1/2) \Delta w_{a1} - \Big \{ i (s_a + B_{1a} + B_{2a})\cdot \nabla w_a + (i/2) (\nabla \cdot s_a) w_a\nn\\
&&\qquad \quad+ (1/2) (s_a + B_{1a} + B_{2a})^2 w_a + {\check B}_{1aS} w_a \big \}_{\geq 1} - ({\check B}_{2a} w_a)_{\geq 2}\ , \\
\label{6.8e}&&R_2 = - t^{-1} \nabla g (w_a)_{\geq 2} = - t^{-1} \nabla g (w_{a1})\ , \\
\label{6.9e}
&&R_3 = - B_1 (w_a)_{\geq 2} = - B_1 (w_{a1})\ , \\
&&R_4 = R_{40} + R_{41}
\label{6.10e}
\eea

\noi where
\bea
\label{6.11e}
&&R_{40} = t F_2 (PB_0 |w_a|^2)\ , \\
&&R_{41} = - {\cal B}_2 (w_a, w_a, s_a + B_{1a} + B_{2a})_{\geq 2} \ .
\label{6.12e}
\eea

\noi We now turn to estimate $(w_a, s_a, B_a)$. We use the spaces $\ddot{H}^k = \dot{H}^1 \cap \dot{H}^k$ (see Section~3) and the notation $v \in (X, f)$ to mean that $v \in {\cal C}(I, X)$ with $\parallel v(t);X\parallel \leq f(t)$ for all $t \in I$, with $I = (0, \tau ]$ for some $\tau$, $0 \leq \tau \leq 1$, with $\tau = 1$ in the present case (see Section 3, especially (\ref{3.1e})).\\

\noi {\bf Lemma 6.1.} {\it Let $w_+ \in H^{k_+}$, $xw_+ \in H^{k_+-1}$ with $k_+ \geq 5 \vee (3 + \beta^{-1})$. Then the components of $w_a$, $s_a$, $B_a$ defined by (\ref{6.1e}) (\ref{6.2e}) (\ref{6.3e}) satisfy the following properties.}
\bea
\label{6.13e}
&&w_{a0} \in \left ( H^{k_+}, 1 \right ) \quad , \quad xw_{a0} \in \left ( H^{k_+-1}, 1\right ) \ ,\\
\label{6.14e}
&&\partial_t w_{a0} \in \left ( H^{k_+-2}, 1 \right ) \quad , \quad \partial_t xw_{a0} \in \left ( H^{k_+ - 3}, 1\right ) \ ,\\
\label{6.15e}
&&B_{1a0} \in \left ( \ddot{H}^{k_++1}, 1 \right )  \quad , \quad {\check B}_{1a0} \in \left ( <x> \dot{H}^1 \cap \dot{H}^2 \cap \dot{H}^{k_+} , t^{-1} \right ) \ ,\\
\label{6.16e}
&&\partial_t B_{1a0} \in \left ( H^{k_+-1}, 1 \right )  \ , \\
\label{6.17e}
&&\partial_t {\check B}_{1a0} \in \left ( <x> \dot{H}^1 \cap \dot{H}^2 \cap \dot{H}^{k_+} , t^{-2} \right ) + \left ( \ddot{H}^{k_+-2} , t^{-1} \right )\ ,  \\
\label{6.18e}
&&s_0 \in  \left ( \ddot{H}^{k_+-1}, 1 - \ell n\ t \right )  \quad , \quad \partial_ts_0 \in  \left ( \ddot{H}^{k_+-1}, t^{-1}\right ) \ , \\
\label{6.19e}
&&B_{2a1} \in \left ( H^{k_++1}, t(1 - \ell n \ t ) \right ) \quad , \quad {\check B}_{2a1} \in \left ( \ddot{H}^{k_++1}, 1 - \ell n\ t \right ) \ ,\\
\label{6.20e}
&&\partial_t B_{2a1} \in \left ( H^{k_++1}, 1 - \ell n \ t  \right ) + \left ( H^{k_+-1}, t (1 - \ell n \ t)  \right )\ , \\
&&\partial_t {\check B}_{2a1} \in \left ( \ddot{H}^{k_++1}, t^{-1} \right ) + \left ( H^{k_+-1}, (1 - \ell n\ t) \right )\ . 
\label{6.21e}
\eea

\noi {\it Let in addition $k_1 = (k_+ - 2) \wedge (k_+ - \beta^{-1})$. Then}
\bea
\label{6.22e}
&&w_{a1} \in \left ( H^{k_1}, t(1 - \ell n \ t )^2 \right ) \quad , \quad xw_{a1} \in \left ( H^{k_1}, t(1 - \ell n \ t )^2\right ) \ ,\\
\label{6.23e}
&&\partial_t w_{a1} \in \left ( H^{k_1}, (1 - \ell n \ t )^2 \right ) \quad , \quad \partial_t xw_{a1} \in \left ( H^{k_1}, (1 - \ell n \ t )^2\right ) \ ,\\
\label{6.24e}
&&B_{1a1} \in \left ( H^{k_1+1}, t(1 - \ell n \ t )^2 \right )  \quad , \quad {\check B}_{1a1} \in \left (  \ddot{H}^{k_1+1} , (1 - \ell n \ t )^2 \right ) \ ,\\
\label{6.25e}
&&\partial_t B_{1a1} \in \left ( H^{k_1+1}, (1 - \ell n \ t )^2 \right )  \ , \\
\label{6.26e}
&&\partial_t {\check B}_{1a1} \in \left ( \ddot{H}^{k_1 +1},  t^{-1} (1 - \ell n \ t )^2\right ) + \left ( H^{k_+-2} , (1 - \ell n \ t )^2 \right ) \ , \\
\label{6.27e}
&&s_1 \in  \left ( H^{k_1}, t(1 - \ell n\ t)^2 \right )  \quad , \quad \partial_ts_1 \in  \left ( H^{k_1}, (1 - \ell n \ t )^2\right ) \ , \\
\label{6.28e}
&&B_{2a2} \in \left ( H^{k_1+1}, t^2(1 - \ell n \ t )^3 \right ) \quad , \quad {\check B}_{2a2} \in \left ( H^{k_1+1}, t(1 - \ell n\ t )^3\right ) \ ,\\
\label{6.29e}
&&\partial_t B_{2a2} \in \left ( H^{k_1+1}, t(1 - \ell n \ t )^3 \right ) + \left ( H^{k_+-2}, t^2(1 - \ell n \ t)^2  \right )\ , \\
&&\partial_t {\check B}_{2a2} \in \left ( H^{k_1+1}, (1 - \ell n \ t )^3\right ) + \left ( H^{k_+-2}, t(1 - \ell n\ t)^2 \right ) \ .
\label{6.30e}
\eea

\vskip 5 truemm

\noi {\bf Proof.} In all the proof, the subscript $a$ will be omitted. The proof uses systematically and without mention the H\"older inequality and the Sobolev inequalities of Lemma 3.1, parts 1 and 2, as well as some trivial commutator identities. In the estimates of the $B_{jk}$'s, it also uses Lemma 3.6, in particular (\ref{3.18e}) (\ref{3.19e}) and the special case $r= 4$ of (\ref{3.20e}), namely
$$\parallel F_j(M)\parallel_4\ \leq \int_1^{\infty} d\nu (\nu - 1)^{-1/2} \nu^{-j+1/4} \parallel M(t/\nu )\parallel_{4/3}$$

\noi which ensures that
\beq
\label{6.31e}
\parallel F_j (M) \parallel_4 \ \leq C\ t^k (1 - \ell n\ t)^{\ell}
\eeq

\noi if
$$\parallel M(t) \parallel_{4/3}  \ \leq C\ t^k (1 - \ell n\ t)^{\ell}$$

\noi provided $j+k > 3/4$, and in particular for $j + k \geq 1$. The estimate (\ref{6.31e}) is used only to show that $F_j(M)$ tends to zero at infinity in space so that an estimate of $\parallel \nabla F_j M\parallel_2$ suffices to prove that $F_j(M) \in \dot{H}^1$. The proof also occasionally uses Lemmas 3.2-4, which will be quoted when appropriate. The $B_{jk}$'s defined by (\ref{6.15e})-(\ref{6.30e}) are given more explicitly by 
\bea
\label{6.32e}
&&B_{10} = - F_1 (Px|w_0|^2)\\
\label{6.33e}
&&{\check B}_{10}  = - t^{-1} F_0 (x \cdot Px |w_0|^2)\\
\label{6.34e}
&&\partial_t  B_{10} = - F_2 (Px\partial_t |w_0|^2)\\
\label{6.35e}
&&\partial_t {\check B}_{10}  = - t^{-1} {\check B}_{10}  - t^{-1} F_1(x\cdot Px \partial_t |w_0|^2)\\
\label{6.36e}
&&B_{21} = t\ F_2 (PN_1) \\
\label{6.37e}
&&{\check B}_{21} = F_1 (x \cdot PN_1) \\ 
\label{6.38e}
&&\partial_t B_{21} = t^{-1} B_{21} + t F_3 (P \partial_t N_1) \\
&&\partial_t {\check B}_{21} = F_2 (x \cdot P \partial_t N_1) 
\label{6.39e}
\eea

\noi where
\bea
\label{6.40e}
&&N_1 = {\rm Im} \ \overline{w}_0 \nabla w_0 - (s_0 + B_{10}) |w_0|^2 \ , \\
\label{6.41e}
&&B_{11} = - F_1 (Px \ 2 {\rm Re} \ \overline{w}_0 w_1)\\
\label{6.42e}
&&{\check B}_{11} = - t^{-1} F_0 (x \cdot Px \ 2 {\rm Re}\ \overline{w}_0w_1)\\
\label{6.43e}
&&\partial_t B_{11} = - F_2 (Px \partial_t \ 2 {\rm Re} \ \overline{w}_0w_1)\\
\label{6.44e}
&&\partial_t {\check B}_{11} = - t^{-1} {\check B}_{11} - t^{-1} F_1 (x \cdot Px \partial_t \ 2 {\rm Re}\  \overline{w}_0w_1)\\
\label{6.45e}
&&B_{22} = t\ F_2 (PN_2) \\
\label{6.46e}
&&{\check B}_{22} = F_1 (x \cdot PN_2) \\
\label{6.47e}
&&\partial_t B_{22} = t^{-1} B_{22} + t \ F_3 (P \partial_t N_2)\\
&&\partial_t {\check B}_{22} = F_2(x \cdot P \partial_t N_2)
\label{6.48e}
\eea

\noi where
\beq
\label{6.49e}
N_2 = 2 {\rm Im} \  \overline{w}_1 \nabla w_0 - (s_0 + B_{10}) 2 {\rm Re} \  \overline{w}_1w_0 - (s_1 + B_{11} + B_{21}) |w_0|^2 \ .
\eeq

We now begin the proof of the estimates. In all the proof, $m$ denotes a real number satisfying suitable conditions.

The properties (\ref{6.13e}) (\ref{6.14e}) of $w_0$ are obvious.

From (\ref{6.32e}) (\ref{3.18e}) we estimate
$$\parallel \omega^{m+1} B_{10} \parallel_2\ \leq I_m \left ( \parallel \omega^m x |w_0|^2 \parallel_2 \right ) \leq C$$

\noi for $0 \leq m \leq k_+$ by Lemma 3.2. Furthermore, $B_{10}$ has $j = 1$, $k = 0$ in (\ref{6.31e}) so that $B_{10} \in L^4$. Similarly, from (\ref{6.33e}) (\ref{3.19e}), we estimate 
$$\parallel \omega^{m+1} {\check B}_{10} \parallel_2\ \leq t^{-1}I_{m-1} \left ( \parallel <x> \omega^m x |w_0|^2 \parallel_2 \right ) \leq C\ t^{-1}$$

\noi for $1 \leq m \leq k_+ - 1$. Furthermore ${\check B}_{10}$ has $j = k = 0$ in (\ref{6.31e}), so that $\omega {\check B}_{10} \in L^4$. Together with the properties of $B_{10}$, this proves (\ref{6.15e}). From (\ref{6.34e}) (\ref{3.18e}), we estimate
$$\parallel \omega^{m} \partial_t B_{10} \parallel_2\ \leq I_m \left ( \parallel \omega^{m-1} x \overline{w}_0\partial_t w_0 \parallel_2 \right ) \leq C$$

\noi for $0 \leq m \leq k_+ - 1$. This proves (\ref{6.16e}). Similarly, from (\ref{6.35e}) (\ref{3.19e}) we estimate
$$\parallel \omega^{m+1} \left ( t \partial_t {\check B}_{10}+  {\check B}_{10}\right ) \parallel_2\ \leq I_{m} \left ( \parallel <x> \omega^m x \overline{w}_0\partial_t w_0\parallel_2 \right ) \leq C$$

\noi for $0 \leq m \leq k_+ - 3$. Furthermore $t \partial_t {\check B}_{10} + {\check B}_{10}$ has $j = 1$, $k = 0$ in (\ref{6.31e}), so that $t \partial_t {\check B}_{10}+  {\check B}_{10} \in L^4$. Together with (\ref{6.15e}), this proves (\ref{6.17e}). 

The properties (\ref{6.18e}) of $s_0$ follow from (\ref{6.15e}) and from the properties of $g$ in Lemma 3.5.

From (\ref{6.36e}) (\ref{3.18e}), we estimate
$$\parallel \omega^{m} B_{21} \parallel_2\ \leq t\ I_m \left ( \parallel \omega^{m-1} PN_1\parallel_2 \right ) \leq C\ t(1 - \ell n\ t )$$

\noi for $0 \leq m \leq k_+ + 1$. Here we have used Lemma 3.3 to prove that $P{\rm Im}\overline{w}_0\nabla w_0$ and $Ps_0 |w_0|^2 = P (\nabla \varphi_0)|w_0|^2$ belong to $H^{k_+}$, using the fact that $s_0 = \nabla \varphi_0$ is a gradient, with properties of $\varphi_0$ that can be read from (\ref{6.18e}). Similarly from (\ref{6.37e}) (\ref{3.19e}) we estimate 
$$\parallel \omega^{m+1}{\check B}_{21}\parallel_2\ \leq I_{m} \left ( \parallel \omega^m x\cdot PN_1 \parallel_2 \right ) \leq C(1 - \ell n\ t)$$

\noi for $0 \leq m \leq k_+$. Here we have used Lemma 3.4 to show that $P \cdot x {\rm Im} \overline{w}_0 \nabla w_0$ and $P\cdot x s_0 |w_0|^2 = P \cdot x \nabla \varphi_0|w_0|^2$ belong to $H^{k_+}$. Furthermore, ${\check B}_{21}$ has $j = 1$, $k = 0$ in (\ref{6.31e}), so that ${\check B}_{21} \in L^4$. This completes the proof of (\ref{6.19e}). \par

From (\ref{6.38e}) (\ref{3.18e}), we estimate
$$\parallel \omega^{m} \left ( \partial_t {B}_{21}- t^{-1}  {B}_{21}\right ) \parallel_2\ \leq t\ I_{m+1} \left ( \parallel \omega^{m-1} P \partial_t N_1\parallel_2 \right ) \ .$$

\noi The contribution of all the terms from $\partial_t N_1$ not containing $\partial_t s_0$ is estimated by $C t(1 - \ell n\ t)$ for $0 \leq m \leq k_+ - 1$, while the contribution of the term containing $\partial_t s_0$ is estimated by $C$ for $0 \leq m \leq k_+ + 1$, again by Lemma 3.3 and the fact that $s_0 = \nabla \varphi_0$. Together with (\ref{6.19e}), this proves (\ref{6.20e}). Similarly, from (\ref{6.39e}) (\ref{3.19e}), we estimate 
$$\parallel \omega^{m} \partial_t {\check B}_{21} \parallel_2\ \leq I_{m} \left ( \parallel \omega^{m -1}x \cdot P \partial_t N_1  \parallel_2 \right )\ .$$

\noi The contribution of the terms from $\partial_tN_1$ not containing $\partial_ts_0$ is estimated by\break\noindent $C(1 - \ell n\ t)$ for $0 \leq m \leq k_+ - 1$, while the contribution of the term containing $\partial_t s_0$ is estimated by $C t^{-1}$ for $1 \leq m \leq k_+ + 1$ by Lemma 3.4 and the fact that $s_0 = \nabla \varphi_0$. Furthermore the term containing $\partial_t s_0$ has $j=2$, $k = -1$ in (\ref{6.31e}) and therefore belongs to $L^4$. This proves (\ref{6.21e}). \par

From the definition of $w_1$ in (\ref{6.3e}), it follows that $w_1$ satisfies (\ref{6.22e}) (\ref{6.23e}). In fact, the basic estimate is that of $\partial_t w_1$. The upper bound $k_1 \leq k_+ - 2$ comes from $\nabla \cdot s_0$ as estimated by (\ref{6.18e}), while the upper bound $k_1 \leq k_+ - \beta^{-1}$ comes from ${\check B}_{10S}$ through the estimate
\begin{eqnarray*}
&&\parallel \omega^{k_+ - \beta^{-1}} {\check B}_{10S} w_0\parallel_2 \  \leq  \ \parallel \omega^{k_+ - \beta^{-1}} {\check B}_{10S}\parallel_2\ \parallel w_0 \parallel_{\infty} \\
&&\qquad \qquad  + \parallel {\check B}_{10S} \parallel_{\infty} \ \parallel \omega^{k_+ - \beta^{-1}} w_0\parallel_2\\
&&\leq \ \parallel \omega^{k_+} {\check B}_{10} \parallel_2 \left ( t \parallel w_0\parallel_{\infty} \ + t^{\beta (k_+ - 3/2)} \parallel \omega^{k_+- \beta^{-1}} w_+ \parallel_2 \right ) \leq C
\end{eqnarray*}

\noi by (\ref{6.15e}). The estimate of $w_1$ is obtained from that of $\partial_t w_1$ by integration in time. The estimate of $x \partial_t w_1$ is the same as that of $\partial_t w_1$, with $x$ absorbed by $w_0$, and the estimate of $xw_1$ follows therefrom by integration in time.\par

From  (\ref{6.41e}) (\ref{3.18e}), we estimate
$$\parallel \omega^{m} {B}_{11} \parallel_2\ \leq I_{m-1} \left ( \parallel \omega^{m-1} x \overline{w}_0 w_1\parallel_2 \right ) \leq C\ t(1 - \ell n\ t)^2$$

\noi for $0 \leq m \leq k_1 + 1$. Similarly, from (\ref{6.42e}) (\ref{3.19e}), we estimate 
$$\parallel \omega^{m+1} {\check B}_{11} \parallel_2\ \leq t^{-1} I_{m-1} \left ( \parallel <x> \omega^{m }x \overline{w}_0 w_1 \parallel_2 \right )\ \leq C(1 - \ell n\ t)^2$$

\noi for $0 \leq m \leq k_1$. Furthermore ${\check B}_{11}$ has $j = 0$, $k=1$ in (\ref{6.31e}) and therefore ${\check B}_{11} \in L^4$. This completes the proof of (\ref{6.24e}).\par

From  (\ref{6.43e}) (\ref{3.18e}), we estimate
$$\parallel \omega^{m} \partial_t {B}_{11} \parallel_2\ \leq I_{m} \left ( \parallel \omega^{-1} x \partial_t (\overline{w}_0 w_1)\parallel_2 \right ) \leq C\ (1 - \ell n\ t)^2$$

\noi for $0 \leq m \leq k_1 + 1$. This proves (\ref{6.25e}).\par

From  (\ref{6.44e}) (\ref{3.19e}), we estimate
$$\parallel \omega^{m} (t \partial_t {\check B}_{11} + {\check B}_{11})\parallel_2\ \leq I_{m-1} \left ( \parallel <x> \omega^{m-1 }x \partial_t (\overline{w}_0 w_1) \parallel_2 \right )\ .$$

\noi The contribution of the term with $\partial_t w_0$ is estimated by $Ct(1 - \ell n\ t)^2$ for $0 \leq m \leq (k_+-2) \wedge (k_1 + 1)$ while the contribution of the term with $\partial_t w_1$ is estimated by $C(1 - \ell n\ t)^2$ for $1 \leq m \leq k_1 + 1$. The latter term has $j = 1$, $k=0$ in (\ref{6.31e}) and therefore belongs to  $L^4$. This completes the proof of (\ref{6.26e}).\par

The properties (\ref{6.27e}) of $s_1$ follow from (\ref{6.24e}) and from the properties of $g$ in Lemma 3.5. \par

From  (\ref{6.45e}) (\ref{3.18e}), we estimate
$$\parallel \omega^{m} {B}_{22} \parallel_2\ \leq t\ I_{m} \left ( \parallel \omega^{m-1} N_2 \parallel_2 \right ) \leq C\ t^2(1 - \ell n\ t)^3$$

\noi for $0 \leq m \leq k_1 + 1$. Similarly, from (\ref{6.46e}) (\ref{3.19e}), we estimate 
$$\parallel \omega^{m} {\check B}_{22} \parallel\ \leq  I_{m-1} \left ( \parallel <x> \omega^{m-1 }N_2 \parallel_2 \right )\ \leq C\ t (1 - \ell n\ t)^3$$

\noi for $0 \leq m \leq k_1+1$. This proves (\ref{6.28e}).\par

From (\ref{6.47e}) (\ref{3.18e}), we estimate
$$\parallel \omega^{m} \left ( \partial_t {B}_{22} - t^{-1} B_{22} \right ) \parallel_2\ \leq t\ I_{m+1} \left ( \parallel \omega^{m-1} \partial_t N_2 \parallel_2 \right )\ .$$

\noi The contribution of the term with $\nabla \partial_t w_0$ is estimated by $Ct^2(1 - \ell n\ t)^2$ for $0 \leq m \leq (k_+-2) \wedge (k_1 + 1)$. The contribution of all the other terms is estimated by  $Ct(1 - \ell n\ t)^3$ for $0 \leq m \leq k_1 + 1$. This proves (\ref{6.29e}).\par

Finally, from (\ref{6.48e}) (\ref{3.19e}), we estimate
$$\parallel \omega^{m} \partial_t {\check B}_{22} \parallel\ \leq  I_{m} \left ( \parallel <x> \omega^{m-1 }\partial_t N_2 \parallel_2 \right )\ .$$

\noi The contribution of the term with $\nabla \partial_t w_0$ is estimated by $Ct(1 - \ell n\ t)^2$ for $0 \leq m \leq (k_+-2) \wedge (k_1 + 1)$. The contribution of all the other terms is estimated by  $C(1 - \ell n\ t)^3$ for $0 \leq m \leq k_1 + 1$. Actually, as compared with $\partial_t B_{22}$, the additional factor $<x>$ is absorbed by $w_1$ with no loss in regularity or decay, and by $w_0$ in $|w_0|^2$ and $\overline{w}_0 \partial_tw_0$ with no loss in decay. This proves (\ref{6.30e}).\par \nobreak
\hfill $\sq$ \par

We summarize the information on $(w_a, s_a, B_a)$ which follows from Lemma 6.1 in the following proposition.\\

\noi {\bf Proposition 6.1.} {\it Let $w_+ \in H^{k_+}$, $xw_+ \in H^{k_+-1}$ with $k_+ \geq 5 \vee (3 + \beta^{-1})$. Let $k_1 = (k_+ - 2) \wedge (k_+ - \beta^{-1})$. Let $(w_a, s_a, B'_a)$ be defined by (\ref{6.1e}) (\ref{6.2e}) (\ref{6.3e}) and $B'_a = B_{1a} + B_{2a}$. Then $(w_a, s_a, B'_a)$ satisfy the following properties}
\beq
\label{6.50e}
w_a \in \left ( H^{k_1} , 1\right ) \quad , \qquad x w_a \in \left ( H^{k_1} , 1\right ) \ ,  
\eeq
\beq
\label{6.51e}
 \partial_t  w_a \in \left ( H^{k_1} , (1- \ell n\ t)^2\right ) \quad , \quad  \partial_t  xw_a \in \left ( H^{k_1\wedge (k_+-3)} , (1- \ell n\ t)^2\right ) \ , 
\eeq
\beq
\label{6.52e}
s_a \in \left ( \ddot{H}^{k_1} , (1- \ell n\ t)\right ) \quad , \quad  \partial_t  s_a \in \left ( \ddot{H}^{k_1} , t^{-1}\right ) \ , 
\eeq
\beq
\label{6.53e}
B'_a \in \left ( \ddot{H}^{k_1+1}, 1\right ) \quad , \quad {\check B}'_a \in \left ( <x> \dot{H}^1 \cap \dot{H}^2 \cap \dot{H}^{k_1 + 1}, t^{-1} \right ) \ , 
 \eeq
\beq
\label{6.54e}
\left \{ \begin{array}{l} \partial_t B'_a \in \left ( H^{k_1+1}, (1- \ell n\ t)^2\right )\\ \\ \partial_t {\check B}'_a \in \left ( <x> \dot{H}^1 \cap \dot{H}^2 \cap \dot{H}^{(k_1 + 1)\wedge (k_+ - 2) }, t^{-2} \right ) \ . \end{array} \right . 
\eeq

\noi {\it Furthermore  ${\check B}'_{aS} \equiv \chi_S  {\check B}_{1a}+  {\check B}_{2a}$ satisfies the estimate}
\beq
\label{6.55e}
\parallel {\check B}'_{aS}\parallel_{\infty} \ + \ t \parallel \partial_t {\check B}'_{aS}\parallel_{\infty}\ \leq C(1 - \ell n\ t) \ .
\eeq
\vskip 5 truemm

\noi {\bf Proof.} The properties and estimates (\ref{6.50e})-(\ref{6.55e}) follow from Lemma 6.1 except for the contribution of ${\check B}_{1a0}$ to (\ref{6.55e}). We estimate
$$\parallel {\check B}_{1a0S}\parallel_{\infty}\ \leq t^{\beta (k_+ - 3/2)} \parallel \omega^{k_+}{\check B}_{1a0}\parallel_2\ \leq C\ t^{\beta /2}$$

\noi by (\ref{6.15e}). Similarly
$$\parallel {\check B}_{1a1S}\parallel_{\infty}\ \leq t^{\beta (k_1 - 1/2)} \parallel \omega^{k_1 + 1}{\check B}_{1a1}\parallel_2\ \leq C\ t^{(1 - \beta /2) \wedge 3\beta /2} (1 - \ell n \ t)^2$$

\noi by (\ref{6.24e}) since
$$\beta k_1 = (\beta k_+ - 2 \beta ) \wedge (\beta k_+ - 1) \geq 1 \wedge 2 \beta \ .$$

The estimate of $\partial_t {\check B}_{1aS}$ is obtained by combining similar arguments with the estimate of $\partial_t {\check B}_{1a}$, taking into account the time derivative of $\chi_S$. \par \nobreak \hfill $\sq$\par

We now turn to the estimates of the remainders. The final result will be that the remainders satisfy the assumption (A3) of Section 5 with 
$$h(t) = t^2 (1 - \ell n\ t)^4\ .$$

\noi We first consider the part not containing $B_0$, namely $R_{11}$, $R_2$ , $R_3$ and $R_{41}$. The estimates for that part follow from or extend Lemma 6.1. The part containing $B_0$ requires different arguments and additional assumptions.\\

\noi {\bf Proposition 6.2.} {\it Let $w_+ \in H^{k_+}$, $xw_+ \in H^{k_+-1}$ with $k_+ \geq 5 \vee (3 + \beta^{-1})$and let $k_1 = (k_+ - 2) \wedge (k_+ - \beta^{-1})$. Then the remainders $R_2$, $R_3$ and $R_{41}$ satisfy the following properties}
\beq
\label{6.56e}
R_2 \in \left ( H^{k_1 + 1}, t(1 - \ell n\ t)^4\right ) \ , 
\eeq
\beq
\label{6.57e}
R_3, t {\check R}_3  \in \left ( H^{k_1 + 1}, t^2(1 - \ell n\ t)^4\right ) \ , 
\eeq
\beq
\label{6.58e}
\partial_t R_3, t \partial_t {\check R}_3  \in \left ( H^{k_1 + 1}, t(1 - \ell n\ t)^4\right ) \ , 
\eeq
\beq
\label{6.59e}
R_{41}, t {\check R}_{41}  \in \left ( H^{k_1 + 1}, t^3(1 - \ell n\ t)^5\right ) \ , 
\eeq
\beq
\label{6.60e}
\partial_t R_{41}, t  \partial_t {\check R}_{41}  \in \left ( H^{k_1 + 1}, t^2(1 - \ell n\ t)^5\right ) \ . 
\eeq

\noi {\it Let in addition $k_+ \geq 2 \beta^{-1}$ and define
\beq
\label{6.61e}
k_2 = \left ( k_1 - 2\right ) \wedge \left ( k_1 + 1 - \beta^{-1}\right )\ .
\eeq
\noi Then the remainder $R_{11}$ satisfies the following properties~:}
\beq
\label{6.62e}
R_{11} , x R_{11} \in  \left ( H^{k_2}, t(1 - \ell n\ t)^4\right ) \ , 
\eeq
\beq
\label{6.63e}
\partial_t R_{11} , x \partial_t R_{11} \in  \left ( H^{k_2}, (1 - \ell n\ t)^4\right ) \ . 
\eeq
\vskip 5 truemm

\noi {\bf Proof.} In all the proof, we omit the subscript $a$.\par

The property (\ref{6.56e}) of $R_2$ follows from (\ref{6.22e}) and from the estimates of $g$ in Lemma 3.5. \par

The properties (\ref{6.57e})  (\ref{6.58e}) of $R_3$ follow from the properties (\ref{6.22e}) (\ref{6.23e}) of $w_1$ by estimates similar to those contained in the proof of Lemma 6.1.\par

We next turn to $R_{41}$ which we rewrite as
$$R_{41} = - t F_2 \left ( P\ {\rm Im} \overline{w}_1 \nabla w_1 - P ((s+ B') |w|^2 )_{\geq 2} \right ) \ . $$

\noi By Proposition 6.1, the space regularity of $w$, $s$, $B'$ is at least $\ddot{H}^{k_1}$. Together with Lemma 3.3, this proves that $R_{41}$ has the regularity of $H^{k_1 + 1}$. The time decay follows from Lemma 6.1, with the worst term being $s_0 |w_1|^2 \leq O (t^2 (1 - \ell n\  t)^5)$. This proves the first part of (\ref{6.59e}). The second part, namely the estimate of  ${\check R}_{41}$, follows from the fact that the additional factor $<x>$ can be absorbed by $w_1$ or $w_0$ with a loss of regularity by zero or one space derivative respectively, and in both cases without any change in the time decay.\par

We next consider the time derivative of $R_{41}$. A time derivative produces at most a loss of a factor $t$ when acting on the various components of $(w, s, B', {\check B}')$. It produces no loss of regularity when acting on $w_1$, $s$ and therefore in the terms generated by $\partial_t w_1$ or $\partial_t s$ in $\partial_t B'$ and $\partial_t {\check B}'$. It produces a loss of regularity by at most two space derivatives when acting on $w_0$ and in the terms generated by $\partial_tw_0$ in $\partial_t B'$ and $\partial_t {\check B}'$. However, that loss of regularity occurs only in terms where the regularity starts from a sufficiently high level, so that it does not affect the final result. This proves (\ref{6.60e}).\par

We finally consider $R_{11}$ which we rewrite as 
$$R_{11} = R'_{11} - \left ( {\check B}_{1S} w\right )_{\geq 1}\ , $$
$$R'_{11} = (1/2) \Delta w_1 - \left \{ i (s+ B')\cdot \nabla w + (i/2) (\nabla \cdot s) w + (1/2) (s + B')^2 w \right \}_{\geq 1} -  \left ( {\check B}_{2} w\right )_{\geq 2}\ . $$

\noi We first consider $R'_{11}$. By Proposition 6.1 the space regularity of $w$, $s$, $B'$ is at least $\ddot{H}^{k_1}$, so that $R'_{11} \in H^{k_1 - 2}$. The time decay follows from Lemma 6.1, with the worst term coming from $s_0^2 w_1 \leq O (t(1 - \ell n\ t)^4)$. Therefore
\beq
\label{6.64e}
R'_{11} \in \left ( H^{k_1-2}, t (1 - \ell n\  t)^4\right ) \ .
\eeq

\noi We next consider
$$ \left ( {\check B}_{1S} w\right )_{\geq 1} = {\check B}_{10S}w_1 + {\check B}_{11S}w_0 + {\check B}_{11S}w_1\ . $$

It follows from Lemma 6.1, especially (\ref{6.24e}), that
\beq
\label{6.65e}
{\check B}_{11} w_1 \in \left ( H^{k_1}, t(1 - \ell n\ t)^4 \right )\ .
\eeq

In particular the use of the $S$ cut off is not needed for that term. Using (\ref{6.15e}), we next estimate
\begin{eqnarray*}
\parallel \omega^{k_1} {\check B}_{10S}w_1\parallel_2&\leq& \parallel \omega^{k_1} {\check B}_{10S}\parallel_2 \ \parallel w_1\parallel_{\infty}\ + \parallel {\check B}_{10S}\parallel_{\infty} \ \parallel \omega^{k_1} w_1\parallel_2 \\
& \leq& C \left ( t^{\beta ( k_+ - k_1)} + t^{\beta (k_+ - 3/2)} \right ) \parallel \omega^{k_+} {\check B}_{10}\parallel_2 t(1 - \ell n\ t)^2\\
&\leq& C\ t(1 - \ell n\ t)^2
\end{eqnarray*}

\noi so that
\beq
\label{6.66e}
{\check B}_{10S} w_1 \in  \left ( H^{k_1}, t(1 - \ell n\ t)^2 \right )\ .
\eeq

\noi Using (\ref{6.24e}), we next estimate
$$\parallel \omega^{k_2} {\check B}_{11S}w_0\parallel_2\ \leq \ \parallel \omega^{k_2} {\check B}_{11S}\parallel_2 \ \parallel w_0\parallel_{\infty}\ + \parallel {\check B}_{11S}\parallel_{r} \ \parallel \omega^{k_2} w_0\parallel_{3/\delta}$$

\noi for $0 < \delta \leq k_2 \leq k_1 + 1$, with $\delta = \delta (r) = 3/2 - 3/r$,
\begin{eqnarray*}
\cdots & \leq& C \left ( t^{\beta ( k_1 + 1 - k_2)} + t^{\beta (k_1+1 - \delta )} \right ) \parallel \omega^{k_1 + 1} {\check B}_{11}\parallel_2 \\
&\leq& C\ t(1 - \ell n\ t)^2
\end{eqnarray*}

\noi for $k_2 \leq k_1 + 1 - \beta^{-1}$, so that 
\beq
\label{6.67e}
{\check B}_{11S} w_0 \in  \left ( H^{k_2}, t(1 - \ell n\ t)^2 \right )\ .
\eeq

It follows from  (\ref{6.64e})-(\ref{6.67e}) that $R_{11}$ satisfies the first part of  (\ref{6.62e}). We now turn to the second part, namely the estimate of $xR_{11}$. Now $x$ can be absorbed either by $w_0$ with a loss of regularity by one space derivative or by $w_1$ without a loss of regularity, in both cases without any change in the time decay. This proves the second part of  (\ref{6.62e}). Finally the time derivative produces the same effects as in $R_{41}$, thereby leading to  (\ref{6.63e}). \par \nobreak \hfill $\sq$ \par

We now turn to the estimates of the parts $R_{10}$ and $R_{40}$ of the remainders containing $B_0$. We shall need the following estimate
\beq
\label{6.68e}
\parallel <x>^{\ell} \partial_t^j (w_a - w_+) ; H^{3-\ell}\parallel \ \leq C\ t^{1-j} (1 - \ell n\ t)^2
\eeq

\noi for $j, \ell = 0, 1$, which holds under the assumptions of Lemma 6.1. That estimate follows from (\ref{6.13e}) (\ref{6.14e}) (\ref{6.22e}) (\ref{6.23e}) and elementary arguments. We consider first $R_{10}$ which we rewrite as
\beq
\label{6.69e}
R_{10} = {\cal R}_1 (w_a)
\eeq

\noi where 
\beq
\label{6.70e}
{\cal R}_1 (v) = - i B_0 \cdot \nabla v - \left ( B_0 (s_a + B'_a) + (1/2) B_0^2 + {\check B}_0 \right ) v\ .
\eeq

We first derive estimates of some parts of $R_{10}$ which can be obtained without additional assumptions on $w_+$ beyond those of Lemma 6.1.\\

\noi {\bf Lemma 6.2.} {\it Let $w_+ \in H^{k_+}$, $xw_+ \in H^{k_+-1}$ with $k_+ \geq 5 \vee (3 + \beta^{-1})$. Let $B_0$ satisfy (\ref{3.32e}) for $2 \leq r \leq \infty$ and $0 \leq j, k \leq 1$. Then the following estimates hold~: }
\beq
\label{6.71e}
\parallel <x> \partial_t {\cal R}_1 (w_a - w_+) \parallel_2 \ \leq C\ t^{1/3} (1 - \ell n\ t)^2 \ ,
\eeq
\beq
\label{6.72e}
\parallel \partial_t {\cal R}_1 (\nabla w_+) \parallel_2 \ \leq C\ t^{-1/2} (1 - \ell n\ t) \ ,
\eeq
\beq
\label{6.73e}
\parallel \nabla  \partial_t {\cal R}_1 (w_a - w_+) \parallel_2 \ \leq C\ t^{-1/2} (1 - \ell n\ t)^3 \ ,
\eeq
\beq
\label{6.74e}
\parallel \left ( (\partial_t B_0) (\nabla (s_a + B_a)\right ) + B_0 \left ( \nabla \partial_t (s_a + B'_a))\right ) w_+ \parallel_2 \ \leq C\ t^{-1/2} (1 - \ell n\ t) \ .
\eeq
\vskip 5 truemm

\noi {\bf Proof.} In all the proof, we omit the subscript $a$. Let $v = w - w_+$. We start from
\bea
\label{6.75e}
&&- \partial_t {\cal R}_1(v) =  (\partial_t B_0 ) \left ( i \nabla v + (s+B) v\right ) + B_0 \left ( \partial_t (s + B') \right ) v + (\partial_t {\check B}_0 ) v \nn \\
&&+ i B_0 \nabla \partial_t v + \left ( B_0 (s+B') + (1/2) B_0^2 + {\check B}_0 \right ) \partial_t v \ .
\eea

\noi We first prove (\ref{6.71e}) and for that purpose we estimate
$$\parallel <x> \partial_t {\cal R}_1(v)\parallel_2\  \leq \ \parallel \partial_t B_0\parallel_2 \ \parallel<x> \nabla v \parallel_{\infty}\ + \ \Big ( \parallel \partial_t B_0 \parallel_2 \ \parallel s+B\parallel_{\infty}$$
$$+\ \parallel B_0 \parallel_2 \ \parallel\partial_t (s+ B') \parallel_{\infty} \ + \ \parallel \partial_t {\check B}_0\parallel_2 \Big ) \parallel<x> v \parallel_{\infty}\ + \ \parallel B_0 \parallel_3 \ \parallel<x> \nabla \partial_t v \parallel_6$$
$$+ \left ( \parallel B_0 \parallel_2 \left ( \parallel s + B' \parallel_{\infty} \ + \ \parallel B_0 \parallel_{\infty} \right ) \ + \ \parallel {\check B}_0 \parallel_2 \right ) \parallel <x> \partial_t v \parallel_{\infty}$$
\beq
\label{6.76e}
  \leq C\ t^{1/3} (1 - \ell n\ t)^2
\eeq

\noi by Proposition 6.1, (\ref{6.68e}) and (\ref{3.32e}). This proves (\ref{6.71e}).

The proof of (\ref{6.72e}) is obtained from (\ref{6.76e}) by omitting $<x>$ and the terms containing $\partial_t v$, replacing $v$ by $\nabla w_+$ in the remaining terms, and using again Proposition 6.1 and (\ref{3.32e}).\par

We next prove (\ref{6.73e}). For that purpose we estimate the gradient of (\ref{6.75e}). We obtain
$$\parallel \nabla \partial_t {\cal R}_1(v) \parallel _2\  \leq \ \parallel \partial_t B_0 \parallel _3 \ \parallel \nabla^2 v \parallel _6\ + \ \Big ( \parallel  \nabla \partial_t B_0 \parallel_2 \ + \ \parallel \partial_t  B_0 \parallel_2\ \parallel s+ B\parallel_{\infty}$$
$$+ \ \parallel B_0 \parallel_2 \ \parallel  \partial_t (s+ B') \parallel_{\infty} \ + \ \parallel \partial_t {\check B}_0 \parallel _2 \Big ) \parallel  \nabla v \parallel _{\infty}$$
$$+\ \Big ( \parallel \nabla \partial_t B_0 \parallel _2\ \parallel  s + B \parallel _{\infty} \ + \ \parallel \partial_t B_0 \parallel_2 \ \parallel \nabla (s+B) \parallel _{\infty} \ + \ \parallel \nabla B_0 \parallel _2 \ \parallel  \partial_t (s+ B') \parallel_{\infty}$$
$$+\ \parallel  B_0\parallel _2\ \parallel \nabla \partial_t (s+B')\parallel _{\infty} \ + \ \parallel \nabla \partial_t {\check B}_0 \parallel_2 \Big )\parallel v\parallel _{\infty} \ + \ \parallel B_0 \parallel _3\ \parallel  \nabla^2 \partial_t v \parallel _6$$
$$ +\ \Big ( \parallel  \nabla B_0 \parallel _2\ + \ \parallel  B_0 \parallel _2 \left ( \parallel  s+ B' \parallel _{\infty} \ + \ \parallel  B_0\parallel _{\infty }\right )  \ + \ \parallel  {\check B}_0 \parallel _2 \Big ) \parallel  \nabla \partial_t v \parallel _{\infty}$$ 
$$ +\ \Big ( \parallel  \nabla B_0 \parallel _2\   \parallel  s+ B \parallel _{\infty} \  + \ \parallel  B_0 \parallel _2\ \parallel  \nabla (s+B')\parallel _{\infty}\ + \ \parallel  \nabla {\check B}_0 \parallel _2 \Big ) \parallel  \partial_t v \parallel _{\infty}$$
$$\leq C\ t^{-1/2} (1 - \ell n\ t )^3$$

\noi by Proposition 6.1,  (\ref{6.68e}) and (\ref{3.32e}).\par

The proof of (\ref{6.74e}) is similar.\par \nobreak \hfill $\sq$\par

The full estimates of $R_{10}$ require additional assumptions on $w_+$ and $B_0$.\\

\noi {\bf Lemma 6.3.} {\it Let $w_+$ and $B_0$ satisfy the assumptions of Lemma 6.2. \par

(1) Let $w_+$ and $B_0$ satisfy in addition
\beq
\label{6.77e}
\parallel  \partial_t B_0 <x> \nabla^k w_+ \parallel _2\ \vee \ t^{-1} \parallel B_0 <x> w_+  \parallel _2\ \vee \ \parallel \partial_t {\check B}_0 <x> w_+ \parallel _2 \ \leq C
\eeq

\noi for $k = 0, 1$. Then
\beq
\label{6.78e}
\parallel  <x> \partial_t  R_{10} \parallel_2\ \leq C(1- \ell n\ t )\ .
\eeq

(2) Let $w_+$ and $B_0$ satisfy in addition
\beq
\label{6.79e}
\parallel  ( \nabla \partial_t B_0) \nabla^k w_+ \parallel _2\ \vee \ t^{-1} \parallel (\nabla B_0)w_+  \parallel _2\ \vee \ \parallel (\nabla \partial_t {\check B}_0) w_+ \parallel _2 \ \leq C\ t^{-1/2}
\eeq

\noi for $k = 0,1$. Then}
\beq
\label{6.80e}
\parallel  \nabla \partial_t  R_{10} \parallel_2\ \leq C\ t^{-1/2} (1- \ell n\ t )^3 \ .
\eeq

\vskip 5 truemm

\noi {\bf Proof.} We again omit the subscript $a$ in the proof.\par

\noi \underline{Part (1).} Using (\ref{6.75e}), we estimate
\bea
\label{6.81e}
&&\parallel <x> \partial_t R_{10} \parallel_2\ \leq \ \parallel<x> \partial_t {\cal R}_1 (w- w_+) \parallel_2\ + \ \parallel \partial_t B_0 <x> \nabla w_+ \parallel_2\nn \\
&&+\ \parallel s+ B \parallel_{\infty} \ \parallel \partial_t B_0 <x> w_+\parallel_2\ + \ \parallel \partial_t (s+ B') \parallel_{\infty} \ \parallel B_0 <x> w_+ \parallel_2\nn \\
&&+\ \parallel \partial_t {\check B}_0 <x> w_+ \parallel_2\ \leq C(1 - \ell n\ t)
\eea

\noi by (\ref{6.71e}) (\ref{6.77e}) and Proposition 6.1.\\

\noi \underline{Part (2).} Using (\ref{6.75e}), we estimate
$$\parallel \nabla  \partial_t R_{10} \parallel_2\ \leq \ \parallel \nabla \partial_t {\cal R}_1 (w- w_+) \parallel_2\ + \ \parallel \partial_t {\cal R}_1(\nabla w_+)\parallel_2$$
$$+\ \parallel \left ( (\partial_t B_0) (\nabla (s+ B)\right )+ B_0 \left ( \nabla \partial_t (s+B'))\right ) w_+  \parallel_{2}$$
$$+ \ \parallel (\nabla \partial_t B_0 ) \nabla  w_+\parallel_2\ + \ \parallel (s+B) \parallel_{\infty}\ \parallel (\nabla \partial_t B_0)w_+\parallel_{2} \ +$$
$$+\ \parallel \partial_t (s+ B') \parallel_{\infty} \ \parallel  (\nabla B_0) w_+ \parallel_2\ + \ \parallel (\nabla \partial_t {\check B}_0) w_+ \parallel_2$$
$$\leq C\ t^{-1/2} (1 - \ell n\ t)^3$$

\noi by (\ref{6.72e}) (\ref{6.73e}) (\ref{6.74e}) (\ref{6.79e}) and Proposition 6.1.\par\nobreak \hfill $\sq$\par

We next turn to $R_{40}$ which we rewrite as
\beq
\label{6.82e}
R_{40} = {\cal R}_4 \left ( |w_a|^2\right )
\eeq

\noi where 
\beq
\label{6.83e}
 {\cal R}_4(f) = t\ F_2 (PB_0f) \ .
\eeq

We first show that one can replace $w_a$ by $w_+$ in the estimates of $R_{40}$ without additional assumptions on $w_+$ beyond those of Lemma 6.1. We shall need the following estimate
\beq
\label{6.84e}
\parallel <x>^{\ell} \partial_t^j (w_a + w_+) ; H^{3 - \ell} \parallel \ \leq C(1 - \ell n\ t)^{2j}
\eeq

\noi with $\ell = 0, 1$, which follows from (\ref{6.13e}) (\ref{6.14e}) (\ref{6.22e}) (\ref{6.23e}).\\

\noi {\bf Lemma 6.4.}  {\it Let $w_+ \in H^{k_+}$, $xw_+ \in H^{k_+-1}$ with $k_+ \geq 5 \vee (3 + \beta^{-1})$. Let $B_0$ satisfy (\ref{3.32e}) for $r= 2$ and $0 \leq j, k \leq 1$. Then the following estimates hold~:
\bea
\label{6.85e}
&&\parallel \nabla^{k+1} \partial_t^j {\cal R}_4(|w_a|^2 - |w_+|^2) \parallel _2\ \vee \ t \parallel  \nabla^{k+1} \partial_t^j {\check {\cal R}}_4(|w_a|^2 - |w_+|^2)\parallel _2\nn \\
&&\qquad \qquad \leq C\ t^{5/2-j-k} (1 - \ell n\ t)^2
\eea

\noi for $0 \leq j, k, j+k \leq 1$.}\\

\noi {\bf Proof.} We define $v_{\pm} = w_a \pm w_+$ so that
$${\cal R}_4 \left ( |w_a|^2 - |w_+|^2\right ) =  {\cal R}_4 \left ( {\rm Re}\  \overline{v}_+ v_-\right ) \ .$$

\noi Using (\ref{3.18e}) (\ref{3.19e}), we estimate
$$\parallel\nabla^{k+1} \partial_t^j {\cal R}_4 \left ( {\rm Re} \ \overline{v}_+ v_- \right ) \parallel_2 \ \leq t \ I_{j+k+1} \left ( \parallel \nabla^k \partial_t^j B_0 \overline{v}_+ v_- \parallel_2\right )$$
$$+ \delta_{j1}\ I_1 \left ( \parallel B_0 \overline{v}_+ v_- \parallel_2 \right ) \ ,$$

$$\parallel\nabla^{k+1} \partial_t^j {\check{\cal R}}_4 \left ( {\rm Re} \ \overline{v}_+ v_- \right ) \parallel_2 \ \leq I_{j+k} \left ( \parallel <x> \nabla^k \partial_t^j B_0 \overline{v}_+ v_- \parallel_2\right )$$

\noi for the relevant values of $j$ and $k$. The result now follows from (\ref{3.32e}) with $r=2$ and from (\ref{6.68e}) (\ref{6.84e}).\par\nobreak \hfill $\sq$ \par

The full estimates of $R_{40}$ require additional assumptions on $w_+$ and $B_0$.\\

\noi {\bf Lemma 6.5.} {\it Let $w_+$ and $B_0$ satisfy the assumptions of Lemma 6.4 and in addition}
\beq
\label{6.86e}
 \parallel\left  ( \nabla^k \partial_t^j B_0\right ) w_+ \parallel_2\ \leq C\ t^{3/2-j-k}
\eeq

\noi {\it for $0 \leq j, k, j+k \leq 1$. Then the following estimates hold~: }
\beq
\label{6.87e}
 \parallel  \nabla^{k+1} \partial_t^j R_{40} \parallel_2\ \vee \ t  \parallel  \nabla^{k+1} \partial_t^j {\check R}_{40} \parallel_2\ \leq C\ t^{5/2 - j - k} (1 - \ell n\ t)^2
\eeq

\noi {\it for $0 \leq j, k, j+k \leq 1$.}\\

\noi {\bf Proof.} By (\ref{6.82e}) and Lemma 6.4, it is sufficient to estimate ${\cal R}_4 (|w_+|^2)$. Using again (\ref{3.18e}) (\ref{3.19e}), we estimate 
$$\parallel\nabla^{k+1} \partial_t^j {\cal R}_4 (|w_+|^2) \parallel_2 \ \leq t \ I_{j+k+1} \left ( \parallel \nabla^k \partial_t^j B_0 |w_+|^2\parallel_2\right )$$
$$+ \delta_{j1}\ I_1 \left ( \parallel B_0 |w_+|^2 \parallel_2 \right ) \ ,$$

$$\parallel\nabla^{k+1} \partial_t^j {\check{\cal R}}_4 (|w_+|^2) \parallel_2 \ \leq I_{j+k} \left ( \parallel <x> \nabla^k \partial_t^j B_0 |w_+|^2 \parallel_2\right )$$

\noi for the relevant values of $j$ and $k$. We next estimate
$$\parallel <x>^{\ell} \left  ( \nabla^k \partial_t^j B_0\right ) |w_+|^2 \parallel_2\ \leq \  \parallel\left  ( \nabla^k \partial_t^j B_0\right ) w_+ \parallel_2\ \parallel <x> w_+ \parallel_{\infty}\ \leq C\ t^{3/2 - j - k}$$

\noi for $\ell = 0, 1$ and 
$$\parallel <x>^{\ell} B_0 \nabla |w_+|^2 \parallel_2\ \leq \  2\parallel B_0 \parallel_2 \ \parallel \nabla w_+ \parallel_{\infty}\ \parallel <x> w_+ \parallel_{\infty} \ \leq C\ t^{1/2}$$

\noi by (\ref{6.86e}) and (\ref{3.32e}). Substituting the last two estimates into the previous ones and using Lemma 6.4 yields (\ref{6.87e}). \par \nobreak \hfill $\sq$ \par

The additional assumptions (\ref{6.77e}) (\ref{6.79e}) (\ref{6.86e}) on $w_+$ and $B_0$ are special cases of the condition
\beq
\label{6.88e}
 \parallel  \left ( \nabla^{k} \partial_t^j B_{0}\right ) <x>^{\ell} \nabla^m w_+  \parallel_2\ \vee \   \parallel  \left ( \nabla^{k} \partial_t^j {\check B}_{0} \right ) <x>^{\ell} \nabla^m w_+ \parallel_2\ \leq C\ t^{3/2 - j - k - \ell /2}
\eeq

\noi with $0 \leq j,k, \ell , m \leq 1$ and $k+ \ell \leq 1$. That condition does not follow from factorized estimates of $B_0$ and $w_+$. In fact for $<x> w_+ \in H^4$, from (\ref{3.32e}) it follows only that
$$\parallel \left  ( \nabla^k \partial_t^j B_0\right ) <x>^{\ell}  \nabla^m w_+ \parallel_2\ \leq \  \parallel  \nabla^k \partial_t^j B_0 \parallel_2\ \parallel <x>^{\ell} \nabla ^m  w_+ \parallel_{\infty}\ \leq C\ t^{1/2 - j - k}$$

\noi and a similar estimate for ${\check B}_0$. This is weaker than (\ref{6.88e}) by a factor $t^{1-\ell /2}$. In order to gain that factor, we impose a support condition on $w_+$ and a decay condition of $B_0$ and ${\check B}_0$ on the support of $w_+$. In fact let $<x>w_+ \in H^3$ and let $\chi_0$ be the characteristic function of the support of $w_+$. Then a sufficient condition to ensure (\ref{6.88e}) is that 
\beq
\label{6.89e}
\parallel \chi_0 \nabla^k \partial_t^j B_0 \parallel_2\ \vee\ \parallel \chi_0 \nabla^k \partial_t^j {\check B}_0 \parallel_2\ \leq\ C\ t^{3/2 - j - k}
\eeq

\noi for $j, k = 0,1$. The support condition that we shall impose on $w_+$ is 
\beq
\label{6.90e}
{\rm Supp}\ w_+ \subset \{x: |\  | x| - 1| \geq \eta \}
\eeq

\noi for some $\eta$, $0 < \eta < 1$. This is the same condition that occurs in \cite{5r} \cite{28r}. Under that condition, it is easy to see that (\ref{6.89e}) holds for compactly supported $(A_+, \dot{A}_+)$. In fact, if 
$${\rm Supp}(A_+, \dot{A}_+) \subset \{x : |x| \leq R \}$$

\noi then by the Huyghens principle
$${\rm Supp}\ A_0 \cup {\rm Supp}\ x\cdot A_0 \subset \{(x,t): |\  | x| - t| \leq R \}$$

\noi so that
$${\rm Supp}\ B_0 \cup {\rm Supp}\ {\check B}_0 \subset \{(x,t): |\  | x| - 1| \leq t\ R \}$$

\noi and the left hand side of (\ref{6.89e}) vanishes for $t \leq \eta /R$. More general assumptions on $(A_+, \dot{A}_+)$ are given in the following lemma.\\

\noi {\bf Lemma 6.6.} {\it Let $w_+$ satisfy the support condition (\ref{6.90e}) for some $\eta$, $0 < \eta < 1$. Let $\chi_R$ be the characteristic function of the set $\{x:|x| \geq R\}$. Let $(A_+, \dot{A}_+)$ satisfy}
\beq
\label{6.91e}
\left \{ \begin{array}{l} 
\parallel \chi_R \nabla^k (x\cdot \nabla )^j A_+ \parallel_2\ \vee\ \parallel \chi_R \nabla^k (x\cdot \nabla)^j x\cdot A_+ \parallel_2\ \leq\ C\ R^{- 1}\\ \\
\parallel \chi_R (x\cdot \nabla )^j \dot{A}_+;L^2 \cap L^{6/5} \parallel\ \vee\ \parallel \chi_R (x \cdot  \nabla )^j x\cdot \dot{A}_+;L^2 \cap L^{6/5}  \parallel\ \leq\ C\ R^{- 1}
\end{array}\right .
\eeq

\noi {\it for $0 \leq j , k \leq 1$ and for all $R \geq R_0$ for some $R_0 > 0$. Then (\ref{6.89e}) holds for $0 \leq j,k \leq 1$ and for all $t\in (0, 1]$.}\\

\noi {\bf Proof.} For $j = 0$ and as regards $B_0$, the result is that of Lemma 5.2, part (2) of \cite{5r} to which we refer for the proof, which is a simple consequence of the Huyghens principle for the wave equation. The case $j = 1$ follows therefrom and from (\ref{3.24e}) (\ref{3.26e}). Finally the result for ${\check B}_0$ follows from that for $B_0$ and from (\ref{3.27e}).\par \nobreak \hfill $\sq$\par

We finally collect the results of this section to show that the asymptotic functions constructed here satisfy the assumptions (A1) (A2) (A3) of Section 5.\\

\noi {\bf Proposition 6.3.} {\it Let $w_+ \in H^{k_+}$, $xw_+ \in H^{k_+-1}$ with $k_+ \geq 5 \vee (3 + \beta^{-1})$. Let $B_0$ satisfy the condition (\ref{3.32e}) for $2 \leq r \leq \infty$ and $0 \leq j, k \leq 1$. Then \par

(1) The asymptotic functions $(w_a, s_a, B_a)$ defined by (\ref{6.1e}) (\ref{6.2e}) (\ref{6.3e}) satisfy the assumptions (A1) (A2).\par

Let in addition $k_+ \geq 2 \beta^{-1}$.\par

(2) Let in addition $B_0$ and $w_+$ satisfy the condition (\ref{6.88e}) for $0 \leq j, k, \ell , m \leq 1$ and $k + \ell \leq 1$. Then the remainders $R_j$ defined by (\ref{2.57e}) satisfy the assumption (A3) with 
$$h(t) = t^2 (1 - \ell n\ t)^4\ .$$

(3) The same result as in Part (2) holds under the assumptions of Lemma 6.6.}\\

\noi {\bf Proof.} Part (1) follows from Proposition 6.1 and from (\ref{3.32e}). Part (2) follows from Proposition 6.2 and Lemmas 6.3 and 6.5. The assumptions (\ref{6.77e}) (\ref{6.79e}) (\ref{6.86e}) of those lemmas are special cases of (\ref{6.88e}). Part (3) follows from Part (2), from (\ref{6.89e}) and from Lemma 6.6. \par \nobreak \hfill $\sq$ \par

\mysection{Final results for the auxiliary system and for the original system}
\hspace*{\parindent} In this section we complete the construction of the wave operators for the system (\ref{2.6e}) (\ref{2.7e}) and we derive asymptotic properties of solutions in their range. For that purpose we first state the main result on the Cauchy problem at $t= 0$ for the auxiliary system (\ref{2.41e}) which follows from Sections 5 and 6. In all this section we take $\beta = 1/2$.\\

\noi {\bf Proposition 7.1.} {\it Let $\beta = 1/2$. Let $X(\cdot )$ be defined by (\ref{3.10e}) with $h(t) =\break\noindent  t^2 (1 - \ell n \ t)^4$. Let $u_+$ be such that $w_+ \equiv Fu_+ \in H^5$, $xw_+ \in H^4$. Let $B_0$ satisfy the conditions (\ref{3.32e}) and (\ref{6.88e}) for $2 \leq r \leq \infty$, for $0 \leq j, k, \ell , m \leq 1$ and $\ell + k \leq 1$. Define $(w_a, s_a, B_a)$ by (\ref{6.1e})  (\ref{6.2e})   (\ref{6.3e}). Then there exists $\tau$, $0 < \tau \leq 1$ such that the auxiliary system  (\ref{2.41e}) has a unique solution $(w, s, B_2)$ such that $\sigma \equiv s- s_a$ satisfies $\sigma (0) = 0$ and such that $(q, G_2) \equiv (w-w_a, B_2 - B_{2a}) \in X((0, \tau ])$. In particular the following estimates hold for all $t \in (0, \tau]$~:
\beq
\label{7.1e}
\parallel \nabla^k \partial_t^j <x>^{\ell} q \parallel_2\ \leq C\ t^{2-j-k/2} (1 - \ell n \ t)^4
\eeq 

\noi for $0 \leq j, \ell \leq 1$ and $0 \leq 2j + \ell + k \leq 3$,
\beq
\label{7.2e}
\parallel \nabla^{k+1}   \partial_t^j G_2 \parallel_2\ \vee\ t \parallel \nabla^{k+1}   \partial_t^j  {\check G}_2 \parallel_2\ \leq C\ t^{5/2-j-k} (1 - \ell n \ t)^4
\eeq 

\noi for $0 \leq j, k, j+k \leq 1$.\par

In addition, the following estimates hold for all $t \in (0, \tau ]$~:
\beq
\label{7.3e}
\parallel \partial_t^j G_2 \parallel_2\ \vee\ t \parallel\partial_t^j {\check G}_2 \parallel_2\ \leq C\ t^{5/2-j} (1 - \ell n \ t)^4\ ,
\eeq 
\beq
\label{7.4e}
\parallel \nabla^{k+1}   \partial_t^j  G_1 \parallel_2\ \vee\ t \parallel \nabla^{k+1}   \partial_t^j {\check G}_1 \parallel_2\ \leq C\ t^{2-j-k/2} (1 - \ell n \ t)^4
\eeq 

\noi for $0 \leq j, k, j+k \leq 1$,
\beq
\label{7.5e}
\parallel \partial_t^j G_1 \parallel_2\ \vee\ t \parallel\partial_t^j {\check G}_1 \parallel_2\ \leq C\ t^{2-j} (1 - \ell n \ t)^4 \ ,
\eeq 
\beq
\label{7.6e}
\parallel \nabla^k  \partial_t^j  \sigma \parallel_2\ \leq C\ t^{2-j-k/2} (1 - \ell n \ t)^4
\eeq 

\noi for $j = 0, 1$ and $0 \leq k \leq 2$.}\\

The solution is actually unique under the conditions on $(q, G_2)$ stated in Proposition 5.2. \\

\noi {\bf Proof.} The result follows from Propositions 5.2 and 6.3 except for the low order estimates (\ref{7.3e}) (\ref{7.5e}). In particular (\ref{7.1e}) is a rewriting of (\ref{5.71e}) (\ref{5.73e}), while (\ref{7.2e}) is a rewriting of (\ref{5.74e}). The estimates (\ref{7.4e}) and (\ref{7.6e}) are a rewriting of (\ref{5.77e})-(\ref{5.81e}). The lower norm estimates (\ref{7.3e}) (\ref{7.5e}) are derived by the same method as in the proof of Lemmas 5.1 and 5.2.\par \nobreak \hfill $\sq$ \par

\noi {\bf Remark 7.1.} The estimates (\ref{7.4e}) for $G_1$ and ${\check G}_1$ actually hold for $j = 0, 1$, and $0 \leq k + 2j \leq 3$ for $G_1$, $0 \leq k + 2j \leq 2$ for ${\check G}_1$. The proof is an extension of that of Lemma 5.1. Unfortunately, a similar extension for $G_2$ and ${\check G}_2$ would require a reinforcement of the assumptions on $B_0$ and $w_+$, because of the explicit dependence of ${\cal B}_2$ and of ${\cal R}_4$ on $B_0$.\\

We now turn to the original system (\ref{2.6e})  (\ref{2.7e}) for $(u, A)$. The first task is to reconstruct the phase $\varphi$. Corresponding to $s_a = s_{a0} + s_{a1}$ defined by (\ref{6.1e}) (\ref{6.2e}) (\ref{6.3e}), we define
\bea
\label{7.7e}
&&\varphi_a = \int_1^t dt' \left ( t{'}^{-1} g (w_{a0} (t')) + {\check B}_{1L} (w_{a0}(t'))\right )\nn \\
&&+ 2 \int_0^t dt' \left (  t{'}^{-1} g (w_{a0} (t'), w_{a1} (t')) + {\check B}_{1L} (w_{a0}(t'), w_{a1} (t'))\right )
\eea

\noi so that $s_a = \nabla \varphi_a$. We shall also need a special term of $\varphi_a$, namely
\beq
\label{7.8e}
\varphi_b = \int_1^t dt' {\check B}_{1}(w_+) = (\ell n\ t) x \cdot B_1(w_+)\ .
\eeq 

\noi The phases $\varphi_a$ and $\varphi_b$ satisfy the following properties. We use again the notation (\ref{3.1e}) as in Lemma 6.1.\\

\noi {\bf Lemma 7.1.} {\it Let $w_+ \in H^5$, $xw_+ \in H^4$. Then} 
\beq
\label{7.9e}
\partial_t \ \varphi_b \in \left ( <x> \dot{H}^1 \cap \dot{H}^2 \cap \dot{H}^5 , t^{-1} \right )\ ,
\eeq 
\beq
\label{7.10e}
<x>^{-1} \partial_t \ \varphi_b \in \left ( \ddot{H}^6, t^{-1}\right )\ ,
\eeq 
\beq
\label{7.11e}
 \varphi_b \in \left ( <x> \dot{H}^1 \cap \dot{H}^2 \cap \dot{H}^5 , 1 - \ell n\ t  \right )\ ,
\eeq
\beq
\label{7.12e}
<x>^{-1} \varphi_b \in \left ( \ddot{H}^6, 1 - \ell n\ t \right )\ ,
\eeq 
\beq
\label{7.13e}
 \partial_t (\varphi_a - \varphi_b) \in \left ( \ddot{H}^4, t^{-1}\right )\ ,
\eeq 
\beq
\label{7.14e}
\varphi_a - \varphi_b \in \left ( \ddot{H}^4, 1- \ell n\ t\right )\ .
\eeq 

\noi {\bf Proof.} The properties of $\varphi_b$ follow from estimates of ${\check B}_1(w_+)$ and $B_1(w_+)$ which are identical with those of ${\check B}_{1a0}$ and $B_{1a0}$ in Lemma 6.1.\par

We rewrite $\varphi_a - \varphi_b$ as 
 $$\varphi_a - \varphi_b = \int_1^t dt' \left ( t{'}^{-1} g(w_{a0}) - {\check B}_{1S} (w_{a0}) + {\check B}_1(w_{a0}) - {\check B}_1(w_+)\right )$$
 $$+ 2 \int_0^{t} dt' \left ( t{'}^{-1} g(w_{a0}, w_{a1}) + {\check B}_{1L} (w_{a0}, w_{a1}) \right ) \ .$$
 
 \noi The contributions of $g$ and of ${\check B}_{1a1}$ are estimated in the same way as in Lemma 6.1. We next estimate
$$\parallel \nabla^k {\check B}_{1S} (w_{a0})\parallel_2 \ \leq \ t^{\beta (5-k)} \parallel \nabla^5 {\check B}_{1} (w_{a0})\parallel_2\ \leq C\ t^{-1 + \beta (5-k)}$$

\noi with $\beta = 1/2$ for $0 \leq k \leq 5$, so that
$${\check B}_{1S} (w_{a0}) \in \left ( H^k, t^{(3-k)/2} \right ) \qquad \hbox{for $0 \leq k \leq 5$}$$

\noi and the contribution of that term satisfies the required properties and estimates. Finally  
 $${\check B}_{1} (w_{a0}) - {\check B}_{1} (w_{+}) = 2 {\check B}_{1} \left ( (U(t) - \1) w_+, w_{a0} + w_+ \right )\ .$$
 
 \noi Now
 $$<x> \left ( U(t) - \1\right ) w_+ \in \left ( H^2, t\right )$$
 
 \noi so that by the same estimates as those of ${\check B}_{1a0}$ and ${\check B}_{1a1}$ in Lemma 6.1 
 $${\check B}_{1} (w_{a0})  - {\check B}_{1}(w_+) \in \left ( \ddot{H}^3, 1 \right ) \cap \left ( \dot{H}^5, t^{-1}\right )$$
 
 \noi and the contribution of that difference also satisfies the required properties and estimates. \par \nobreak \hfill $\sq$ \par
 
 Let now $(w, s, B_2)$ be a solution of the auxiliary system (\ref{2.41e}) as obtained in Proposition 7.1, and let $(q, \sigma) = (w-w_a, s-s_a)$ and $B = B_0 + B_1 + B_2$ with $B_1$ defined by (\ref{2.31e}) (\ref{2.32e}). We define
 \beq
\label{7.15e}
\psi = \int_0^t dt' \left \{ t{'}^{-1} \left ( g (q, 2w_a + q) + g(w_{a1})\right ) + {\check B}_{1L}(q, 2w_a+q) + {\check B}_{1L}(w_{a1})  \right \} (t')
\eeq 

\noi so that by  (\ref{2.52e})  (\ref{2.53e}) and  (\ref{6.8e})  (\ref{6.9e}), $\nabla \psi = \sigma$. From Proposition 7.1, especially  (\ref{7.6e}) and an $L^4$ estimate similar to those in Lemma 6.1, it follows that $\psi \in {\cal C}((0, \tau ], \ddot{H}^3)$ and that $\psi$ satisfies the estimates 
\beq
\label{7.16e}
\parallel \nabla^{k+1} \partial_t^j \psi \parallel_2 \ \leq C\ t^{2-j-k/2} (1 - \ell n\ t)^4
\eeq 

\noi for $j = 0,1$ and $0 \leq k \leq 2$. Finally we define $\varphi = \varphi_a + \psi$ with $\varphi_a$ defined by (\ref{7.7e}) so that $\nabla \varphi = s$ and $\varphi$ satisfies (\ref{2.38e}).\par

We can now define the modified wave operator for the MS system in the form (\ref{2.6e}) (\ref{2.7e}). We start from the asymptotic data $(u_+, A_+, \dot{A}_+)$ for $(u, A)$. We define $w_+ = Fu_+$, we define $B_0$ by (\ref{2.9e}) (\ref{2.17e}), namely 
\beq
\label{7.17e}
A_0 (t) = (\cos \omega t ) A_+ + \omega^{-1}( \sin \omega t )\dot{A}_+ = - t^{-1} D_0 (t) B_0 (1/t) \ .
\eeq 

\noi We define $(w_a, s_a, B_a)$ by (\ref{6.1e}) (\ref{6.2e}) (\ref{6.3e}). We solve the auxiliary system (\ref{2.41e}) by Proposition 7.1. We reconstruct the phase $\varphi = \varphi_a + \psi$ as explained above. We finally substitute $(w, \varphi , B_2)$ into (\ref{2.16e}) (\ref{2.33e}) (\ref{2.17e}), thereby obtaining a solution $(u, A)$ of the system (\ref{2.6e}) (\ref{2.7e}) defined for large time. The modified wave operator is the map $\Omega : (u_+, A_+, \dot{A}_+) \to (u,A)$ thereby obtained. \par

We now turn to the study of the asymptotic properties of $(u, A)$ and in particular of its convergence to its asymptotic form $(u_a, A_a)$ defined in analogy with (\ref{2.16e}) (\ref{2.17e}) by
\beq
\label{7.18e}
u_a(t) = M(t) D(t) \exp \left ( i \varphi_a (1/t)\right ) \overline{w}_a (1/t)\ ,
\eeq 
\beq
\label{7.19e}
A_a (t) = - t^{-1} D_0(t) \ B_a(1/t) = A_0(t) - t^{-1} D_0(t) \left ( B_{1a} + B_{2a} \right ) (1/t) \ .
\eeq 

\noi The properties of $u$ are best expressed in terms of $\widetilde{u}$ and $\widetilde{u}_a$ defined by 
\beq
\label{7.20e}
\widetilde{u}(t) = U(-t) u(t) \quad , \qquad \widetilde{u}_a(t) = U(-t) u_a(t)
\eeq

\noi so that
\beq
\label{7.21e}
\widetilde{u}(t) = M(t)^* F^* \exp (i \varphi (1/t)) \overline{w}(1/t)\ ,
\eeq 
\beq
\label{7.22e}
\widetilde{u}_a(t) = M(t)^* F^* \exp (i \varphi_a (1/t)) \overline{w}_a(1/t)\ .
\eeq 

In order to translate the properties of $(w, \varphi )$ into properties of $u$, we need the following commutation relations
\beq
\label{7.23e}
M(t)^* F^*\nabla = - ix M(t)^* F^* \ ,
\eeq
\beq
\label{7.24e}
M(t)^* F^* x = - i M(t)^* \nabla F^* = (- i \nabla + x/t) M(t)^* F^* \ ,
\eeq 
\beq
\label{7.25e}
M(t)^* F^* i \partial_t = \left ( i \partial_t +\left  ( 2t^2\right )^{-1} x^2 \right ) M(t)^* F^*
\eeq

\noi so that
\beq
\label{7.26e}
M(t)^* F^* \overline{i\partial_t v} (1/t) = \left ( i t^2 \partial_t + (1/2) x^2 \right ) M(t)^* F^* \overline{v} (1/t)\ .
\eeq

We shall need in addition the following lemmas, which we state in terms of a general function $h$ as considered in Section 3, although we shall use them only in the case $h(t) = t^2 (1 - \ell n\ t)^4$.\\

\noi {\bf Lemma 7.2.} {\it Let $0 < \tau \leq 1$, let $(q, 0)  \in X((0, \tau ])$ so that
\beq
\label{7.27e}
\parallel \nabla^{k} \partial_t^j x^{\ell} q\parallel_2 \ \leq C\ t^{-j-k/2} \ h(t)
\eeq 

\noi for $0 \leq j, \ell \leq 1$ with $0 \leq 2j+k+\ell \leq 3$ and for all $t\in (0, \tau ]$. Let $\theta \in$\break\noindent ${\cal C}^1((0, \tau ],  <x>L^{\infty} \cap \dot{H}^2 \cap \dot{H}^3)$ and let $\theta$ satisfy 
\beq
\label{7.28e}
\parallel \nabla^{k+1} \partial_t^j \theta \parallel_2 \ \leq C\ t^{1/4-j-k/2}
\eeq 

\noi for $(j, k) = (0,1), (0,2)$ and $(1,1)$, and
\beq
\label{7.29e}
\parallel <x>^{-1} \partial_t \theta \parallel_{\infty} \ \leq C\ t^{-1}
\eeq 

\noi for all $t \in (0, \tau ]$. Then $q \exp (-i \theta )$ satisfies the same estimates as (\ref{7.27e}) for $q$ with the exception of the case $(j, k, \ell) = (1, 0,1)$. If in addition 
\beq
\label{7.30e}
\parallel \partial_t \theta \parallel_{\infty} \ \leq C\ t^{-1}
\eeq

\noi then also the latter estimate holds.}\\

\noi {\bf Proof.} We consider the various cases successively. The result follows from the assumptions, from Sobolev inequalities and from the estimates below. \par

\noi $j = \ell = 0$.
$$\parallel  q \exp (- i \theta ) \parallel _2 \ = \ \parallel  q\parallel _2\ \leq C\ h \ ,$$
$$\parallel \nabla q \exp (-i \theta) \parallel _2\ \leq\ \parallel  \nabla q \parallel _2\ + \ \parallel  q\parallel _3 \ \parallel \nabla \theta \parallel _6\ \leq C \ t^{-1/2}\ h\ ,$$
\begin{eqnarray*}
\parallel \nabla^2 q \exp (- i \theta ) \parallel_2 &\leq& \parallel  \nabla^2 q \parallel _2\ + \ \parallel  \nabla q\parallel _3 \ \parallel \nabla \theta \parallel _6\ +\  \parallel  q\parallel _6 \ \parallel \nabla \theta \parallel _6^2\\
&&+\  \parallel  q\parallel _{\infty} \ \parallel \nabla^2 \theta \parallel _2\ \leq C \ t^{-1}\ h\ ,
\end{eqnarray*}
\begin{eqnarray*}
&&\parallel \nabla^3 q \exp (- i \theta ) \parallel_2 \ \leq \ \parallel  \nabla^3 q \parallel _2\ + \ \parallel  \nabla^2 q\parallel _3 \ \parallel \nabla \theta \parallel _6\ +\  \parallel  \nabla q\parallel _6 \ \parallel \nabla \theta \parallel _6^2\\
&&+\  \parallel  \nabla q\parallel _{\infty} \ \parallel \nabla^2 \theta \parallel _2\ +\ \parallel  q\parallel _{\infty} \left (  \parallel \nabla \theta \parallel _6^3 \ + \ \parallel \nabla^3 \theta \parallel _2\right )  \\
&&+\  \parallel  q\parallel _{6} \ \parallel \nabla \theta \parallel _6\  \parallel \nabla^2 \theta \parallel _6 \ \leq C \ t^{-3/2}\ h\ .
\end{eqnarray*}

\noi $j = 0$, $\ell = 1$. It suffices to replace $q$ by $xq$ in the first three estimates above.\par

\noi $j = 1$, $\ell = 0$.
$$\parallel  \partial_t (q \exp (- i \theta ) ) \parallel _2 \ \leq \ \parallel  \partial_t q  \parallel _2 \ +\ \parallel  <x> q  \parallel _2\ \parallel  <x>^{-1} \partial_t \theta  \parallel _{\infty} \ \leq C \ t^{-1}\ h\ ,$$
\begin{eqnarray*}
&&\parallel  \nabla \partial_t (q \exp (- i \theta ) ) \parallel _2 \ \leq \ \parallel  \nabla \partial_t q  \parallel _2 \ +\ \parallel  \partial_t q  \parallel _3\ \parallel  \nabla \theta  \parallel _{6}\\
&&+\ \parallel  <x> \nabla q  \parallel _2\ \parallel  <x>^{-1} \partial_t \theta  \parallel _{\infty}\ +\ \parallel  q \parallel _3\ \parallel  \nabla \partial_t \theta \parallel _6\\
&&\parallel  <x>q \parallel _3\ \parallel  \nabla  \theta  \parallel _6\ \parallel  <x>^{-1} \partial_t \theta  \parallel _{\infty}    \ \leq C \ t^{-3/2}\ h\ .
\end{eqnarray*}

\noi $j = \ell = 1$.
$$\parallel  \partial_t (xq \exp (- i \theta ) ) \parallel _2 \ \leq \ \parallel  \partial_t xq  \parallel _2 \ +\ \parallel  xq  \parallel _2\ \parallel \partial_t \theta  \parallel _{\infty} \ \leq C \ t^{-1}\ h\ .$$
\nobreak \hfill $\sq$ \par

\noi {\bf Lemma 7.3.} {\it Let $0 < \tau \leq 1$, let $v \in {\cal C}((0, \tau ], H^3) \cap {\cal C}^1((0, \tau ], H_3^1 \cap L^{\infty })$, $<x> v \in {\cal C}((0, \tau ], H^2) \cap {\cal C}^1((0, \tau ], L^3)$ and let $v$ satisfy
 \beq
\label{7.31e}
\left \{ \begin{array}{l} \parallel v; H^3 \parallel\ \vee \ \parallel <x> v; H^2 \parallel\ \leq C\cr \cr\parallel \partial_t v; H_3^1 \cap  L^{\infty }\parallel\ \vee \ \parallel <x> \partial_t v \parallel_3 \ \leq C\ t^{-1} \end{array} \right .
\eeq

\noi for all $t \in (0, \tau]$. Let $\psi \in {\cal C}((0, \tau ], \ddot{H}^3)$ satisfy 
\beq
\label{7.32e}
\parallel \nabla^{k+1} \partial_t^j  \psi \parallel_2 \ \leq C\ t^{-j-k/2}\ h(t)
\eeq 

\noi for $j = 0, 1$ and $0 \leq k \leq 2$, for all $t\in (0, \tau ]$. Then $q_1 \equiv v (\exp (- i \psi ) - 1)$ satisfies the estimates (\ref{7.27e}) for the same values of $j$, $k$, $\ell$. }\\

\noi {\bf Proof.} We consider the various cases successively. We estimate\par

\noi $j = \ell = 0$.
$$\parallel q_1 \parallel_2 \ \leq \ \parallel v \parallel_3\ \parallel \psi \parallel_6\ \leq C\ h\ ,$$
$$\parallel \nabla q_1 \parallel_2 \ \leq \ \parallel \nabla v \parallel_3\ \parallel \psi \parallel_6\  + \ \parallel v \parallel_{\infty}\ \parallel \nabla \psi \parallel_2\ \leq C\ h\ ,$$
\begin{eqnarray*}
\parallel \nabla^2 q_1 \parallel_2 &\leq& \parallel \nabla^2 v \parallel_2\ \parallel \psi \parallel_{\infty}\  + \ \parallel \nabla v \parallel_{6}\ \parallel \nabla \psi \parallel_3\\
&&+ \ \parallel v \parallel_{\infty}\left (  \parallel \nabla^2 \psi \parallel_2\ + \ \parallel \nabla \psi \parallel_4^2 \right )\\
&\leq& C\ h\left ( t^{-1/4} + t^{-1/2} \right ) +  h^2 t^{-3/4} \leq C\ h\ t^{-1/2}\ ,
\end{eqnarray*}
\begin{eqnarray*}
&&\parallel \nabla^3 q_1 \parallel_2 \ \leq \ \parallel \nabla^3 v \parallel_2\ \parallel \psi \parallel_{\infty}\  + \ \parallel \nabla^2 v \parallel_{6}\ \parallel \nabla \psi \parallel_3\\
&&+ \ \parallel \nabla v \parallel_{\infty}\left (  \parallel \nabla^2 \psi \parallel_2\ + \ \parallel \nabla \psi \parallel_4^2 \right )\\
&&+\  \parallel v \parallel_{\infty} \left (  \parallel \nabla^3 \psi \parallel_2\ + \ \parallel \nabla^2 \psi \parallel_3\ \parallel \nabla \psi \parallel_6\ + \ \parallel \nabla \psi \parallel_6^3 \right ) \\
&&\leq C\ h\left ( t^{-1/4} + t^{-1/2}+ t^{-1} \right ) +  h^2\left (  t^{-3/4} + t^{-5/4} \right ) + h^3\ t^{-3/2}  \leq C\ h\ t^{-1}\ .
\end{eqnarray*}

\noi $j= 0$ $\ell = 1$. It suffices to replace $v$ by $xv$ in the first three estimates above.\par

\noi $j = 1$, $\ell = 0$.
$$\parallel \partial_t q_1 \parallel_2 \ \leq \ \parallel \partial_t  v \parallel_3\ \parallel \psi \parallel_6\  + \ \parallel v \parallel_{3}\ \parallel \partial_t  \psi \parallel_6\ \leq C\ h \ t^{-1}\ ,$$
\begin{eqnarray*}
\parallel \nabla  \partial_t  q_1 \parallel_2 &\leq& \parallel \nabla \partial_t  v \parallel_3\ \parallel \psi \parallel_{6}\  + \ \parallel \nabla v \parallel_{3}\ \parallel  \partial_t  \psi \parallel_6\ +\ \parallel  \partial_t  v \parallel_{\infty} \ \parallel  \nabla \psi \parallel_2\\
&&+ \ \parallel v \parallel_{\infty}\left (  \parallel \nabla \partial_t   \psi \parallel_2\ + \ \parallel \partial_t   \psi \parallel_6\  \parallel \nabla \psi \parallel_3\right )\\
&\leq& C\ h\  t^{-1} + h^2\ t^{-5/4}  \leq C\ h\ t^{-1}\ .
\end{eqnarray*}

\noi $j = 1$, $\ell = 1$. It suffices to replace $v$ by $xv$ in the first of the previous estimates.\par \nobreak \hfill $\sq$ \par

In order to translate the properties of $B$ into properties of $A$, we need the following commutation relation 
\beq
\label{7.33e}
\nabla^k S^j A(t) = (-)^{j+1}\ t^{-1-k} D_0 (t) \left ( \nabla^k (t \partial_t)^j B\right ) (1/t)
\eeq 

\noi where $S = t \partial_t + x \cdot \nabla + 1$ (see Section 3, especially (\ref{3.28e})).\\

\noi {\bf Proposition 7.2.} {\it Let $\beta = 1/2$. Let $u_+$ be such that $w_+ = Fu_+ \in H^5$, $xw_+ \in H^4$ and such that $w_+$ satisfies the support condition (\ref{6.90e}). Let $A_+, \dot{A}_+$ satisfy (\ref{3.29e})  (\ref{3.30e}) and  (\ref{6.91e}) for $0 \leq j, k \leq 1$. Define $(w_a, s_a, B_a)$ by  (\ref{6.1e})  (\ref{6.2e})  (\ref{6.3e})  (\ref{2.43e}) and $(\varphi_a, u_a, A_a)$ by (\ref{7.7e}) (\ref{7.18e}) (\ref{7.19e}). Let $(w, s, B_2)$ be the solution of the auxiliary system (\ref{2.41e}) obtained in Proposition 7.1, let $\varphi = \varphi_a + \psi$ with $\psi$ defined by (\ref{7.15e}), let $B = B_0 + B_1 + B_2$, let $(u, A)$ be defined by (\ref{2.16e})  (\ref{2.17e}) and let $\widetilde{u}$ be defined by (\ref{7.20e}). Let $T = \tau^{-1}$ and $I = [T, \infty )$. Then \par

(1)\qquad \qquad \qquad  $x^k  \partial_t^j  \nabla^{\ell}  \widetilde{u} \in {\cal C}(I,L^2)$

\noi for $0 \leq j, \ell , j+ \ell \leq 1$ and $0 \leq 2j+k+\ell \leq 3$, and $\widetilde{u}$ satisfies the following estimates for the same values of $j$, $k$, $\ell$ and for all $t\in I$~:
\beq
\label{7.34e}
\parallel x^k \partial_t^j \nabla^{\ell}( \widetilde{u} - \widetilde{u}_a)\parallel_2\ \leq C\ t^{-2-j+k/2} (1 + \ell n\ t)^4 \ .
\eeq

\noi Furthermore $\partial_t \nabla U(-t) \exp (- i \varphi_b (1/t, x/t))u(t) \in {\cal C}(I,L^2)$ and the following estimate holds for all $t \in I$~: 
\beq
\label{7.35e}
\parallel \partial_t \nabla U(-t) \exp (- i \varphi_b (1/t, x/t))(u(t) - u_a(t))\parallel_2\ \leq C\ t^{-3} (1 + \ell n\ t)^4\ .
\eeq 

\noi Finally the following estimate holds
\beq
\label{7.36e}
\parallel x^{\ell} (u-u_a)\parallel_r\ \leq C\ t^{-2 + \ell - \delta (r)/2} (1 + \ell n\ t)^4
\eeq 

\noi for $\ell = 0,1$, for $2 \leq r \leq \infty$ and for all $t\in I$, with $\delta (r) = 3/2 - 3/r$. \par

(2) $A \in {\cal C}(I, \dot{H}^1 \cap \dot{H}^2)$, $x \cdot A \in {\cal C}(I, \dot{H}^2)$, $SA \in {\cal C}(I, {H}^1)$, $Sx \cdot A \in {\cal C}(I, <x>\dot{H}^1)$, where $S = t \partial_t + x \cdot \nabla + 1$. Furthermore $A-A_a$, $x\cdot (A-A_a) \in {\cal C}(I, {H}^2)$, $S(A-A_a)$, $Sx\cdot (A - A_a) \in {\cal C}(I, {H}^1)$ and the following estimates hold for all $t\in I$~:
\beq
\label{7.37e}
\parallel S^j(A-A_a)\parallel_2\ \vee\ t^{-1}\parallel S^jx\cdot (A-A_a)\parallel_2 \ \leq C\ t^{-3/2}  (1 + \ell n\ t)^4\ ,
\eeq
\bea
\label{7.38e}
&&\parallel \nabla^{k+1} S^j(A-A_a) \parallel_2\ \vee\ t^{-1}\parallel \nabla^{k+1} S^jx\cdot (A-A_a)\parallel_2 \nn\\
&&\qquad \qquad \leq C\ t^{-5/2-k/2}  (1 + \ell n\ t)^4\ ,
\eea  

\noi for $0 \leq j, k, j+k \leq 1$.}\\

\noi {\bf Proof.}

\noi \underline{Part (1).} The existence of the solution $u$ with the regularity stated follows from the existence and regularity of $w$ obtained in Proposition 7.1, using Lemma 6.6 which ensures (\ref{6.88e}), from the reconstruction of the phase $\varphi = \varphi_a + \psi$ performed above and from the regularity properties of $\varphi_a$ and $\psi$, through the change of variables in the form (\ref{7.21e}). The details of the proof follow from the estimates to be given below. We now turn to the derivation of the estimates. We first consider
\bea
\label{7.39e}
v&=& w \exp (-i \varphi ) - w_a \exp (-i \varphi_a)\nn \\
&=& \left ( w (\exp (-i \psi ) - 1 ) + (w - w_a)\right ) \exp (- i \varphi_a)\ .
\eea 

\noi It follows from (\ref{7.16e}) that $\psi$ satisfies the assumption (\ref{7.32e}) of Lemma 7.3 for the relevant $h$, while $w$ obviously satisfies (\ref{7.31e}), so that by that Lemma,\break\noindent $q_1 = w (\exp (-i \psi ) - 1)$ satisfies the same estimates (\ref{7.27e}) as $q = w- w_a$ for the relevant $h$. We now show that the assumptions on $\theta$ made in Lemma 7.2 are satisfied by $\varphi_a$ and/or $\varphi_a - \varphi_b$. For that purpose we use Lemma 7.1. In fact from (\ref{7.11e}) (\ref{7.14e}) it follows that
$$\parallel \varphi_a ; \dot{H}^2 \cap \dot{H}^3\parallel\ \leq C(1 - \ell n\ t)$$

\noi so that $\varphi_a$ satisfies (\ref{7.28e}) for $j = 0$ and $k = 1,2$. From (\ref{7.9e}) (\ref{7.13e}), it follows that
$$\parallel \partial_t \varphi_a ; \dot{H}^2 \parallel\ \leq C\ t^{-1}$$

\noi so that $\varphi_a$ satisfies (\ref{7.28e}) with $j=k=1$. From (\ref{7.13e}) it follows that 
$$\parallel \partial_t (\varphi_a - \varphi_b) ; \dot{H}^1 \cap \dot{H}^2\parallel\ \leq C\ t^{-1}$$

\noi so that $\varphi_a - \varphi_b$ satisfies (\ref{7.30e}). Finally from (\ref{7.10e}) it follows that 
$$\parallel <x>^{-1} \partial_t \varphi_b  ; \dot{H}^1 \cap \dot{H}^2\parallel\ \leq C\ t^{-1}$$

\noi so that $\varphi_b$ and therefore also $\varphi_a$ satisfy (\ref{7.29e}). It follows now from the previous results and from Lemma 7.2 that $v$ defined by (\ref{7.39e}) satisfies the estimates
\beq
\label{7.40e}
\parallel \nabla^k \partial_t^j x^{\ell}v\parallel_2\ \leq C\ t^{-j-k/2}\ h(t) 
\eeq

\noi for $0 \leq j, \ell , j+ \ell \leq 1$ and $0 \leq 2j+k+ \ell \leq 3$, and 
\beq
\label{7.41e}
\parallel  \partial_t  x  v \exp (i \varphi_b) \parallel_2\ \leq C\ t^{-1}\ h(t) \ .
\eeq

\noi We now derive the estimates  (\ref{7.34e})  (\ref{7.35e}). From  (\ref{7.21e})  (\ref{7.22e}) it follows that 
\beq
\label{7.42e}
\left ( \widetilde{u} -  \widetilde{u}_a\right ) (t) = M(t)^* F^* \overline{v} (1/t)\ .
\eeq

\noi From (\ref{7.40e}) (\ref{7.42e}) it follows immediately that 
\beq
\label{7.43e}
\parallel x^k ( \widetilde{u} -  \widetilde{u}_a) (t)\parallel_2\ = \ \parallel \nabla^k v (1/t) \parallel_2\ \leq C\ t^{k/2} \ h(1/t)
\eeq

\noi for $0 \leq k \leq 3$. From (\ref{7.40e}) (\ref{7.42e}) and the commutation relation (\ref{7.26e}), it then follows that
$$\parallel x^k \left ( it^2 \partial_t + (1/2) x^2 \right ) ( \widetilde{u} -  \widetilde{u}_a) (t)\parallel_2\ = \ \parallel \nabla^k (\partial_t v) (1/t) \parallel_2\ \leq C\ t^{1+k/2} \ h(1/t)$$

\noi for $0 \leq k\leq 1$, which together with (\ref{7.43e}) implies
\beq
\label{7.44e}
\parallel x^k  \partial_t ( \widetilde{u} -  \widetilde{u}_a) (t)\parallel_2\ \leq C\ t^{-1+k/2} \ h(1/t)\ .
\eeq

\noi From (\ref{7.40e}) (\ref{7.42e}) and the commutation relation (\ref{7.24e}), it then follows that
$$\parallel x^k( i\nabla - x/t  ) ( \widetilde{u} -  \widetilde{u}_a) (t)\parallel_2\ = \ \parallel \nabla^k x v (1/t) \parallel_2\ \leq C\ t^{k/2} \ h(1/t)$$

\noi for $0 \leq k \leq 2$, which together with (\ref{7.43e}) implies
\beq
\label{7.45e}
\parallel x^k  \nabla ( \widetilde{u} -  \widetilde{u}_a) (t)\parallel_2\ \leq C\ t^{k/2} \ h(1/t)\ .
\eeq

\noi Collecting (\ref{7.43e})-(\ref{7.45e}) for the relevant $h$ yields (\ref{7.34e}). We now turn to (\ref{7.35e}). We define
\beq
\label{7.46e}
\widetilde{u}_-(t) = U(-t) \exp \left ( - i \varphi_b (1/t, x/t)\right ) (u(t) - u_a(t))\ .
\eeq

\noi From a minor variation of (\ref{7.21e}) (\ref{7.22e}), it follows that 
\beq
\label{7.47e}
\widetilde{u}_-(t) = M(t)^* F^* \ \overline{v \exp (i \varphi_b)} (1/t)\ .
\eeq

\noi From (\ref{7.41e}) (\ref{7.47e}) and the commutation relations (\ref{7.24e}) (\ref{7.26e}), it follows that
\begin{eqnarray*}
\parallel \left ( i t^2 \partial_t + (1/2) x^2 \right ) (i \nabla - x/t) \widetilde{u}_-(t) \parallel_2 &=&  \parallel x \left ( \partial_t v \exp (i \varphi_b)\right ) (1/t) \parallel_2\\
&\leq& C\ t\ h(1/t) \ .
\end{eqnarray*}

\noi From this estimate and from the analogues of  (\ref{7.43e})-(\ref{7.45e})  for $\widetilde{u}_-$, which can be derived in the same way, it follows that 
\beq
\label{7.48e}
\parallel \partial_t \nabla \widetilde{u}_-(t) \parallel_2\ \leq C\ t^{-1} \ h(1/t)
\eeq

\noi which reduces to (\ref{7.35e}) for the relevant $h$. Finally using the relation 
$$(u-u_a)(t) = M(t) \ D(t)\ \overline{v(1/t)}$$

\noi we estimate
\begin{eqnarray*}
&&\parallel x^{\ell} (u-u_a) (t) \parallel_r\ = \ t^{-\delta + \ell} \parallel x^{\ell} v (1/t) \parallel_r\\
&&\leq \ t^{-\delta + \ell} \parallel x^{\ell} v (1/t) \parallel_2^{1- \delta /2}\ \parallel \nabla^2 x^{\ell} v (1/t) \parallel_2^{\delta /2}\ \leq C \ t^{- \delta /2 + \ell } \ h(1/t)
\end{eqnarray*}

\noi for $\ell = 0, 1$ and $0 \leq \delta = \delta (r) = 3/2 - 3/r \leq 3/2$, by the use of (\ref{7.38e}). That estimate reduces to (\ref{7.36e}) for the relevant $h$.\\

\noi \underline{Part (2).} The existence and regularity properties of $A$ follow from Proposition 7.1 as regards $A - A_a$, while the regularity properties of $A_a$ follow from Lemma 3.7 and Proposition 6.1, through the change of variables (\ref{2.17e}). We now turn to the estimates (\ref{7.37e}) (\ref{7.38e}). From the commutation relation (\ref{7.33e}), it follows that 
\beq
\label{7.49e}
\nabla^k S^j (A-A_a) (t) = (-)^{j+1} \ t^{-1-k} \ D_0(t) \left ( \nabla^k (t \partial_t )^j G \right ) (1/t)\ ,
\eeq
\beq
\label{7.50e}
\nabla^k S^j x\cdot (A-A_a) (t) = (-)^{j+1} \ t^{-1-k} \ D_0(t) \left ( \nabla^k (t \partial_t )^j {\check G} \right ) (1/t)\ .
\eeq

\noi The estimates (\ref{7.37e}) (\ref{7.38e}) follow immediately from (\ref{7.2e})-(\ref{7.5e}) and (\ref{7.49e}) (\ref{7.50e}).\par \nobreak \hfill $\sq$\par 

\section*{Acknowledgments}
\hspace*{\parindent}
Part of this work was done while one of the authors (G. V.) was visiting the Department of Mathematics of the New York University (DMNYU), New York, USA, the Institut des Hautes Etudes Scientifiques (IHES), Bures-sur-Yvette, France and the Laboratoire de Physique Th\'eorique (LPT), Universit\'e de Paris-Sud, Orsay, France. He is very grateful to Professor Jalal Shatah, Chairman of the DMNYU, to Professor Jean-Pierre Bourguignon, Director of the IHES and to Professor Hendrik-Jan Hilhorst, Director of the LPT, for the kind hospitality extended to him.

\end{document}